\documentclass[10pt,a4paper]{article}
\usepackage[utf8]{inputenc}
\usepackage{amsmath}
\usepackage{amsfonts}
\usepackage{amssymb}
\usepackage{hyperref} 
\usepackage{graphicx}
\usepackage{color}
\usepackage[french]{babel}
\usepackage{url}
\usepackage{fancyhdr}
\usepackage{multicol}
\author{Abdelkader BENAISSAT et Mouad MOUTAOUKIL}
\title{\fbox{Corps valués locaux}}
\begin{document}
\vspace*{\stretch{0.2}}
\begin{center}
{\large Abdelkader BENAISSAT}
 \end{center}
\begin{center}
{\large  Mouad MOUTAOUKIL}
\end{center}
\vspace*{\stretch{0.25}}
\begin{center}
\includegraphics[scale=0.8]{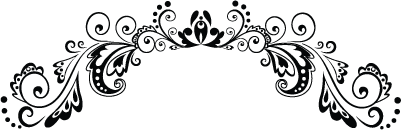}  
\end{center}
	\begin{center}
	{\LARGE {\Huge Corps valués locaux}}
		\subparagraph*{}
		{\LARGE {\huge Et applications à la théorie des nombres}}
	\end{center}
\begin{center}
\includegraphics[scale=1]{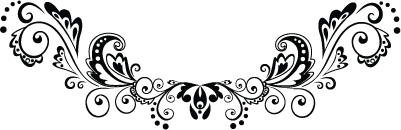}  
\end{center}
\paragraph*{}
\begin{center}
Faculté des sciences Dhar El Mehraz - Département de mathématiques
\end{center}
\begin{center}
Université Sidi Mohamed Ben Abdellah, Fès - Maroc 	
\end{center}

\vspace*{\stretch{1}}

\newpage
\vspace*{\stretch{1}}
\begin{center}
{\LARGE \textbf{\textcolor{blue}{Résumé}}}
\end{center}
\subparagraph*{} Beaucoup de sujets mathématiques de recherche très actifs de nos jours touchent les concepts de corps valués et corps valués locaux, surtout le corps des nombres $p$-adiques $\mathbb{Q}_{p}$ et le corps des séries formelles de Laurent $F((X))$. Les corps valués locaux sont une notion située un peu dans la frontière entre la théorie des nombres, l'algèbre et la topologie. Ils utilisent beaucoup de définitions et théorèmes -plus ou moins avancés- d'algèbre et topologie générales. Progressivement, on ira du général au local, des corps valués aux corps valués locaux, dont on verra quelques applications surtout en théorie des nombres élémentaire et algébrique.  

\paragraph*{Mots clés :} corps locaux, théorie des nombres, valuation, nombres p-adiques, vulgarisation.
\subparagraph*{}
\begin{center}
{\LARGE \textbf{\textcolor{blue}{Abstract}}}
\end{center}
\subparagraph*{}
Many active mathematical research topics nowadays include the concepts of valued fields and local fields, especially the local field of $p$-adic numbers $\mathbb{Q}_{p}$ and the field of formal Laurent series $F((X))$. Local fields are a notion situated in the boundary between number theory, algebra and topology. They use many definitions and theorems - more or less advanced - of general algebra and topology. Gradually, we will go from the general to the local, from the valued fields to the local fields, of which we will discuss some applications, especially in elementary and algebraic number theory.


\vspace*{\stretch{1}}
\newpage

\tableofcontents
\newpage

\vspace*{\stretch{1}}
\begin{center}
{\LARGE \textbf{\textcolor{blue}{Introduction générale}}}
\end{center}

\subparagraph*{} La notion de valeur absolue, qu'on connaît depuis le collège et qu'on généralise au module en terminale, peut être encore plus généralisée, à tout corps, introduisant ainsi les corps valués, qu'on passera au filtre de la localisation pour récolter un fruit très utile et intéressant : les corps valués locaux. Ce travail se présente dans le cadre de notre projet de fin d 'études en vue de l'obtention d'une licence en mathématiques, qui traite justement de "corps valués locaux".
\addcontentsline{toc}{part}{Introduction générale} 
\subparagraph*{} Après un travail acharné et passionné, incluant des nuits blanches et beaucoup de documentation, nous décidâmes de commencer notre projet par une partie préliminaire introduisant les notions et outils d'algèbre générale et de topologie nécessaires à la bonne compréhension des deux parties suivantes et qui constituent le cœur de notre sujet. 
\subparagraph*{} La première partie présentera les corps valués et leurs propriétés, commençant par l'application valeur absolue et terminant par quelques propriétés topologiques, en passant par les notions de complétude et surtout complétion.

\subparagraph*{} La deuxième partie est le cœur et corps de notre sujet. Ayant désormais les moyens nécessaires à définir une telle notion qu'est les corps valués locaux, on clôturera en présentant les extensions de ces corps ainsi que quelques unes de leurs structures. 
\subparagraph*{} La troisième et dernière partie consistera à concrétiser un peu les notions discutées en introduisant quelques applications et se terminera en revenant un peu au plus simple en vulgarisant le corps local des nombres p-adiques d'une façon accessible même aux lycéens, ou encore aux collégiens précoces.

\vspace*{\stretch{1}}
\newpage

\vspace*{\stretch{1}}
\begin{center}
{\LARGE \textbf{\textcolor{blue}{Positionnement historique}}}
\end{center}

\subparagraph*{} Les notions relatives aux corps valués virent le jour, pour la plupart,
\addcontentsline{toc}{part}{Positionnement historique}
 au cours du $20^{eme}$ siècle, grâce à une génération de mathématiciens influencée par les idées de \textsc{Fréchet} et de \textsc{Riesz} sur la topologie, et par celles de \textsc{Steinitz} sur l'algèbre. Cette génération va savoir rendre assimilables et mettre à leur vraie place les travaux de \textsc{Hensel}. Dès 1913, \textsc{Kürschak} définit de façon générale la notion de valeur absolue et reconnaît l'importance des valeurs absolues ultramétriques (donnant l'exemple par la valeur p-adique). \textsc{Ostrowski} va ensuite déterminer toutes les valeurs absolues sur le corps des rationnels $\mathbb{Q}$.
\subparagraph*{} Entre 1920 et 1935, cette théorie va encore avancer avec une étude plus détaillée des valeurs absolues non nécessairement discrètes et l'introduction d'une notion plus générale qui est la valuation par \textsc{Krull}. Des études plus profondes s'ensuivirent, traitant surtout de la structure des corps valués complets et des anneaux locaux complets.
\vspace*{\stretch{1}}
\newpage

\vspace*{\stretch{1}}
\begin{center}
\includegraphics[scale=0.8]{2.png}  
\end{center}
	\begin{center}
	{\LARGE {\Huge Préliminaires}}
		
	\end{center}
\begin{center}
\includegraphics[scale=1]{3.png}  
\end{center}
\addcontentsline{toc}{part}{Préliminaires}
\vspace*{\stretch{1}}
\newpage
\subparagraph*{}
\paragraph*{Introduction} Avant de passer au vif de notre sujet, nous avons jugé nécessaire de présenter cette partie préliminaire, qui constitue un rappel et un complément sur des notions d'algèbre et topologie générales dont on aura besoin pour une présentation correcte des corps valués locaux. Cependant, ce rappel/complément est loin de cerner tous les concepts qui seront abordés, le lecteur est donc supposé familier avec les définitions des notions basiques d'algèbre et de topologie, et les propriétés élémentaires de ces structures, qui ont été vues et revues pendant notre programme de licence, et introduites bien avant. 
 
\section*{1. Algèbre générale}
\addcontentsline{toc}{section}{1. Algèbre générale}
\subsection*{1.1. Compléments sur les groupes }
\addcontentsline{toc}{subsection}{1.1. Compléments sur les groupes }
La théorie des groupes constitue une partie fondamentale de l'algèbre générale. La notion de groupe quotient, dont on aura besoin pour notre sujet, est apparue, pour la première fois, chez \textsc{Jordan}, mais c'est \textsc{Hölder} qui a introduit l'expression "quotient des groupes $ G $ et $ H $" en 1889.

\paragraph*{1.1.1. Groupes quotients}
\subparagraph*{Construction du quotient d'un groupe :}
Soit $(G,.)$ un groupe et $(H,.)$ un sous-groupe de $G$. On considère la relation $\Re$ définie sur $G$ par : \begin{center}
$x \Re y \Leftrightarrow x^{-1}.y \in H$
\end{center}
On vérifie facilement que $\Re$ est une relation d'équivalence sur $G$. On note $G/H$ l'ensemble des classes d'équivalence de la relation $\Re$ sur $G$. C'est l'ensemble quotient du groupe $G$ suivant $H$. Sous des conditions sur $H$, on peut définir sur l'ensemble $G/H$ une loi de groupe. 
\subparagraph*{Sous-groupe distingué :}
On dit que le sous-groupe $H$ de $G$ est \textbf{distingué} ou normal, si pour tout $g$ dans $G$ et tout $h$ dans $H$ on a : $ g.h.g^{-1} \in H$.
\\On note $H \triangleleft G$ le fait que $H$ soit normal dans $G$.
\subparagraph*{Construction du groupe quotient :} 
Si $H$ est un sous-groupe distingué de $G$, l'ensemble $G/H$ muni de la loi interne induite de $G$ sur $G/H$ a une structure de groupe.   
\subparagraph*{Exemple :} 
\subparagraph*{$\ast$} Soit $n \in \mathbb{N} $. Si on définit sur $\mathbb{Z}$ la relation d'équivalence $\sim$ suivante : $x \sim y \Longleftrightarrow n/x-y $.\\
On considère le groupe additif $ (\mathbb{Z},+) $ et son sous-groupe $n\mathbb{Z}$. L'ensemble des classes d'équivalence de la relation $\sim$, $\mathbb{Z}/n\mathbb{Z}$, muni de la loi induite $+$ est un groupe commutatif. 
\subparagraph*{$\ast$} On peut définir sur $\mathbb{Z}/n\mathbb{Z}$ une loi multiplicative. L'ensemble $\mathbb{Z}/n\mathbb{Z} - \lbrace0\rbrace$ muni de cette loi n'est pas toujours un groupe. Une conséquence directe du théorème de Bézout est que $(\mathbb{Z}/n\mathbb{Z} - \lbrace0\rbrace, \times)$ est un groupe si et seulement si $n$ est un nombre premier. 
\paragraph*{1.1.2. Propriétés}
\subparagraph*{$\ast$} $G/G$ est un groupe trivial, réduit à l'élément neutre de $G$.
\subparagraph*{$\ast$} Si $G$ est commutatif ou cyclique, il en est de même pour $G/H$. \\
Rappelons qu'il n'existe, à isomorphisme près, qu'un seul groupe cyclique infini qui est $ (\mathbb{Z}, +)$, et qu'un seul groupe d'ordre $n$, qui est $\mathbb{Z}/n\mathbb{Z}$.

\subsection*{1.2. Compléments sur les anneaux et les corps }
\addcontentsline{toc}{subsection}{1.2. Compléments sur les anneaux et les corps }
La théorie des anneaux et des corps a été développée à partir de la fin du 19ème siècle, notamment sous l'influence de David \textsc{Hilbert} et \textsc{Emmy} \textsc{Noether}. Elle a joué un rôle central dans le développement des mathématiques du 20 ème siècle,  et elle a, jusqu'à nos jours, des applications en mathématiques, en cryptographie et en physique.

\paragraph*{1.2.1. Définitions initiales}
\subparagraph*{$\ast$} Un corps commutatif $(K,+,\times)$ est un anneau commutatif non nul pour lequel tout élément non nul admet un inverse.\\
i.e. $(K,+)$ est un groupe abélien dont l'élément neutre est noté $0$, \\
$(K-\lbrace0\rbrace,\times)$ est également un groupe abélien et son élément neutre est noté $1$, \\
$\times$ est distributive par rapport à $+$.  

\subparagraph*{$\ast$} Soit $P$ une partie de $K$. Si $P$ est un sous groupe de $(K,+)$ et $P-\lbrace0\rbrace$ muni de $\times$ un sous groupe de $(K-\lbrace0\rbrace,\times)$, alors $P$ est un sous corps de $K$. 
\subparagraph*{$\ast$} Le plus petit sous-corps d'un corps $K$ contenant la partie $P$ est appelée sous-corps de $K$ engendré par $P$.
\subparagraph*{$\ast$} Un corps est dit \textbf{premier} s'il ne contient aucun sous-corps strict. 
\subparagraph*{$\ast$} Si $K$ est un corps, le sous-corps engendré par $1$ est un corps premier. C'est le \textbf{sous-corps premier} de $K$. 
\paragraph*{1.2.2. Rappel et complément sur les idéaux}
\subparagraph*{Définition :} Une partie $I$ d'un anneau $A$ est appelée \textbf{idéal à gauche} (resp. \textbf{à droite}) de $A$, lorsque : \\
 - $I$ est un sous-groupe additif de A\\
 - Pour tout $a$ de $A$ et tout $x$ de $I$, $ax \in I$ (resp. $xa \in I$).
\subparagraph*{$\ast$} On appelle \textbf{idéal bilatère} de $A$, toute partie qui est simultanément idéal à gauche et à droite. 
\subparagraph*{$\ast$} Si $I$ est un idéal bilatère, le groupe quotient $A/I$ peut être muni d'une structure d'anneau, qu'on appelle \textbf{anneau quotient}.  
\subparagraph*{Idéal premier :} Un idéal $I$ de $A$ est dit premier si et seulement si l'anneau quotient $A/I$ est un anneau intègre.
\subparagraph*{Idéal maximal :} Un idéal $I$ de $A$ est dit maximal si et seulement si $I$ est contenu exactement dans deux idéaux, $A$ et lui-même.
\subparagraph*{$\ast$} Dans un anneau commutatif, un idéal bilatère est maximal si et seulement si l'anneau quotient est un corps.
\subparagraph*{$\ast$} Tout idéal maximal est donc premier.
\subparagraph*{$\ast$ Prolongement d'une fonction :} Soient $f$ et $g$ deux fonctions définies respectivement sur $A$ et $C$ à valeurs dans $B$ et $D$ (respectivement) tel que $A \subseteq C$, $f$ coïncide avec $g$ sur $A$ et $B \subseteq D$, on dit alors que $g$ est un prolongement de $f$. 

\subparagraph*{$\ast$ Anneau intègre :} On dit qu'un anneau $A$ est intègre (ou que $A$ est un anneau d'intégrité) s'il est commutatif, non réduit à $\{0\}$, et si le produit de deux éléments non nuls de $A$ est non nul. 

\subparagraph*{$\ast$ Anneau local et corps résiduel d'un anneau :} Un anneau local est un anneau possédant un unique idéal maximal $m(A)$. Le quotient d'un anneau local par son idéal maximal $A/m(A)$ s'appelle le corps résiduel de $A$, on le note souvent $\kappa(A)$. 
\subparagraph*{$\ast$ Anneau principal :} Un anneau principal est un anneau 
intègre dont tout idéal est principal (un idéal principal est un idéal engendré par un seul élément).
\subparagraph*{$\ast$ Anneau noethérien :} On dit qu'un anneau $A$ est noethérien si toute suite croissante d'idéaux de $A$ est stationnaire \footnote{rappelons qu'une suite est stationnaire s'il existe un rang à partir duquel tous les termes de la suite sont égaux}(ou, ce qui revient au même, si tout idéal de $A$ est engendré par un nombre fini d'éléments). 
\subparagraph*{$\ast$ Corps des fractions d'un anneau :} Le corps des fractions d'un anneau intègre $A$ est le plus petit corps commutatif (à isomorphisme près) contenant $A$. 
\subparagraph*{$\ast$ Relation de domination entre anneaux locaux :}
Soient $A$ et $B$ deux anneaux locaux. On dit que $B$ domine $A$ si $A$ est un sous-anneau de $B$ et si $m(A) = A \cap m(B)$. On notera que si $K$ est un anneau, la relation "$B$ domine $A$ " est une relation d'ordre dans l'ensemble des sous-anneaux locaux de $K$. 
\subparagraph*{$\ast$ Anneau intégralement clos :} Un élément $x$ d'un anneau contenant $A$ est dit entier sur $A$ s'il vérifie une équation dite " de dépendance intégrale " : 
\begin{center}
$(*)$ $x^n + a_1x^{n-1} + ... + a_n = 0$, avec $a_i \in A$ $\forall i \in \{1,...,n\}$
\end{center} 
On dit que $A$ est intégralement fermé dans un anneau $B$ le contenant si tout élément de $B$ qui est entier sur $A$ appartient à $A$. On dit que $A$ est intégralement clos s'il est intègre et s'il est intégralement fermé dans son corps des fractions. 
\paragraph*{1.2.3. Extension de corps} 
\subparagraph*{Définition d'une extension de corps :}
Soit $K$ un corps, une extension $L$ de $K$ est un corps tel que $K$ soit sous-corps de $L$. On note $L/K$ le fait que $L$ soit une extension de $K$.\\
Si $L$ est une extension de $K$, il est muni d'une structure de K-espace vectoriel via la multiplication.
\subparagraph*{Degré d'une extension :} Si $L$ est de dimension finie sur $K$, on note $[L:K]$ la dimension du K-espace vectoriel $L$. C'est un entier $>0$ qu'on appelle \textbf{degré de l'extension} $L$ sur $K$. On dit dans ce cas que $L$ \textit{est fini sur} $K$.
\subparagraph*{$\ast$} Si $L$ est fini sur $K$ et $M$ est fini sur $L$ , alors $M$ est fini sur $K$ et on a : $[M:K]=[M:L].[L:K]$ 
\subparagraph*{Élément algébrique - élément transcendant :} Un élément de $L$ est dit \textbf{algébrique} s'il est racine d'un polynôme non nul à coefficients dans $K$, sinon il est \textbf{transcendant}.\\
Plus formellement, si $L/k$ est une extension et $\alpha$ un élément de $L$. On définit un morphisme d'anneaux (également de K-espaces vectoriels) $ \varphi : K[T]\rightarrow L $ par $ P\mapsto P(\alpha) $,\\
- Si $\varphi$ est injectif, on dit que $\alpha$ est transcendant sur $K$.\\
- Si $\varphi$ n'est pas injectif, on dit que $\alpha$ est algébrique sur $K$.   
\subparagraph*{Extension algébrique :} Une extension $L/K$ est dite algébrique si tout élément de $L$ est algébrique sur $K$. 
\subparagraph*{Clôture algébrique :} 
\subparagraph*{$\ast$} Un corps commutatif $K$ est dit \textbf{algébriquement clos} si tout polynôme de degré supérieur ou égal à $1$, à coefficients dans $K$, admet (au moins) une racine dans $K$.
\subparagraph*{$\ast$} Une extension $\overline{K}$ d'un corps commutatif $K$ est dite \textbf{clôture algébrique} si $\overline{K}/K$ est algébrique et $\overline{K}$ est algébriquement close.
\subparagraph*{$\ast$} Tout corps $K$ possède une clôture algébrique, et elle est unique à isomorphisme près (non trivial).


\section*{2. Topologie générale}
\addcontentsline{toc}{section}{2. Topologie générale}

Ce qu'on appelle aujourd'hui topologie générale est l'étude mathématique qualitative des \textit{lieux} et des \textit{relations spatiales} : elle introduit, entre autres, les notions de proximité, frontière, localité, continuité, ainsi que leurs liens mutuels. C'est \textsc{Riemann} qui est généralement considéré comme étant à l'origine de la topologie. 

\subsection*{2.1. Structure topologique }
\addcontentsline{toc}{subsection}{2.1. Structure topologique}
\paragraph*{2.1.1. Ensemble ouvert}

\subparagraph{Définition (Topologie et ensemble ouvert) :}

On appelle structure topologique (ou topologie) sur un ensemble $X$ une structure constituée par la donnée d'un ensemble $O$ de parties de $X$ possédant les propriétés suivantes (dites axiomes des structures topologiques) :
\begin{itemize}
\item[-]($O_1$): Toute réunion d'ensembles de $O$ est un ensemble de O.
\item[-]($O_2$): Toute intersection finie d'ensembles de $O$ est un ensemble de $O$.
\end{itemize}
Les ensembles de $O$ sont appelés ensembles ouverts de la structure topologique définie par $O$ sur $X$.
\subparagraph{$\ast$} On appelle \textbf{ensemble fermé} le complémentaire d'un ensemble ouvert de $X$.
\subparagraph{$\ast$} On appelle \textbf{espace topologique} un ensemble $X$ muni d'une structure topologique $O$ et on note $(X,O)$.

\subparagraph{Remarque :} Certains exposés ajoutent inutilement aux axiomes $O_1$ et $O_2$ l'axiome selon lequel $X$ et $\emptyset $  sont des éléments de $O$, qu'on peut déduire des deux premiers comme suit :
\begin{itemize}
\item[-] L'axiome ($O_1$) implique en particulier que la réunion de la partie vide de $O$, (c'est-à-dire l'ensemble vide) appartient à $O$. \\
En effet pour $(U_i)_i{}_\in{}_I $ une famille de sous ensembles de $O$  tel que $\displaystyle \bigcup_{i\in I} U_i \in O $ quel que soit $I$ alors que pour $I = \emptyset$ on a $\displaystyle \bigcup_{i\in I} U_i = \emptyset$, donc $\emptyset \in O$. 
\item[-] L'axiome ($O_2$) implique que l'intersection de la partie vide de $O$, (c'est-à-dire l'ensemble $X$) appartient à $O$.\\
En effet pour $(U_i)_i{}_\in{}_I $ une famille de sous ensembles de $O$  tel que $\displaystyle \bigcap_{i\in I} U_i \in O $ pour $I$ fini alors si on prend $I = \emptyset$ on a $\displaystyle \bigcap_{i\in I} U_i = X$, donc $X \in O$.  
\end{itemize}

\paragraph{Définition (Recouvrement ouvert) :}
\begin{itemize}
\item[-] Un recouvrement de $X$ est une famille $(U_i)_i{}_\in{}_I $ de $X$  telle que $X =\displaystyle \bigcup_{i\in I} U_i $, si de plus $I$ est un ensemble fini, on dit que $(U_i)_{}i_\in{}_I $ est un recouvrement fini de $X$.
\item[-] Soit $(U_i)_i{}_\in{}_I $ un recouvrement de $X$. Si $\exists J \subseteq I$ tel que $X =\displaystyle \bigcup_{j\in J} U_j $, on dit que $(U_j)_j{}_\in{}_J $ est un sous-recouvrement de $X$. Si de plus, $J$ est fini, on dit alors que $(U_j)_j{}_\in{}_J $ est un sous-recouvrement fini de $X$. 
\item[-] Un recouvrement d'une partie A d'un espace topologique $X$ est dit ouvert si tous les $U_i$ sont ouverts dans $X$.
	\end{itemize}
\paragraph{2.1.2. Exemples de topologies :} Un ensemble $X$ peut généralement être muni de plusieurs topologies distinctes. Parmi celles-ci on peut citer :
\subparagraph{$\ast$ Topologie grossière et topologie discrète :}
La topologie grossière est celle qui comporte le moins d'ouverts : $O = \{X,\emptyset\}$.\\
À l'autre extrémité, la topologie discrète est celle pour laquelle toute partie de $X$ est ouverte : $O = P(X)$, ensemble de toutes les parties de $X$.

\subparagraph{$\ast$ Topologie des espaces métriques :} La notion d'espace métrique fut introduite en 1906 par M. \textsc{Fréchet}, et développée quelques années plus tard par F. \textsc{Hausdorff}. Elle acquit une grande importance après 1920, d'une part grâce aux travaux fondamentaux de S. \textsc{Banach} sur les espaces normés et leurs applications à l'analyse fonctionnelle, de l'autre en raison de l'intérêt que présente la notion de valeur absolue en Arithmétique et en Géométrie algébrique (où notamment la complétion par rapport à une valeur absolue se montre très féconde).
\\
Une distance (ou métrique) sur un ensemble $X$ est une application\\
$d:X \rightarrow \mathbb{R_+}$ possédant, pour tous $ x, y, z \in X $, les propriétés suivantes :
\begin{itemize}
\item[-]$(M_1)$ $d(x, y) = 0 \Leftrightarrow x = y .$
\item[-]$(M_2)$ $d(x, y) = d(y, x) .$
\item[-]$(M_3)$ $d(x, z) \leqslant d(x, y) + d(y, z) $  (\textbf{Inégalité triangulaire})
\end{itemize}

Muni de la distance $d$, $X$ est appelé \textbf{espace métrique}, et on note un tel espace $(X,d)$. Le nombre réel positif $d(x,y)$ est appelé la \textbf{distance entre} $x$ \textbf{et} $y$ dans $X$.\\

Dans un un espace métrique $(X,d)$ , et pour $a \in X$ et $r > 0$. On appelle boule ouverte de centre $a$ et de rayon $r$ l'ensemble :  $B(a, r) = \{ x \in X, d(a, x) < r \}$\\
On peut définir une topologie sur un espace métrique $(X, d)$ comme étant l'ensemble des boules ouvertes de $X$.
\subparagraph{$\ast$ Topologie de $\mathbb{Q}$ et $\mathbb{R}$ :}

Soit $X$ un ensemble = $\mathbb{Q}$ ou $\mathbb{R}$.\\

On appelle intervalle ouvert toute partie de $X$ de l'un des types suivants :
\begin{itemize}
\item[-]$]a,b[=\{ c \in X,a<c<b \}$, $\forall(a,b)\in X^{2}$.
\item[-]$]a,+\infty[ =\{ c\in X, c > a \}$, $\forall a\in X$.
\item[-]$]-\infty,a[ = \{ c \in X ,c < a \}$, $ \forall a\in X$.
\end{itemize}
\subparagraph*{} On appelle ouvert de $X$ toute réunion d'intervalles ouverts. Les ouverts de $X$ définissent une topologie dite \textbf{topologie ordonnée}.

\subsection*{2.2. Voisinages}
\addcontentsline{toc}{subsection}{2.2. Voisinages}
\paragraph*{2.2.1. Définitions}
\subparagraph{$\ast$ Voisinages :} Dans un espace topologique $X$, on appelle voisinage d'une partie A de $X$ tout ensemble qui contient un ensemble ouvert contenant A.\\Les voisinages d'une partie $\{ x \}$ réduite à un seul point s'appellent aussi voisinages du point $x$.\\
On note $V(x)$ la famille des voisinages de $x$.

\subparagraph{$\ast$ Base de voisinages :}
Soient $(X,O)$ un espace topologique et $x \in X$. On appelle \textbf{système fondamental de voisinages} de $x$ ou \textbf{base de voisinages} de $x$, toute famille $B(x)$ de voisinages de $x$ telle que pour tout voisinage $V$ de $x$, il existe $W \in B(x)$ tel que $W \subseteq V$.
\paragraph*{2.2.2. Exemple : Base de voisinages sur la droite rationnelle $\mathbb{Q}$}
\subparagraph{} Sur la droite rationnelle $\mathbb{Q}$, l'ensemble des intervalles ouverts contenant un point $x$ est un système fondamental de voisinages de ce point. Il en est de même pour l'ensemble des intervalles ouverts $]x-\dfrac{1}{n} , x+\dfrac{1}{n}[$, et pour l'ensemble des intervalles fermés $[x-\dfrac{1}{n} , x+\dfrac{1}{n}]$ avec $n \in \mathbb{N^*}$ ou une suite infinie strictement croissante à valeurs dans $\mathbb{N^*}$.
\footnote{On a des résultats analogues pour la droite numérique de $\mathbb{R}$.}

\subsection*{2.3. Adhérence}
\addcontentsline{toc}{subsection}{2.3. Adhérence}
\paragraph{Définition :}
Dans un espace topologique $(X,O)$, on dit qu'un point $x$ est adhérent à un ensemble $A$ lorsque tout voisinage de $x$ rencontre $A$. L'ensemble des points adhérents à $A$
s'appelle adhérence de $A$ et se note $\overline{A}$ .

\subsection*{2.4. Densité}
\addcontentsline{toc}{subsection}{2.4. Densité}
\paragraph{Définition (par l'adhérent) :}
On dit qu'une partie A d'un espace topologique $(X,O)$ est dense dans $X$ \footnote{ou encore est partout dense, lorsqu'il n'en résulte pas de confusion sur $X$} si $\overline{A} = X$.\\
Autrement dit, si pour toute partie ouverte non vide $U$ de $X$, $ U\cap A $ est non vide.

\paragraph{Définition (par les suites) :}
On dit qu'une partie $A$ d'un espace topologique $(X,O)$ est dense dans $X$ si et seulement si pour tout $x$ de $X$ il existe une suite $(x_n)_{n \in \mathbb{N}} \subset A$ qui converge vers $x$. 

\subsection*{2.5. Séparabilité}
\addcontentsline{toc}{subsection}{2.5. Séparabilité}
\paragraph{Définition (espace topologique séparé) :}
Un espace topologique $(X,O)$ est dit séparé \textbf{ou un espace de \textsc{Hausdorf}} s'il vérifie la propriété suivante, appelée \textbf{axiome de \textsc{Hausdorff}:}
\footnote{Certains mathématiciens définissent un espace séparé par six axiomes, alors qu'on peut les résumer à celui de \textsc{Hausdorff}}\\ 
(H): Pour tout couple $(x,y)$ de points distincts $(x \neq y)$ de $X$, il existe un voisinage de $x$ et un voisinage de $y$ disjoints. 

\paragraph{Exemples :}
\subparagraph{La droite des nombres rationnels :} La droite rationnelle $\mathbb{Q}$ est séparée, car si $x$, $y$ sont deux nombres rationnels tels que $x < y$, et $z$ un nombre rationnel tel que $x < z < y$, les voisinages respectifs $]-\infty,z[$ et $]z,+\infty[$ de $x$ et $y$ ne se rencontrent pas.
\subparagraph{Les espaces métriques :} Tout espace métrique $(X, d)$ est séparé : soient en effet $x$ et $y$ deux points distincts de $X$; alors $B(x,\dfrac{d(x, y)}{2})$ et $B(y,\dfrac{d(x,y)}{2} )$ sont des voisinages, respectivement de $x$ et de $y$, disjoints. 

\subsection*{2.6. Complétude}
\addcontentsline{toc}{subsection}{2.6. Complétude}
\paragraph{Définition (Suite de \textsc{Cauchy}) :}
Soit (X, d) un espace métrique.\\
Une suite $(x_n)_{n \geqslant_0}$ dans $(X, d)$ est dite suite de Cauchy si pour tout $\epsilon > 0$, il existe $n_0 \in \mathbb{N}$ tel que pour tout $n, m \in \mathbb{N}$ vérifiant $n \geqslant n_0$ et $m \geqslant n_0$, on ait $d(x_m, x_n) < \epsilon$.\\
Il revient au même de dire que pour tout $\epsilon > 0$, il existe $n_0 \in \mathbb{N}$ tel que pour tout $n \in \mathbb{N}$ vérifiant $n \geqslant n_0$ et pour tout $p \in \mathbb{N}$, on a $d(x_{n+p}, x_n) < \epsilon$.
\paragraph{Définition (Espace complet) :}
Soit $ (X, d) $ un espace métrique. On dit que $(X, d)$ est \textbf{complet} si toute suite de \textsc{Cauchy} dans $(X, d)$ est convergente.

\subsection*{2.7. Complétion}
\addcontentsline{toc}{subsection}{2.7. Complétion}
La complétion est une notion indispensable pour notre sujet, c'est, en de simples mots, l'extension d'un espace uniforme $E$ (souvent un espace métrique pour nous) en un autre espace dans lequel s'envoie $E$, qui est complet et qui est minimal parmi ceux-ci.
\subparagraph{Construction :} Un espace métrique $E$ non complet possède au moins une suite de Cauchy qui ne converge pas. Le complété de $E$ est obtenu en ajoutant les limites des suites de Cauchy.
Par exemple, $\mathbb{Q}$, muni de la métrique de distance, n'est pas complet car il existe des suites de Cauchy qui ne convergent pas, en l'occurrence, 1, 1.4, 1.41, 1.414, ... ne converge pas parce que $\sqrt{2}$ n'est pas rationnel. Le complété de $\mathbb{Q}$ est $\mathbb{R}$. Remarquons que le complété dépend de la métrique. Par exemple, pour tout $p$ premier, $\mathbb{Q}$ peut être muni de la norme p-adique, et alors le complété de l'ensemble des rationnels est l'ensemble des nombres p-adiques.\\
Techniquement parlant, le complété de $E$ est l'ensemble des limites des suites de \textsc{Cauchy}. $E$ est contenu dans cet ensemble, en considérant les suites constantes.\\
Lorsque $E$ est séparé - c'est en particulier le cas si $E$ est métrique - $E$ est un sous-espace de son complété.

\subsection*{2.8. Compacité}
\addcontentsline{toc}{subsection}{2.8. Compacité}
\paragraph{Définition (Espace topologique compact) :}
Soit $X$ un espace topologique séparé. On dit que $X$ est compact s'il vérifie la propriété dite \textbf{propriété de Borel-Lebesgue :}\\
De tout recouvrement ouvert de $X$, on peut extraire un sous-recouvrement fini.\\
Autrement dit, pour toute famille d'ouverts $(U_i)_{i\in I}$ de $X$ telle que $X =\displaystyle \bigcup_{i\in I} U_i $, il existe un sous-ensemble fini $J$ de $I$ tel que $X =\displaystyle \bigcup_{j\in J} U_j $.

\paragraph{Proposition :}
\subparagraph*{} Soit $(X,O)$ un espace topologique séparé. Les propriétés suivantes sont équivalentes :
\begin{itemize}
\item[$(i)$] $X$ est compact.
\item[$(ii)$]De toute famille de fermés de $X$ dont l'intersection est vide, on peut extraire une sous-famille finie dont l'intersection est vide.
\item[$(iii)$] Toute famille de fermés de $X$ dont toute sous-famille finie est d'intersection non vide, est elle-même d'intersection non vide.
\end{itemize}

\paragraph{Propriétés:}
\begin{itemize}
\item[-] Si $(K_i)_{i \in I}$ est une famille de parties compactes de $X$, alors l'intersection \textbf{quelconque} $\displaystyle \bigcap_{i\in I} K_i $ est compacte, toute réunion \textbf{finie} de parties compactes est également compacte. 
\item[-] Si $X$ est compact et $A$ est fermée dans $X$, alors $A$ est compacte.
\item[-] Tout intervalle fermé et borné de $\mathbb{R}$ est compact (\textbf{Théorème de \textsc{Heine}}). 
\end{itemize}

\paragraph{Définition (Espace topologique localement compact) :}
Soit $X$ un espace topologique séparé. On dit que $X$ est localement compact si tout point de $X$ admet un voisinage compact.
\paragraph{Exemples :} 
\begin{itemize}
\item[a -] Tout espace compact est localement compact. 
\item[b -] La droite réelle $\mathbb{R}$ et plus généralement, l'espace $\mathbb{R}^n$, sont localement compacts.
\item[c -] Tout espace muni de la topologie discrète est localement compact. En particulier, $\mathbb{Z}$ (muni de sa topologie de sous-espace de $\mathbb{R}$) est localement compact. 
\item[d -] L'ensemble des rationnels $\mathbb{Q}$, pour sa topologie de sous-espace de $\mathbb{R}$, n'est pas localement compact.
\end{itemize}
\paragraph{Propriétés :}
\begin{itemize}
\item[-]Les parties ouvertes et les parties fermées d'un espace localement compact sont localement compactes.
\item[-] Dans un espace séparé, l'intersection de deux parties localement compactes est localement compacte. Par contre, la réunion de deux parties localement compactes n'est pas en général localement compacte.
\end{itemize} 
\paragraph{Définition (Espace topologique précompact) :} Un espace métrique est dit précompact si son complété est 
compact. 
\paragraph*{Proposition :}Un espace métrique est précompact si et seulement si, pour tout $\varepsilon > 0$, cet espace peut être recouvert par un nombre fini de boules ouvertes 
de rayon $\varepsilon$. 

\subsection*{2.9. Compactification}
\addcontentsline{toc}{subsection}{2.9. Compactification}
\paragraph{2.9.1. Espace compactifié }
\paragraph{Définition (Compactifié d'un espace topologique) :}
Une compactification d'un espace topologique $X$ est la donnée d'un couple $(\widehat{X}, f)$ constitué d'un espace compact $\widehat{X}$ et d'un homéomorphisme $f$ de $X$ sur un sous-ensemble dense de $\widehat{X}$.\\
On identifie fréquemment l'espace $X$ avec $f(X) \subset \widehat{X} $. Et on dit simplement que $\widehat{X}$ est un compactifié de $X$.
Si $X$ est déjà compact, l'espace $\widehat{X}$ est homéomorphe à $X$, et dans ce cas, il ne sert à rien de parler de compactifié de $X$.
\paragraph{2.9.2. Compactification d'\textsc{Alexandroff}}
\paragraph{Ingrédients pour construire un espace compactifié d'\textsc{Alexandroff} :\\}
Soit $(X, O )$ un espace localement compact. Ajoutons à $X$ un nouveau point, noté $w$ ou $\infty$ et appelé point à l'infini, et considérons l'ensemble :\\
$X_\infty = X \cup \{\infty\}$.\\
On dit souvent que $ \{ \infty \}$ est le \textbf{point à l'infini} de $X_\infty $, et que $X_\infty $ résulte de $X$ par adjonction d'un point à l'infini.\\
Soit $O_\infty $ l'ensemble des parties de $X_\infty $ défini par : une partie $ U$ de $X_\infty $ appartient à $O_\infty $ si  $U$ appartient à $T$ , ou bien si $U$ est le complémentaire dans $X_\infty $ d'un compact de $(X, O )$. On peut vérifier que $O_\infty $ est une topologie sur $X_\infty $.\\
Il est clair que l'espace $X_\infty $, muni de la topologie $T_\infty $, est séparé.\footnote{On peut construire l'espace $X_\infty $ pour n'importe quel espace topologique $X$, mais dans ce cas, $X_\infty $ n'est pas toujours séparé. En fait, $X_\infty $ est séparé si et seulement si $X$ est localement compact.}
\paragraph{Définition (compactifié d'\textsc{Alexandroff}) :}
Soit $(X, O )$ un espace localement compact. L'espace compact $(X_\infty , O_\infty )$ est appelé le \textbf{compactifié d'Alex-\\
androff} de $(X, O )$.

\subsection*{2.10. Compacité des espaces quotients}
\addcontentsline{toc}{subsection}{2.10. Compacité des espaces quotients}
\paragraph{2.10.1. Définitions :}
\paragraph{Définition (Application propre) :}
Soient $X, Y$ des espaces topologiques, avec $X$ séparé et $f: X \rightarrow Y$ une application.\\
On dit que $f$ est propre si $f$ est continue et fermée (envoie les fermés du premier espace vers les fermés du second) et si pour tout $y \in Y$, $f^{-1}(\{y\})$ est une partie compacte de $X$.
\paragraph{Définition (Saturé pour une relation d'équivalence) :}
Soit $X$ un ensemble et soit A une partie de $X$, l'ensemble des points de $X$ qui sont équivalents à un point de A est appelé \textbf{le saturé} de A pour la relation d'équivalence.\\
 Le saturé de A s'écrit alors $\pi^{-1}(\pi(A))$ où $\pi$ est la projection canonique de $X$ sur $X/\Re$.
\paragraph{Définition (Relation d'équivalence ouverte/fermée) :}
soient $\Re$ une relation d'équivalence ouverte (resp. fermée) dans un espace topologique $(X,O)$, $f$ l'application canonique $f:X \rightarrow X/\Re$, A une partie de $X$. Supposons que l'une des deux conditions suivantes soit vérifiée :
\begin{itemize}
\item[a - ]A est ouvert (resp. fermé) dans $X$.
\item[b - ]A est saturé pour $\Re$.
\end{itemize}
Sous ces conditions, la relation $\Re_A$ induite sur A est ouverte (resp. fermée) et l'application canonique de $A/\Re_A$ sur $f(A)$ est un homéomorphisme.
\paragraph{Définition (Graphe d'une relation d'équivalence) :}

Soient $\Re$ une relation d'équivalence sur un espace topologique $(X,O)$, on appelle graphe de $\Re$ l'ensemble $G(\Re) = \{ (x, y) \in X^2 ; x \Re y \}$.

\paragraph{Définition (Séparation d'un espace quotient) :}

Cherchons des conditions pour qu'un espace quotient $X/\Re$ soit séparé (auquel cas on dit que la relation d'équivalence $\Re$ est séparée). En premier lieu, si $X/\Re$ est séparé, les ensembles réduits à un point dans $X/\Re$ sont fermés, donc toute classe d'équivalence suivant $\Re$ est fermée dans $X$. Cette condition nécessaire n'est pas suffisante; la définition des ensembles ouverts dans $X/\Re$ donne la
condition nécessaire et suffisante suivante :
\subparagraph{}
Pour que $X/\Re$ soit séparé, il faut et il suffit que deux classes d'équivalence distinctes dans $X$ soient respectivement contenues dans deux ensembles ouverts saturés sans point commun. 
\subparagraph{$\ast$} Autre condition plus maniable : Pour qu'un espace quotient $X/\Re$ soit séparé, il est nécessaire que le graphe $G(\Re)$ de $\Re$ soit fermé dans $X^2$. Cette condition est suffisante lorsque la relation $\Re$ est ouverte.
\subparagraph{$\ast$} Dans un espace $X$ régulier (\textit{séparé et l'ensemble des voisinages fermés d'un point quelconque de $X$ est un système fondamental de voisinages de ce point}), toute relation d'équivalence à la fois ouverte et fermée est séparée.

\paragraph{2.10.2. Espaces quotients compacts}
\subparagraph*{}

Soient $X$ un espace compact, $\Re$ une relation d'équivalence dans $X$, $G(\Re)$ son graphe dans $X \times X$ et $q : X \rightarrow X/\Re$ l'application quotient. Les propriétés suivantes sont équivalentes.
\begin{itemize}
\item[$(i)$] L'espace topologique quotient $X/\Re$ est séparé.
\item[$(ii)$] Le graphe $G(\Re)$ est fermé dans $X \times X $.
\item[$(iii)$] La relation $\Re$ est fermée.
\item[$(iv)$]  L'application quotient $q$ est propre.

\end{itemize}
En outre, lorsque une de ces propriétés est vérifiée, alors l'espace $X/\Re$ est compact.

\paragraph{2.10.3. Espaces quotients localement compacts}
\subparagraph*{}Soient $X$ un espace compact, $\Re$ une relation d'équivalence dans $X$, $G(\Re)$ son graphe dans $X^2$ et $q : X \rightarrow X/{\Re} $ l'application quotient.
\subparagraph*{} Soient $X_{\infty}$ = $X \cup \{\infty\}$ le compactifié d'\textsc{Alexandroff} de $X$; et $\Re  _{\infty }$ la relation d'équivalence dans $X$ dont le graphe est $G(\Re ) \cup \{ -\infty ,+ \infty \}$. Les propriétés suivantes sont équivalentes :
\begin{itemize}
\item[$(i)$] L'application quotient $q$ est propre.
\item[$(ii)$] Le saturé pour $\Re$ de toute partie compacte de $X$ est un ensemble compact.
\item[$(iii)$] La relation $\Re_{\infty}$ est fermée.
\item[$(iv)$]  La restriction à $G(\Re )$ de l'application $(x, y)\mapsto y$ de $X^2$ dans $X$ est propre.
\item[$(v)$] La relation $\Re $ est fermée et les classes suivant $\Re$ sont compactes.
\end{itemize}
En outre, lorsqu'une de ces propriétés est vérifiée, alors l'espace $X/\Re$ est localement compact.

\paragraph{2.10.4. Corollaire et proposition :}
\subparagraph*{$\ast$} Soient $X$ un espace séparé, $Y$ un espace topologique, $f: X \rightarrow Y$ une application propre. Pour que $X$ soit compact (resp. localement compact), il faut et il suffit que $f(X)$ soit compact (resp. localement compact), et il suffit que Y soit compact (resp. localement compact).
\subparagraph*{$\ast$}
Soient $X$ un espace localement compact, $\Re$ une relation d'équivalence ouverte et séparée dans $X$, $f$ l'application canonique $X \rightarrow X/\Re$. Alors $X/\Re$ est localement compact, et pour toute partie compacte $K'$ de $X/\Re$, il existe une partie compacte $K$ de $X$ telle que $f(K) = K'$.
\newpage 

\vspace*{\stretch{1}}
\begin{center}
\includegraphics[scale=0.8]{2.png}  
\end{center}
	\begin{center}
		{\LARGE {\huge Première partie :} }\\
		\textit{  }\\
		\textit{  }\\
		{\LARGE {\Huge Corps valués et propriétés}}
\begin{center}
\includegraphics[scale=1]{3.png}  
\end{center}
\addcontentsline{toc}{part}{Première partie : Corps valués et propriétés}
	\end{center}
\vspace*{\stretch{1}}
\newpage

\subparagraph*{}
\section*{1. La valeur absolue}
\addcontentsline{toc}{section}{1. La valeur absolue}
\subsection*{1.1. Corps valué}
\addcontentsline{toc}{subsection}{1.1. Corps valué}
\paragraph{Définition (\textbf{Valeur absolue archimédienne}) :}
\subparagraph{} On appelle valeur absolue sur un corps $K$ une application $ x \rightarrow | x | $ de $K$
dans $\mathbb{R_+}$, satisfaisant aux conditions suivantes :
\begin{itemize}
\item[$VA_1$] $\forall x \in K$, $|x| = 0  \Leftrightarrow   x = 0 $.
\item[$VA_2$] $\forall x,y \in K$, $|xy| = |x|.|y|$ .
\item[$VA_3$] $\forall x,y \in K$, $|x + y| \leqslant |x| + |y|$ .
\end{itemize}
En outre, une telle application est appelée \textit{valeur absolue archimédienne}. 

\paragraph{Définition (\textbf{Valeur absolue non-archimédienne}) :}
\subparagraph{} Une valeur absolue sur un corps $K$ est dite non-archimédienne si elle vérifie de plus l'inégalité suivante :
  \begin{center}
  $\forall x,y \in K$, $|x + y| \leqslant \max \{|x|;|y|\}$    (\textit{Inégalité ultra-métrique})
  \end{center}

\paragraph{Remarque :}
\subparagraph{$\ast$} D'après $VA_2$, on a $|x| = |1|.|x|$, et comme il existe d'après $VA_1$ au moins un $x$ tel que $|x| \ne 0$, on a $|1|=1$; on en tire $1=|-1|^2$, d'où également $\mid -1 \mid$ = 1, et par suite on en conclut que :
\begin{center}
 $ \forall (x, y)\in K$ on a $|x - y| \leqslant |x| + |y|$ 
\end{center}
\subparagraph*{$\ast$} Soit $|.|$ une valeur absolue sur un corps $K$. L'inégalité triangulaire peut être généralisée (par récurrence) à $n$ éléments de $K$ par :
\begin{center}
$|{\sum \limits_{\underset{}{i=1}}^n x_i}| \leqslant \sum \limits_{\underset{}{i=1}}^n |x_i|$
\end{center}

\subparagraph*{$\ast$} Soit $|.|$ une valeur absolue sur un corps $K$. L'inégalité ultramétrique peut être généralisée (encore par récurrence) à $n$ éléments de $K$ par :
\begin{center}
$|{\sum \limits_{\underset{}{i=1}}^n x_i}| \leqslant \max \limits_{\underset{}{i\in \{1,...,n\}}} \{|x_i|\}$
\end{center}

\subparagraph*{$\ast$} Soit $|.|$ une valeur absolue sur un corps $K$. On a le résultat suivant :
\begin{center}
$\forall (x, y) \in K^2$, $||x| - |y|| \leqslant |x - y|$.
\end{center}
On l'utilise souvent pour montrer la continuité de la fonction valeur absolue.\\
On note $Val(K)$ (resp. $Val^u(K)$) l'ensemble des valeurs absolues (resp. des valeurs absolues ultramétriques) de $K$.

\paragraph{Définition (\textbf{Corps valué}) :}
\subparagraph{} On appelle corps valué un corps $K$ muni de la structure définie par la donnée d'une valeur absolue sur $K$. Et on le note $(K,+,.,| |)$.\\

Si la valeur absolue est archimédienne, $(K,+,.,| |)$ est appelé corps valué archimédien.\\

Si la valeur absolue est non-archimédienne, $(K,+,.,| |)$ est appelé corps valué non-archimédien.
\paragraph{Exemples :}
\subparagraph{$\ast $ Valeur absolue impropre / triviale :}
Soit $K$ un corps quelconque ; pour tout $ x \in K$, posons $\mid x \mid$ = 1 si
$x \neq 0 $, et $\mid 0 \mid$ = 0 ; l'application $x \rightarrow \mid x \mid$ ainsi définie est une valeur absolue sur $K$, dite \textbf{valeur absolue impropre}.

Tout corps est donc valué car il possède au moins la valeur absolue triviale.

\subparagraph{Remarque :}
Si, dans un corps valué, $x$ est racine de l'unité, on a $\vert x \vert = 1$, car si $x^n = 1$ pour un entier $n > 0$, on en tire $\vert x \vert^n = 1$ d'où $\vert x \vert = 1$. En particulier, la seule valeur absolue sur un corps fini est la valeur absolue impropre, puisque tout élément $\neq 0$ du corps est racine de l'unité.

\subparagraph{$\ast $ Valeur absolue d'un nombre réel/rationnel :} La fonction définie de $\mathbb{R}$ ou $\mathbb{Q}$ vers $\mathbb{R}^+$ (ou resp. $\mathbb{Q}^+$) :
\begin{center}
$|x| = \left\{ \begin{array}{rcl}
x & \mbox{si} & x \geqslant 0 \\ - x & \mbox{si} & x \leqslant 0 \end{array}\right.$
\end{center}
est bel et bien une valeur absolue sur $\mathbb{R}$ (resp. $\mathbb{Q}$).
\subparagraph{$\ast $ Valeur absolue définie par une valuation :}
Sur un corps $K$, une valuation réelle est une fonction $v$, définie dans $K^*$, à
valeurs dans $\mathbb{R}$, satisfaisant aux conditions suivantes:
\begin{itemize}
\item[$a -$] Pour $x \in K^*$, $y \in K^*$, $v(xy)= v(x)+ v(y)$.
\item[$b -$] si en outre $x + y \neq 0$, $v(x + y) \geqslant inf(v(x),v(y))$.
\end{itemize}
Si $a$ est un nombre réel quelconque $> 1$, on définit alors sur $K$ une valeur absolue en posant
$\mid x \mid$ = $a^-{}^v{}^({}^x{}^)$ pour $x \neq 0$, et $\mid x \mid$ = 0 pour $x=0$. 
\subparagraph*{}
En effet, de la relation $v(xy) = v(x) + v(y)$ pour
$x \neq 0$ et $y \neq 0$, on déduit la relation $\mid xy \mid$ = $\mid x \mid.\mid y \mid$ pour ces valeurs de $x$ et $y$, et cette relation est trivialement vérifiée si l'un des éléments $x, y$ est nul, de même, de $v(x + y) \geqslant inf(v(x), v(y))$ pour $x \neq 0$, $y \neq 0$ et $x + y \neq 0$, on déduit
\begin{center}
 $\mid x +y \mid \leqslant \sup( \mid x \mid, \mid y \mid ) \leqslant \mid x \mid + \mid y \mid \ $
 \end{center} 
Et ces inégalités sont encore vérifiées si l'un des éléments $x$, $y$, $x + y$ est nul. 
\subparagraph{$\ast $ Valeur absolue $p$-adique sur $\mathbb{Q}$ :}
C'est un cas particulier de l'exemple précédent, si $v_{p}(x)$ est la valuation \textit{p-adique} sur le corps $\mathbb{Q}$ des nombres rationnels
(Exposant de $p$ dans la décomposition de $x$ en produit de facteurs premiers), la valeur absolue correspondante $\mid x \mid_p$ = $p^{-v_p(x)}$ est dite \textbf{valeur absolue p-adique sur le corps $\mathbb{Q}$}.
\\
Plus précisément :\\

Soit $p$ un nombre premier. On considère la fonction $\mid\mid_p$ définie par $\forall x \in \mathbb{Q}^*$,\\

$ \mid x \mid_p = \left\{ \begin{array}{rcl}
p^{-k} & \mbox{si}
& x = p^k\dfrac{a}{b} \\ 0 & \mbox{si} & x = 0 \\
\end{array}\right.$
avec $(a,p) = (b,p) = 1$, où $(a,b)$ désigne le PGCD de $a$ et $b$.

\paragraph*{Lemme :}

L'application $x \rightarrow \mid x\mid_p$ est une valeur absolue sur $\mathbb{Q}$. Elle vérifie, en outre, l'inégalité (ultramétrique) :
 \begin{center}
 $\forall x,y \in \mathbb{Q}$, $\mid x + y\mid_p \leqslant max(\mid x \mid_p , \mid  y\mid_p)$.     $ (\ast) $
 \end{center}
Cette valeur absolue $\mid x\mid_p$ est appelée \textbf{valeur absolue p-adique de $\mathbb{Q}$}.\\

\paragraph*{Preuve :}
\subparagraph*{} Par définition de $\mid x\mid_p$, l'égalité $VA_1$ est satisfaite.\\
Si $x = p^k\dfrac{a}{b} $ et $y = p^{k'}\dfrac{d}{c} $ avec $(a,p) = (b,p) = (c,p) = (d,p) = 1$ alors $xy = p^{k+k'} \dfrac{ac}{bd}$ et $(ac,p) = (bd,p)$, ce qui nous donne
$|xy|_p = p^{-(k+k')} = p^{-k_p -k'} = |x|_p.|y|_p$ et démontre l'égalité $VA_2$.\\
D'autre part, en supposant que $k' \leqslant k$ (sinon on échange $x$ et $y$), on a 
\begin{center}
$|y|_p = max_p(|x|_p , |y|_p)$ et $x + y = p^k(\dfrac{a}{b} + p^{k' -k} \dfrac{d}{c} ) = p^k(\dfrac{ac + dbp^{k'-k}}{bc})$
\end{center}

Puisque $ac + dbp^{k'-k}$ est un entier, il existe deux entiers $r$ et $a'$ tels que $ac + dbp^{k'-k}$ = $p^ra'$ avec $(a',p) = 1$ et $r > 0$.\\
En outre, $p$ est premier à $b$ et $c$ donc $p$ est premier à $bc$, ce qui nous donne l'inégalité \begin{center}
$|x + y|_p = p^{-(k+r)} \leqslant p^{-k} = max(|x|_p , |y|_p) \leqslant |x|_p + |y|_p$
\end{center} 
qui démontre les inégalités $VA_3$. et $ (\ast) $
\begin{flushright}
$\blacksquare$
\end{flushright}
\subsection*{1.2. Topologie des corps valués}
\addcontentsline{toc}{subsection}{1.2. Topologie des corps valués}
\subsubsection*{1.2.1. Métrique définie par une valeur absolue}
Soit $( K , |.|)$ un corps valué, alors on peut munir $K$ d'une métrique $d$ dite
associée à la valeur absolue sur $K$ définie de la façon suivante :
\begin{center}
$\forall(x, y)\in K^2$, $d(x, y) =\mid x - y\mid$ 
\end{center}
En effet:
\subparagraph*{$(M_1)$:}  $d(x, y) = 0 \Leftrightarrow \mid x - y \mid = 0$.
d'où d'après $VA_1$: $ |x - y| = 0  \Leftrightarrow   x = y $
\subparagraph*{$(M_2)$:} $d(x, y) = \mid x - y \mid = \mid y - x \mid = d(y, x) .$ d'après la première Remarque.
\subparagraph*{$(M_3)$:} $d(x, z) = \mid x - y \mid =\mid x - z + z - y \mid \leqslant \mid x - z \mid + \mid z - y \mid = d(x, y) + d(y, z)$
\subparagraph*{}
Et alors $(K, d)$ est un espace métrique.
Si la valeur absolue de $K$ est non-archimédienne, alors la distance $d$ satisfait l'inégalité ultramétrique suivante :
\begin{center}
 $\forall (x, y, z) \in K^3, d (x, z) \leqslant max(d(x, y) , d(y, z))$
 \end{center}
\subsubsection*{1.2.1. Topologie d'un corps valué }
Soient $(K, |.|)$ un corps valué, $x \in K$ et $\epsilon > 0$. L'ensemble :
\begin{center}
 $B_o(x,\epsilon) = \{y \in K : |y - x| < \epsilon \}$
 \end{center} est la boule ouverte de centre $x$ et de rayon $\epsilon$. Et l'ensemble :
\begin{center}
 $B_f(x,\epsilon) = \{y \in K : |y - x| \leqslant \epsilon \}$
 \end{center} est la boule fermée de centre $x$ et de rayon
$\epsilon$.
\subparagraph*{} Pour $x \in K$ et $\epsilon > 0 $ l'ensemble de toutes les boules ouvertes $B_o(x, \epsilon)$ est une \textit{base de la topologie} métrique définie sur $K $.

\paragraph{Définition (\textbf{Corps topologique}) :}
\subparagraph*{} Soit $(K, |.|)$ un corps valué, alors muni de la topologie définie par la distance $d$, $K$ est appelé \textbf{corps topologique}. C'est-à-dire les applications suivantes sont continues :
\begin{multicols}{2}
\begin{itemize}
\item[1 -]$\begin{array}{ccccc}
f_1 & : & K^2 & \to & K \\
 & & (x,y) & \mapsto & x+y \\
\end{array}$
\item[2 -]$\begin{array}{ccccc}
f_2 & : & K^2 & \to & K \\
 & & (x,y) & \mapsto & x.y \\
\end{array}$
\item[3 -]$\begin{array}{ccccc}
f_3 & : & K^2 & \to & K \\
 & & (x,y) & \mapsto & x^{-1} \\
\end{array}$
\item[4 -]$\begin{array}{ccccc}
f_4 & : & K^2 & \to & K \\
 & & (x,y) & \mapsto & -x \\
\end{array}$
\end{itemize}
\end{multicols}

\subsection*{1.3. Équivalence des valeurs absolues\footnote{Dans cette sous-partie, $p$ désignera systématiquement un nombre premier.}}
\addcontentsline{toc}{subsection}{1.3. Équivalence des valeurs absolues}
\paragraph*{Définition (\textbf{Deux valeurs absolues équivalentes}) :}
\subparagraph*{} Deux valeurs absolues sur un corps $K$ sont dites équivalentes si, et seulement si, elles définissent la même topologie (naturelle) sur $K$. Autrement dit si, et seulement si, leurs distances associées respectives induisent la même topologie sur $K$.

\paragraph*{Remarque :} 
\subparagraph*{} La relation "\textit{être équivalent}" est une relation d'équivalence sur l'ensemble des valeurs absolues d'un corps.

\paragraph*{Lemme :}
\subparagraph*{} Soit $K$ un corps et $|.|_1$ , $|.|_2$ deux valeurs absolues sur $K$.
Les valeurs absolues $|.|_1$ et $|.|_2$ sont équivalentes ssi pour toute suite $(x_n)_{n\in \mathbb{N}}$ de $K$ \begin{center}
$(|x_n|_1 \rightarrow 0$ lorsque $n \rightarrow + \infty)
\Leftrightarrow
(|x_n|_2 \rightarrow 0$ lorsque $n \rightarrow + \infty)$
\end{center}

\paragraph*{Preuve :}
\subparagraph*{} – Supposons que les valeurs absolues $|.|_1$ et $|.|_2$ sont équivalentes.\\
Soit $(x_n)_{n>0}$ une suite de $K$ convergeant vers 0 pour la distance $d_1$. Alors, pour tout voisinage ouvert $V$ de 0 (pour la distance $d_1$), il existe un rang $n_0$ tel que $\forall n>n_0$, $x_n \in V$. Or tout ouvert pour $d_2$ est un ouvert de $d_1$, donc pour tout voisinage ouvert $V$ de 0 (pour la distance $d_2$), il existe un rang $n_0$ tel que $\forall n>n_0$, $x_n \in V$ ce qui démontre que $(x_n)_{n>0}$ converge vers 0 pour $d_2$.
\subparagraph*{} – Supposons que pour toute suite $(x_n)_{n\in \mathbb{N}}$ de $K$ 
\begin{center}
$(|x_n|_1 \rightarrow 0$ lorsque $n \rightarrow + \infty)
\Leftrightarrow
(|x_n|_2 \rightarrow 0$ lorsque $n \rightarrow + \infty)$
\end{center}

Démontrons que les ouverts pour $d_1$ sont les ouverts pour $d_2$ revient à démontrer que les fermés pour $d_1$ sont
les fermés de $d_2$ (le complémentaire d'un ouvert est un fermé et vice-versa). La caractérisation séquentielle
des fermés montre que $F$ est fermé ssi pour toute suite $(x_n)_{n>0}$ d'éléments de $F$ convergeant vers $x$ dans $K$ pour la distance $d_1$ alors $x \in F$ .\\
Soit $F$ un fermé pour la distance $d_1$ et soit $(x_n)_{n>0}$ une suite d'éléments de $F$ convergeant vers $x \in K$ pour la distance $d_2$. Alors \begin{center}
$d_{||_2}(x_n,x) = |x_n - x|_2 \rightarrow 0 \Leftrightarrow |x_n - x|_1 = d_{||_1}(x_n,x) \rightarrow 0$.
\end{center}
On en déduit que la suite $(x_n)_{n>0}$ converge vers $x$ dans $K$ pour la distance $d_1$ et, puisque $F$ est fermé pour la distance $d_1$, $x \in F$.\\
L'ensemble $F$ est donc fermé pour la distance $d_2$. En échangeant les rôles de $d_1$ et $d_2$, on conclut.
\begin{flushright}
$\blacksquare$
\end{flushright}
\paragraph*{Théorème :} 
\subparagraph*{} Soient $|.|_1$ et $|.|_2$ deux valeurs absolues sur $K$, alors $|.|_1$ et $|.|_2$, sont équivalentes ssi il existe un réel positif $a$ tel que
$\forall x \in K$, $|x|_1 = |x|^a_2$

\paragraph*{Preuve :}
\subparagraph*{} L'implication réciproque est évidente grâce au lemme précédent.\\
Pour l'implication directe, soit $x$ un élément de $K$ tel que $|x|_1 < 1$.\\
La suite $(x^n)_{n\in \mathbb{N}}$ converge vers 0 $(|x^n|_1 = |x|^n_1)$ dans $(K,d_{||_1})$ donc elle converge vers 0 dans $(K,d_{||_2})$ c'est-à-dire
$|x|^n_2 = |x^n|_2 \rightarrow 0$ lorsque $ n \rightarrow +\infty$
d'où $|x|_2 < 1$. En échangeant le rôle joué par les deux valeurs absolues, on obtient que
\begin{center}
 $\forall x \in K, (|x|_1 < 1) \Leftrightarrow (|x|_2 < 1)$
\end{center}
ensuite en remplaçant $x$ par $\dfrac{1}{x} (x \neq 0)$, on obtient
\begin{center}
$\forall x \in K, (|x|_1 > 1) \Leftrightarrow (|x|_2 > 1)$
\end{center}
et par conséquent
\begin{center}
$\forall x \in K, (|x|_1 = 1) \Leftrightarrow (|x|_2 = 1)$
\end{center}
Ainsi si $|.|_1$ est la valeur absolue triviale, on en déduit que $|.|_2$ est également la valeur triviale.\\
Supposons $|.|_1$ ne soit pas triviale : il existe $x_0 \in \mathbb{K}$ tel que $|x_0|_1 > 1$ (donc $|x_0|_2 > 1$) ce qui implique qu'il existe
\begin{center}
$a \in \mathbb{R}^+$ tel que $|x_0|_1 = |x_0|^a_2$  $(a =  \dfrac{ln |x_0|_1}{ln |x_0|_2} > 0)$
\end{center}
Soit $x \in K$ tel que $|x|_1 > 1$. Considérons le réel $b$ pour lequel $|x|_1 = |x_0|^b_1$. Pour tout rationnel $\dfrac{p}{q} > b$, on a les équivalences suivantes :

\begin{center}
$|x|_1 < |x_0|^{\dfrac{p}{q}}_1 \Leftrightarrow |x^q|_1 < |x^p_0|_1 \Leftrightarrow \dfrac{|x^q|_1}{|x^p_0|_1}
< 1 \Leftrightarrow |x^q|_2 < |x^p_0|_2 \Leftrightarrow |x|_2 < |x_0|^{\dfrac{p}{q}}$
\end{center}
En faisant tendre $\dfrac{p}{q}$ vers $b$ dans $\mathbb{R}$, on obtient que $|x|_2 \leqslant |x_0|^b_2$.
En appliquant le même raisonnement à un rationnel $\dfrac{p}{q} < b$ puis en passant à la limite, on obtient que $|x|_2 \geqslant |x_0|^b_2$ ce qui nous fournit l'égalité
\begin{center}
$|x|_2 = |x_0|^b_2 = |x_0|^{\dfrac{b}{a}}_1= |x|^{\dfrac{1}{a}}_1 \Rightarrow |x|_1 = |x|^a_2$
\end{center}
valable pour tout élément $x$ de $K$ tel que $|x|_1 > 1$. En remplaçant $x$ par $\dfrac{1}{x}$ et en utilisant la multiplicativité des valeurs absolues, on en déduit que
\begin{center}
$\forall x \in K$, tel que $|x|_1 \neq 1, |x|_1 = |x|^a_2$ .
\end{center}
Soit $x \in K$ tel que $|x|_1 = 1$. L'élément $\dfrac{x}{x_0}$ qui vérifie $\mid \dfrac{x}{x_0}\mid_1 = \dfrac{|x|_1}{|x_0|_1} = \dfrac{1}{|x_0|_1} < 1$ donc on a :
\begin{center}
$\mid \dfrac{x}{x_0}\mid_1 =\mid \dfrac{x}{x_0}\mid^a_1 \Leftrightarrow |x|_1 = |x|^a_2$
\end{center}
(car $|x_0|_1 = |x_0|^a_2$), ce qui nous permet d'affirmer que :
\begin{center}
$\forall x \in K, |x|_1 = |x|^a_2$
\end{center}
\begin{flushright}
$\blacksquare$
\end{flushright}

\paragraph*{Corollaire :}
\subparagraph*{} Deux valeurs absolues $|.|_p$ et $|.|_l$ (respectivement $p$-adique et $l$-adique) sont équivalentes si, et seulement si, $p = l$.

\paragraph*{Preuve :}
\subparagraph*{} La réciproque est triviale. Pour l'implication directe, il suffit de considérer la suite $(p_n)_{n>0}$. Elle converge vers 0 pour $|.|_p$ car : 
\begin{center}
 $|p^n|_p = p^{-n} \rightarrow 0$ quand $n \rightarrow +\infty$
 \end{center}
 et si $p \neq l$, elle ne converge pas vers 0 pour $|.|_l$ car :
\begin{center}
 $|p^n|_l = 1 \nrightarrow 0$
 \end{center} 
\begin{flushright}
$\blacksquare$
\end{flushright} 
 
\paragraph*{Proposition :}
\subparagraph*{} Soient $|.|_1$, $|.|_2$, ... et $|.|_n$ des valeurs absolues non triviales et non équivalentes deux à deux sur un corps $K$, alors : \begin{center}
Il existe $x \in K$ tel que $|x|_1> 1$ et $(\forall i \geqslant 2)$ $|x|_i< 1$
\end{center} 
\paragraph*{Preuve :}
\subparagraph*{} Par récurrence sur $n$ :\\
Pour $n = 2$ : $|.|_1$ n'est pas équivalente à $|.|_2$
\begin{center}
donc $(\exists (a, b) \in K^2)$ :
$\left\{ \begin{array}{rcl}
\mid b \mid_1 < 1 & \mbox{et}
& \mid b \mid_2 \geqslant 1 \\ \mid a \mid_1 \geqslant 1 & \mbox{et} & \mid a \mid_2 < 1 \\
\end{array}\right.$
\end{center}
Il suffit donc de prendre $x = \dfrac{b}{a}$.\\
Pour $n > 2$, on choisit $a \in K$ tel que $| a |_1 > 1$ et $(\forall i \in \{2, 3, ..., n - 1 \})$ $| a |_i< 1$ et $b \in K$ tel que $| b |_1> 1$ et $| b |_n < 1$. On a alors :
\begin{itemize}
\item[$\ast$]$1^{er}$ cas : Si $ | a |_n< 1 $, alors $x = a$ convient.
\item[$\ast$]$2^{\text{ém}}$ cas : Si $|a|_n= 1$, alors

$a^mb \xrightarrow[m\to +\infty]{} $
$\begin{cases}
+\infty & \text{pour } \mid . \mid_1 \\
b & \text{pour } \mid . \mid_n \\
 0 & \text{pour } \mid. \mid_i \text{; } 2 \leqslant i \leqslant n - 1
\end{cases}$

Il suffit de prendre $x = a^mb$ (pour $m$ assez grand).
\item[$\ast$]$3^{\text{ém}}$ cas : Si $| a |_n> 1$, alors 

$\dfrac{a^mb}{1 + a^m} \xrightarrow[m\to +\infty]{}$ 
$\begin{cases}
b & \text{pour}  \mid . \mid_n \text{et} |.|_1 \\
0 & \text{pour}  \mid . \mid_i \text{; } 2 \leqslant i \leqslant n - 1   \\
\end{cases}$

\end{itemize}
Ainsi $\dfrac{a^mb}{1 + a^m} $ convient pour $m$ assez grand. 
\begin{flushright}
$\blacksquare$
\end{flushright}

\paragraph*{Théorème (\textbf{Approximation faible d'\textsc{Artin-Whaples}}) : } 
\subparagraph*{} Soit $K$ un corps muni des
valeurs absolues $|.|_1$, $|.|_2$, ... et $|.|_n$ non triviales et non équivalentes deux à deux , alors : \begin{center}
$\forall (x_1, x_2,..., x_n) \in K^n$ et $\forall \epsilon > 0$, $\exists x \in K$ : $(\forall i \in \{1, 3, ..., n\})$ $| x - x_i |_i < \epsilon$
\end{center}
(C'est à dire que la diagonale de $K^n$ est partout dense dans $K^n$).

\paragraph*{Preuve :}
\subparagraph*{} $(\forall i \in \{1, 3, ..., n\}) (\exists a_i \in K) :| a_i |_i > 1$ et $(\forall j \in \{1, 3, ..., n \} \backslash \{i\})$ $| a_i |_j< 1$.\\
Soit\begin{center}
 $b_m$ = $\sum \limits_{\underset{}{j=1}}^{n} \dfrac{a_j^m}{1 + a_j^m}x_j$
 $\xrightarrow[m\to +\infty]{} x_i  $  pour $|.|_i$ 
\end{center}
et $x = b_m$ convient pour $m$ assez grand.
\begin{flushright}
$\blacksquare$
\end{flushright}
\paragraph*{Lemme :}
\subparagraph*{} Une valeur absolue sur un corps $K$ est ultramétrique si et seulement si elle est bornée sur l'anneau premier.
\paragraph*{Démonstration :}
\subparagraph*{} La nécessité est triviale.\\
Inversement :\\
Pour tout $(x, y) \in K^2$ et pour tout $n \in \mathbb{N}^*$.\\
On peut écrire (formule du binôme de Newton et inégalité triangulaire):
\begin{center}
$|(x + y)^n|\leqslant \sum\limits_{k =0}^{n} 
| C_n^kx^ky^{n-k} | \leqslant  (n + 1) C (max (| x |, | y|))^n$
\end{center}
où $C$ est telle que $\forall m \in \mathbb{Z}$, $| m.1_K |\leqslant C$, d'où :
\begin{center}
$| x + y |\leqslant (n + 1)^{\dfrac{1}{n}}C^{\dfrac{1}{n}} max (| x |, | y |)$.
\end{center}
Il suffit de faire tendre $n$ vers $+\infty$ pour obtenir le résultat.
\begin{flushright}
$\blacksquare$
\end{flushright}

\paragraph*{Premier théorème d’\textsc{Ostrowski} :}
\subparagraph*{} Sur $\mathbb{Q}$, il n'existe à l'équivalence près qu'une seule valeur absolue non ultramétrique. Elles sont toutes de la forme : $| . |=| . |^{\rho}_{\infty}$ où $| . |_\infty$ est la valeur absolue usuelle et $\rho \in ]0, 1]$.\\
\paragraph*{Preuve :}
\subparagraph*{} Soit $| . |$ une valeur absolue non ultramétrique sur $\mathbb{Q}$, alors pour tout entier naturel $n$ on a :
$| n | \leqslant | 1 | + | 1 | +...+ | 1 | \leqslant n$.\\
D'où $(\forall n \in \mathbb{Z})$ $| n |\leqslant| n |_{\infty}$.\\
Soient $m$ et $n$ deux entiers supérieurs strictement à 1, alors il existe\\
$(a_0, a_1, ..., a_s) \in \mathbb{N}^{s+1} $ tel que $m = a_0 + a_1n + ... + a_sn^s$ (1) \\
Et $(\forall i \in \{0, 1, ..., n\})$ $0 \leqslant a_i < n$ et $a_s \neq 0$ (et si $m < n$, $m = a_0$) (l'écriture de $m$ dans la base $n$).\\
On a : $m \geqslant n^s$, alors :
\begin{center}
 $s \leqslant \dfrac{ln(m)}{ln(n)}$ (2) et 
$\mid m \mid \leqslant \mid a_0 \mid + \mid a_1 \mid\mid n \mid +...+ \mid a_s \mid . (\mid n \mid)^s$\\
\end{center}
donc :
$\mid m \mid \leqslant n(1 + s) max(1, (\mid n \mid)^s)$ (d'après (1) et le fait que $0 \leqslant a_i < n$)\\
alors :
\begin{center}
$|m|\leqslant n \left(1+\dfrac{\ln (m)}{\ln (n)} \right) \max \left(1, (|n|)^{\dfrac{\ln(m)}{\ln (n)}}\right)$ d'après (2)
\end{center}
Si l'on substitue dans cette formule $m^t$, $t \in \mathbb{N}^*$, à $m$ et l'on prend la racine $t^{ieme}$ on a :
$|m| \leqslant \sqrt[t]{n \left(1+t\dfrac{\ln (m)}{\ln (n)}\right)} \max \left(1,(|n|)^{\dfrac{\ln (m)}{\ln (n)}}\right)$
Comme $(\forall \alpha \in [0, +\infty[)$, $\lim\limits_{n \to +\infty}\sqrt[n]{(1 + \alpha n)} = 1$ ; on a :
\begin{center}
$|m| \leqslant \max \left(1,(|n|)^{\dfrac{\ln (m)}{\ln (n)}} \right)$ (3)
\end{center}
$\mid . \mid$ étant non ultramétrique, elle n'est pas bornée sur $\mathbb{Z}$ (d'après le lemme précédent), il existe $n_0 \in \mathbb{N}$ tel que $\mid n_0 \mid > 1$ ; pour $m = n_0$, on a :
\begin{center}
 $(\forall n \in \mathbb{N}^* \backslash \{1\}) \mid n \mid > 1$ 
\end{center}
de sorte que (3) devient :
\begin{center}
$\mid m \mid \leqslant (\mid n \mid)^{\dfrac{ln(m)}{ln(n)}}$
\end{center}
c'est à dire que $\dfrac{ln (\mid m \mid)}{ln (m)} \leqslant \dfrac{ln (\mid n \mid)}{ln (n)}$.
Par symétrie, ces deux membres ont la même valeur $\rho : 0 < \rho \leqslant 1$. D'où $\ln (|m|) = \rho \ln (m)$ ou encore $| m| = m^\rho$;\\
On obtient : $| - m| = |m| = m^{\rho}$ 
et 
$|\dfrac{m}{n}| = \left( \dfrac{m}{n} \right)^\rho$
de sorte que $|.|$ est équivalente à $|.|^{\rho}_{\infty}$ sur $\mathbb{Q}$ $(0 < \rho \leqslant 1)$.
\begin{flushright}
$\blacksquare$
\end{flushright}

\subsection*{1.4. Caractérisations des valeurs absolues}
\addcontentsline{toc}{subsection}{1.4. Caractérisations des valeurs absolues}
\subparagraph*{} Pour une application $f$ de $K$ dans $\mathbb{R}^+$, et un nombre réel $A > 0$, nous noterons $(U_A)$ la relation :
\begin{center}
 $\forall x,y \in k$ : $ f(x+y)\leqslant A. \sup(f(x), f(y))$ 
\end{center}
\subparagraph*{} Nous noterons $\vartheta(K)$ l'ensemble des applications $f$ de $K$ dans $\mathbb{R}^+$ vérifiant les deux axiomes $(VA_1)$ et $(VA_2)$ de la valeur absolue et pour lesquelles il existe un $A > 0$ (dépendant de $f$) tel que $(U_A)$ soit vraie. 
\subparagraph*{} On remarquera que si  $f \in \vartheta(K)$, on a, en faisant $x = 1$, $y = 0$ dans $(U_A)$, $1 = f(1+0)\leqslant A.sup (f(1),f(0)) = A$, donc $1 \leqslant A$ .
\paragraph*{Proposition 1 :}
\subparagraph*{} Pour qu'une application $f$ de $K$ dans $\mathbb{R}^+$ vérifiant $(VA_1)$ et $(VA_2)$ appartienne à $\vartheta(K)$, il faut et il suffit que :
\begin{center}
 $\forall x\in \{ x \in  K$; $ f(x)\leqslant 1 \}$ : $f(1 + x)$ soit borné.
 \end{center} 
\paragraph*{Démonstration :}
\subparagraph*{$\Rightarrow$ :} Soit $x \in K$, si $f$ vérifie $(U_A)$, On a :
\begin{center} 
$f(1+x)\leqslant A.sup (f(1),f(x))$ 
\end{center}
et si $f(x)\leqslant 1$, on aura :\begin{center}
$sup (f(1),f(x)) = f(1) =1$
\end{center}
Donc :
\begin{center}
$f(1+x) \leqslant A$
\end{center}
Donc :
 \begin{center}
$\forall x\in \{ x \in  K$; $ f(x)\leqslant l \}$ : $f(1 + x)$ est borné
\end{center}
\subparagraph*{$\Leftarrow$ :} Inversement, supposons que $f(x + 1) \leqslant A$ pour les $x \in K$ tels que $f(x) \leqslant 1$ (ce qui entraîne $A \geqslant f(1) = 1$); alors :\\
Si $x = 0$ ou $y = 0$ :
\begin{center}
la condition $(U_A)$ est vérifiée;
\end{center}
Si au contraire $x\neq 0$ et $y \neq 0$, on peut par exemple supposer $f(y) \leqslant f(x)$, donc, d'après $(VA_2)$ : \begin{center}
$f(yx^{-1}) \leqslant 1$
\end{center}
et par suite $f(1 + yx^{-1})$, ce qui donne, en vertu de $(VA_2)$ :
\begin{center}
$f(x + y)f(x)^{-1} \leqslant A$
\end{center}
D'où 
\begin{center}
$f(x + y) \leqslant A.f(x) \leqslant A.\sup(f(x), f(y))$. 
\end{center}
\begin{flushright}
$\blacksquare$\\
\end{flushright}

\subparagraph*{} Si $f$ est une valeur absolue sur $K$, on a $f(n.1)\leqslant n$ par récurrence sur l'entier $n>0$ à partir de $(VA_3)$; réciproquement : 
\paragraph*{Proposition 2 :}
\subparagraph*{} Soit $f$ une application de $K$ dans $\mathbb{R}^+$ appartenant à $\vartheta(K)$; s'il existe $C>0$ tel que $f(n.1)\leqslant C.n$ pour tout entier $n > 0$, $f$ est une valeur absolue sur $K$.\\

\paragraph*{Démonstration : }
\subparagraph*{} Par récurrence sur $r > 0$, on déduit de ($U_A$) la relation 
\begin{center}
$f(x_1+x_2+...+x_{2^r}) \leqslant A^r \sup\limits_{i\in \{1,...,2^r\}} f(x_i)$ (1)
\end{center}
pour toute famille ($x_i$) de $2^r$ éléments de $K$. Posons $n = 2^{r}-1$; pour tout $x \in K$, on déduit de (1) 
\begin{center}
$(f(1 + x))^n = f((1 + x)^n) = f\left( \sum \limits_{\underset{}{i=0}}^n \binom{n}{i} x^i \right) \leqslant A^r \sup \left( f \left( \binom{n}{i} \right)\left( f(x) \right)^i \right)$ 
\end{center}
\begin{center}
$\leqslant CA^r  \sum \limits_{\underset{}{i=0}}^n \binom{n}{i} \left( f(x) \right)^i = CA^r(1+f(x))^n$
\end{center}
car $f\left( \binom{n}{i} \right) \leqslant C \binom{n}{i}   $ on a donc

\begin{center}
$f(1+x) \leqslant C^{1/n}A^{r/n} (1+f(x))$.
\end{center}

Faisons tendre $n$ vers $+\infty$, il vient $f(1 + x) \leqslant 1 + f(x)$ pour tout $x \in K $; appliquons cette inégalité en remplaçant $x$ par $xy^{-1}$ (pour $y \neq 0$) et en tenant compte de ($VA_2$), on obtient la relation ($VA_3$), ce qui prouve la proposition. 
\begin{flushright}
$\blacksquare$
\end{flushright}
\paragraph*{Corollaire 1 :}
\subparagraph*{} Pour qu'une application $f$ de $K$ dans $\mathbb{R}^{+}$ soit une valeur absolue, il faut et il suffit qu'elle vérifie les conditions ($VA_1$), ($VA_2$) et ($U_2$). 

\paragraph*{Démonstration :}
\subparagraph*{} C'est nécessaire, car ($VA_3$) entraîne 
\begin{center}
 $f(x + y)\leqslant f(x) + f(y) \leqslant 2 sup(f(x), f(y))$.
 \end{center} 
\subparagraph*{} Inversement, supposons que $f$ vérifie ($VA_1$), ($VA_2$) et ($U_2$); pour tout entier $n > 0$, soit $r$ le plus petit entier tel que $2^r > n$ ; si dans (1) on remplace $A$ par 2, les $x_i$ d'indice $i \leqslant n$ par 1 et les $x_i$ d'indice $i > re$ par 0, on obtient $f(n.1) \leqslant 2^r \leqslant 2n$; on peut alors appliquer la proposition 2 avec $C = 2$, donc $f$ est une valeur absolue. 
\begin{flushright}
$\blacksquare$
\end{flushright}
\paragraph*{Corollaire 2 :}
\subparagraph*{} Pour qu'une application $f$ de $K$ dans $\mathbb{R}^+$ appartienne à $\vartheta(K)$, il faut et il suffit qu'elle soit de la forme $g^t$, où $t > 0$ et $g$ 
est une valeur absolue, sur $K$. 
\subparagraph*{} En effet, dire que $f$ vérifie ($U_A$) équivaut à dire que $f^s$ vérifie ($U_{A^S}$), comme il existe $s>0$ tel que $A^s \leqslant 2 $, le corolaire 1 montre que pour une telle valeur de $s$, $f$ est une valeur absolue.
\begin{flushright}
$\blacksquare$
\end{flushright}
\paragraph*{Remarque :}
On peut donner une autre définition d'une valeur absolue ultramétrique on se basant sur les derniers résultats :
\subparagraph*{} On dit qu'une application $f$ de $K$ dans $\mathbb{R}^+$ est une valeur absolue ultramétrique si elle vérifie les conditions ($VA_1$), ($VA_2$) et ($U_1$) (ce qui entraîne évidemment que $f$ est une valeur absolue). 

\paragraph*{Proposition 3 :}
\subparagraph*{} Soit $f$ une application de $K$ dans $\mathbb{R}^+$. Les propriétés suivantes sont équivalentes :
\begin{itemize}
\item[a) -] $f$ est une valeur absolue ultramétrique.
\item[b) -] Il existe une valuation $\upsilon$ de $K$, à valeurs dans $\mathbb{R}$ , et un nombre, réel $a$ tels que $0 < a < 1$ et $f = a^{\upsilon}$. 
\item[c) -] $f$ appartient à $\vartheta(K)$ et l'on a $f(n.1) \leqslant 1$ pour tout entier $n> 0$.
\item[d) -] Pour tout $s > 0$, $f^s$ est une valeur absolue. 
\end{itemize} 
\paragraph*{Démonstration : }
\subparagraph*{} Pour tout nombre réel $c$ tel que $0 < c < 1$, l'application $ t \mapsto c^t $ est un isomorphisme du groupe ordonné $\mathbb{R}$ (muni de l'ordre opposé à l'ordre usuel) sur le groupe ordonné $\mathbb{R}^*$ ; cela montre l'équivalence de a) et b) \footnote{On parlera de façon plus détaillée de la relation entre les valeurs absolues et les valuations dans le début de la partie suivante.}. Il est clair que a) implique c); c) entraîne d), car on déduit de c) que $(f(n.l))^s \leqslant 1 \leqslant n$ pour tout entier $n > 0$ et la proposition 2 montre que $f^s$ est une valeur absolue. Enfin d) entraîne a) : en effet, si $f^s$ est une valeur absolue, elle vérifie ($U_2$), donc $f$ vérifie ($U_{2^s}$) pour tout $s > 0$, et par suite aussi ($U_1)$ en faisant tendre $s$ vers $\infty$. 
\begin{flushright}
$\blacksquare$
\end{flushright}
\paragraph*{Corollaire 3 :}
\subparagraph*{} Si $K$ est un corps (pas nécessairement commutatif) de caractéristique $p > 0$, toute fonction de $\vartheta(K)$ est une valeur absolue ultramétrique. 

\paragraph{Démonstration :}
\subparagraph*{} En effet, tout élément $z = n.1$ ($n$ entier positif ) non nul appartient au sous-corps premier $\mathbb{F}_p$ de $K$, donc vérifie la relation $z^{p-1} = 1$, ce qui entraîne $f(z) =1$ et l'on peut appliquer la proposition 3, c). 
\begin{flushright}
$\blacksquare$
\end{flushright}
\paragraph*{Remarque :}
\subparagraph*{} Étant donné un nombre réel $c$ tel que $0< c< 1$, les formules 
\begin{center}
 $f(x) = c^{\upsilon(x)}$ , $\upsilon(x) = log_c f(x)$
 \end{center} 
établissent donc une correspondance biunivoque entre valeurs 
absolues ultramétr-iques sur $K$ et valuations de $K$ à valeurs réelles.\\ 
A la valeur absolue impropre correspond la valuation impropre. Soient $\upsilon_1(x)$, $\upsilon_2(x)$ deux valuations 
de $K$ à valeurs réelles, et $f_1$, $f_2$ les valeurs absolues  
correspondantes; pour que $\upsilon_1(x)$ et $\upsilon_2(x)$ soient équivalentes, il faut et il suffît que $f_1$ et $f_2$ le soient : en effet, dire que $\upsilon_1(x)$ et $\upsilon_2(x)$ sont équivalentes revient à dire que les relations $\upsilon_1(x) \geqslant 0$ et $\upsilon_2(x) \geqslant 0$ sont équivalentes, ou encore que les relations $f_1(x) \leqslant 1$ et $f_2(x) \leqslant 1$ sont équivalentes; il suffit donc d'appliquer les propositions vues précédemment.

\section*{2. Complétude et complétion d'un corps valué }
\addcontentsline{toc}{section}{2. Complétude et complétion d'un corps valué }
\subsection*{2.1. Construction du complété d'un corps valué}
\addcontentsline{toc}{subsection}{2.1. Construction du complété d'un corps valué}
\paragraph*{ Corps valué complet : }
\subparagraph*{} Soit $(K,\vert.\vert)$ un corps valué. On dit que $K$ est complet si l'espace métrique $(K,d)$ est complet. \\
i.e. Toute suite de Cauchy dans $(K,d)$ est convergente dans $K$. 
\paragraph*{ Complété d'un corps valué : } Soit $(K,\vert.\vert)$ un corps valué. \\
Comme mentionné dans la partie préliminaire, on peut obtenir le complété de $K$ en ajoutant les limites des suites de Cauchy de $K$, $K$ est également contenu dans son complété en considérant les suites constantes, on verra que cette construction, dont l'origine est très simple, devient plutôt compliquée quand traduite en langage mathématique. Plus formellement, on pose :
\subparagraph*{$\ast$} $C(K)=\lbrace A=(a_{n})_{n\geq1} \in K^{\mathbb{N^{*}}} /$ $\lim\limits_{\substack{n \rightarrow +\infty \\ m\rightarrow +\infty }} \vert a_{n} - a_{m} \vert = 0 \rbrace $ l'ensemble des suites de Cauchy de $K$, qui est un anneau unitaire,
\subparagraph*{$\ast$} $ M(K)=\lbrace A=(a_{n})_{n\geq1} \in K^{\mathbb{N^{*}}} /$ $\lim\limits_{\substack{n \rightarrow +\infty  }} \vert a_{n} \vert = 0 \rbrace $ qui est un idéal de $C(K)$,
\subparagraph*{$\ast$} Et $\widehat{K}=C(K)/M(K)$ l'anneau quotient.

\paragraph*{Proposition :} 
\subparagraph*{} Soit $(K,\vert.\vert)$ un corps valué. L'ensemble $\widehat{K}$ est un corps. La valeur absolue de $K$ se prolonge de façon unique sur $\widehat{K}$ et on la note également $\vert.\vert$. De plus, $(\widehat{K},\vert.\vert)$ est complet et il est unique à un isomorphisme près. \\
$\widehat{K}$ est appelé \textbf{complété de $K$}.
\paragraph*{Démonstration :}
\subparagraph*{$\ast$} On montre d'abord que $C(K)$ est un anneau unitaire et que $M(K)$ est un idéal de $C(K)$. \\
Toute suite $A=(a_{n})_{n\geq1} \in C(K)$ est \textit{bornée}. En effet,\\
Il existe $n_{1}$ entier tel que $\forall m,n \geq n_{1}$, on a : \begin{center}
$ \vert a_{n}-a_{m} \vert < 1$ ;
\end{center} 
En particulier, \begin{center}
$ \vert a_{n}-a_{n_{1}} \vert < 1$, $ \forall n>n_{1} $.
\end{center} 
D'où \begin{center}
$ \vert a_{n} \vert < 1+ \vert a_{n_{1}} \vert $,  $ \forall n>n_{1} $ et $ \alpha = \sup_{n\geq 1} \vert a_{n} \vert \leq \max(\vert a_{1} \vert,..., \vert a_{n_{1}} \vert, 1+\vert a_{n_{1}} \vert)$
\end{center}
$C(K)$ est un \textit{anneau unitaire}, d'unité $(1,1,...,1,...)$.\\
Il suffit de vérifier que $C(K)$ est un sous-anneau de l'anneau produit $K^{\mathbb{N^{*}}}$.\\
Si $A=(a_{n})_{n\geq1}$ et $B=(b_{n})_{n\geq1}$ deux éléments de $C(K)$, on a pour $ n\geq1 $ et $ m\geq 1$ :
\begin{center}
$a_{n}b_{n}-a_{m}b_{m}=(a_{n}-a_{m})b_{n}+a_{m}(b_{n}-b_{m})$
\end{center}
ainsi :
\begin{center}
$\vert a_{n}b_{n}-a_{m}b_{m} \vert \leq \beta \vert a_{n}-a_{m}\vert + \alpha\vert b_{n}-b_{m} \vert$ 
\end{center}
Où :
\begin{center}
$\alpha=\sup_{n\geq 1} \vert a_{n} \vert$ et $\beta=\sup_{n\geq 1} \vert b_{n} \vert$.
\end{center}
On en déduit que :
\begin{center}
$\lim\limits_{\substack{n \rightarrow +\infty \\ m\rightarrow +\infty }} \vert a_{n}b_{n} - a_{m}b_{m} \vert = 0$,
\end{center}
C'est-à-dire $AB \in C(K)$.\\
Soit maintenant $A$ dans $M(K)$, puisque :\begin{center}
$ \vert a_{n}-a_{m} \vert \leq \vert a_{n} \vert + \vert a_{m} \vert $
\end{center}
On a :
\begin{center}
$\lim\limits_{\substack{n \rightarrow +\infty \\ m\rightarrow +\infty }} \vert a_{n} - a_{m} \vert = 0$
\end{center}
Ce qui implique que $M(K)\subset C(K)$. On vérifie facilement que $M(K)$ est un \textit{idéal} de $C(K)$.

\subparagraph*{$\ast$} L'anneau quotient $\widehat{K}$ est un corps.\\
Considérons $A=(a_{n})_{n\geq1} \in C(K)\setminus M(K)$. Puisque \begin{center}
$\vert \vert a_{n} \vert - \vert a_{m} \vert \vert \leq  \vert a_{n} - a_{m} \vert $
\end{center}
la suite $(|a_{n}|)_{n\geq1}$ est de Cauchy dans $ \mathbb{R^{+}} $, donc :
\begin{center}
 $\lim\limits_{\substack{n \rightarrow +\infty}} \vert a_{n}\vert =\alpha > 0$.
\end{center}
Ainsi, il existe $n_{0}\geq1$ et pour tout $n \geq n_{0}$, on a :
\begin{center}
$\vert \vert a_{n} \vert - \alpha \vert <\dfrac{\alpha}{2}$,
\end{center} 
Donc :
\begin{center}
$\dfrac{\alpha}{2}<\vert a_{n} \vert$, $\forall  n\geq n_{0} $.
\end{center}
Soit $ B=(b_{n})_{n\geq1} $, la suite définie par $b_{n}=1$ pour $1\leq n < n_{0}$ et $b_{n}=a_{n}$ pour $n \geq n_{0}$. Alors \begin{center}
$\vert b_{n} - b_{m} \vert = \vert a_{n} - a_{m} \vert$, pour $n, m \geq n_{0}$ et $B\in C(K)$.
\end{center}
De plus\begin{center}
 $ b_{n}^{-1} \in K $, $\forall n\geq 1$ 
\end{center}
Avec 
\begin{center}
$\vert b_{n}^{-1} - b_{m}^{-1} \vert = \vert b_{n}^{-1}\vert  \vert b_{m}^{-1} \vert \vert b_{n} - b_{m} \vert \leq \dfrac{4}{\alpha^{2}}\vert b_{n} - b_{m} \vert$. 
\end{center}
Il vient que 
\begin{center}
$\lim\limits_{\substack{n \rightarrow +\infty \\ m\rightarrow +\infty }} \vert b_{n}^{-1} - b_{m}^{-1} \vert = 0$ et $C=(b_{n}^{-1})_{n\geq1} \in C(K)$. 
\end{center}
D'autre part, comme $\vert a_{n} - b_{n} \vert = 0$, $\forall n \geq n_{0}$, on a :
\begin{center}
$A-B \in M(K)$
\end{center}
Posons $\widehat{a}$ la classe de $A$ modulo $M(K)$. Puisque $BC=1$, on a $\widehat{b}\widehat{c}=\widehat{bc}=\widehat{1}$. Mais $\widehat{b}=\widehat{a}$, il vient que si    $ A \in C(K) \setminus M(K)$, sa classe $\widehat{a}$ est inversible dans $\widehat{K}$. Mais la classe $\widehat{a}$ d'un élément $A$ de $C(K)$ est non nulle dans $\widehat{K}$ si et seulement si $A\notin M(K)$. On en déduit que tout élément non nul de $\widehat{K}$ est inversible et $\widehat{K}$ est un corps.      

\subparagraph*{$\ast$} Extension de la valeur absolue.\\
Soient $A=(a_{n})_{n\geq1} \in C(K)$ et $\widehat{a}$ sa classe dans $\widehat{K}$. Posons $\vert \widehat{a} \vert = \lim\limits_{\substack{n \rightarrow +\infty  }} \vert a_{n} \vert $ alors $\vert \widehat{a} \vert$ ne dépend pas de $ A\in \widehat{a}$. En effet, si $A, B \in \widehat{a}$, on a :
\begin{center}
$A-B \in M(K)$ 
\end{center}
Et comme \begin{center}
$\vert a_{n} \vert \leq \vert a_{n}-b_{n} \vert + \vert b_{n} \vert $ et $\vert b_{n} \vert \leq \vert b_{n}-a_{n} \vert + \vert a_{n} \vert $
\end{center} 
On obtient \begin{center}
$\vert \widehat{a} \vert = \lim\limits_{\substack{n \rightarrow +\infty  }} \vert a_{n} \vert = \lim\limits_{\substack{n \rightarrow +\infty  }} \vert b_{n} \vert = \vert \widehat{b} \vert $
\end{center}
On vérifie alors que 
\begin{center}
$\vert \widehat{a} \vert = \lim\limits_{\substack{n \rightarrow +\infty  }} \vert a_{n} \vert $ définit une valeur absolue sur $\widehat{K}$.
\end{center}

\begin{flushright}
$\blacksquare$
\end{flushright}

Considérons l'injection canonique $ i : K\rightarrow \widehat{K} $ définie par $i(a)=(a,...,a...)$. On a $ \vert \widehat{i(a)} \vert = \vert a \vert$. 

\subparagraph*{$\ast$} $(\widehat{K},\vert.\vert)$ est complet.\\
Soit $(\widehat{a}(n))_{n\geq1} $ une suite de Cauchy. Considérons $A(m)=(a_{m}(n))_{m\geq1} \in \widehat{a}(n)$. Pour tout $ \varepsilon>0 $, il existe $n_{\varepsilon}$ tel que pour $n, q \geq n_{\varepsilon}$, on a :
\begin{center}
$\vert \widehat{a}(n)-\widehat{a}(q) \vert = \lim\limits_{\substack{m \rightarrow +\infty  }} \vert a_{m}(n)-a_{m}(q) \vert < \varepsilon$.
\end{center}
Mais il existe $m_{\varepsilon}$ tel que $\forall m\geq m_{\varepsilon}$, on a :
\begin{center}
$\vert{a}_{m}(n)-{a}_{m}(q)\vert < \varepsilon$, pour tous $n, q \geq n_{\varepsilon}$. 
\end{center}
Ainsi : 
\begin{center}
$\forall m\geq \max(m_{\varepsilon}, n_{\varepsilon})$, $\forall n\geq n_{\varepsilon}$, on a $\vert{a}_{m}(n)-{a}_{m}(m)\vert < \varepsilon$.
\end{center} 
La suite $A=(a_{m}(m))_{m\geq1}$ appartient à $ C(K)$. Si $\widehat{a}$ est la classe de $A$, on a :
\begin{center}
 $\lim\limits_{\substack{n \rightarrow +\infty  }} \vert \widehat{a}(n)-\widehat{a} \vert = \lim\limits_{\substack{n \rightarrow +\infty \\ m\rightarrow +\infty }} \vert a_{m}(n) - a_{m}(m) \vert \leq \varepsilon $, pour tout $\varepsilon > 0$
 \end{center}
Donc \begin{center}
$\lim\limits_{\substack{n \rightarrow +\infty  }} \vert \widehat{a}(n)-\widehat{a} \vert = 0$ 
\end{center}
C'est-à-dire
\begin{center}
 $ \widehat{a}=\lim\limits_{\substack{n \rightarrow +\infty  }} \widehat{a}(n) $ et $\widehat{K}$ est complet.
\end{center} 
\begin{flushright}
$\blacksquare$
\end{flushright}

\paragraph*{Proposition :} 
\subparagraph*{} $K$ est dense dans $\widehat{K}$ qui est unique à un isomorphisme près.
 
\paragraph*{Démonstration :} 
\subparagraph*{$\ast$} Soit $\widehat{a} \in \widehat{K}$ la classe de $A=(a_{n})_{n\geq1} \in C(K)$, pour tout $ \varepsilon > 0 $, il existe $n_{\varepsilon}$ tel que pour $m, n \geq n_{\varepsilon}$, on a :
\begin{center}
$\vert a_{m}-a_{n} \vert < \varepsilon $
\end{center} 
Ainsi :
\begin{center}
$\vert \widehat{a}-\widehat{i(a_{n})} \vert = \lim\limits_{\substack{m \rightarrow +\infty  }} \vert a_{m}-a_{n} \vert \leq \varepsilon $.
\end{center} 
Il vient que :
\begin{center}
$\lim\limits_{\substack{n \rightarrow +\infty  }} \vert \widehat{a}-i(\widehat{a_{n}}) \vert = 0 $ et $K$ est dense dans $\widehat{K}$.
\end{center}
\subparagraph*{$\ast$} Soient $(L,\vert.\vert_{1})$ un corps valué complet, $ j : K\rightarrow L $ un morphisme isométrique de corps. Alors $j$ se prolonge de façon unique en un morphisme     isométrique :
\begin{center}
$ \widehat{j} : \widehat{K}\rightarrow L $,
\end{center} 
En posant pour $ \widehat{a} \in \widehat{K}$, 
\begin{center}
$\widehat{j}(\widehat{a})=\lim\limits_{\substack{n \rightarrow +\infty  }} j(a_{n})$.
\end{center}
Si de plus $j(K)$ est dense dans $L$, on a :
\begin{center}
$\widehat{j}(\widehat{K)}=L$
\end{center} 
Dans ces conditions $\widehat{K}$ et $L$ sont isométriquement isomorphes et $\widehat{K}$ est unique à un isomorphisme isométrique près.

\begin{flushright}
$\blacksquare$
\end{flushright}

\subsection*{2.2. Exemples }
\addcontentsline{toc}{subsection}{2.2. Exemples }

\paragraph*{$\mathbb{Q}$ muni de la valeur absolue usuelle se complète en $\mathbb{R}$ :}
\subparagraph*{} Le corps des nombres réels peut être construit de plusieurs façons équivalentes dont deux sont les plus importantes : la construction par les coupures de Dedekind, imaginée par Richard \textsc{Dedekind}, et la complétion de $\mathbb{Q}$ pour la valeur absolue usuelle. Ce qui implique, entre autres et d'après ce qui précède, que $\mathbb{R}$ est complet et que $\mathbb{Q}$ est dense dans $\mathbb{R}$.

\paragraph*{$\mathbb{Q}$ muni de la valeur absolue $p$-adique se complète en $\mathbb{Q}_{p}$ :}
\subparagraph*{} On peut munir $\mathbb{Q}$, comme vu précédemment, de la valeur absolue $p$-adique, définie à partir de la valuation $p$-adique; compléter ce corps entraînera la construction du \textbf{corps des nombres $p$-adiques}, qu'on notera $\mathbb{Q}_{p}$ et qu'on utilisera tout au long de notre projet. Cette complétion a été l'origine de l'émergence de toute une branche de mathématiques : \textbf{l'analyse $p$-adique.}
\subparagraph*{Définitions : } 
\subparagraph*{ $\ast$ Anneau des entiers.} Soit $(K,\vert .\vert)$ un corps valué et $B(0,1)$ la boule unité fermée de centre $0$ et de rayon $1$ définie comme suit :\\
 $B(0,1)=\lbrace x \in K : \vert x\vert \leq 1 \rbrace$.\\
$B(0,1)$ est un sous-anneau de $K$ qui contient l'unité, il est appelé \textbf{l'anneau des entiers de $K$}. 
\subparagraph*{ $\ast$}
$\mathbb{Q}_{p} = (\mathbb{Q},\vert . \vert_p)$ $\widehat{} $  est appelé le corps des nombres $p$-adiques. Son anneau des entiers $\mathbb{Z}_{p} = \lbrace a \in \mathbb{Q}_{p} / \vert a \vert_{p}\leq 1 \rbrace$ est l'anneau des entiers $p$-adiques.      
\subparagraph*{Théorème : Développement de \textsc{Hensel} }
\subparagraph*{$(i)$} Tout élément $a$ de $\mathbb{Q}_p$ admet un développement unique sous forme de série convergente dans $\mathbb{Q}_p$ :
\begin{center}
$a = \sum\limits_{n \geqslant j_0} a_np^n$
\end{center}
où $ j_{0}=\upsilon_{p}(a) $ et $ 0\leqslant a_n \leqslant p-1$, pour tout entier $n\geqslant j_0$.
\subparagraph*{$(ii)$} L'anneau des entiers $p$-adiques $\mathbb{Z}_{p}$ est égal à l'ensemble des séries de la forme 
$a = \sum\limits_{n\geqslant 0} a_np^n$ et $\mathbb{N}$ est dense dans $\mathbb{Z}_p$.

\subsection*{2.3. Résultats et conséquences }
\addcontentsline{toc}{subsection}{2.3. Résultats et conséquences }
\paragraph*{Proposition 1 : }
\subparagraph*{} Soit $(K, \vert . \vert)$ un corps valué non archimédien, alors $ \widehat{K} $ muni du prolongement de $ \vert . \vert $ est aussi non archimédien.
\paragraph*{Démonstration :} 
\subparagraph*{} Soit $(\widehat{x}, \widehat{y}) \in \widehat{K}^{2} $, alors :
\begin{center}
 $\widehat{x}=\lim\limits_{\substack{n \rightarrow +\infty }} x_{n} $ et $\widehat{y}= \lim\limits_{\substack{n \rightarrow +\infty }} y_{n} $ où $ (x_{n})_{n} $ et $ (y_{n})_{n} \in K$;
\end{center}
D'où :
\begin{center}
$ \vert \widehat{x}+\widehat{y} \vert = \lim\limits_{\substack{n \rightarrow +\infty }} \vert x_{n}+y_{n} \vert  $
\end{center}
Et :
\begin{center}
$\lim\limits_{\substack{n \rightarrow +\infty }} \vert x_{n}+y_{n} \vert  \leq \lim\limits_{\substack{n \rightarrow +\infty }}\max(\vert x_{n}\vert , \vert y_{n} \vert)$
\end{center}
Or :
\begin{center}
$ \lim\limits_{\substack{n \rightarrow +\infty }}\max(\vert x_{n}\vert , \vert y_{n} \vert) \leq \max(\lim\limits_{\substack{n \rightarrow +\infty }} \vert x_{n}\vert , \lim\limits_{\substack{n \rightarrow +\infty }} \vert y_{n} \vert)$ 
\end{center}
Et on a :
\begin{center}
$ \max(\lim\limits_{\substack{n \rightarrow +\infty }} \vert x_{n}\vert , \lim\limits_{\substack{n \rightarrow +\infty }} \vert y_{n} \vert) = \max(\vert\widehat{x}\vert , \vert\widehat{y}\vert) $
\end{center}
D'où :
\begin{center}
$ \vert \widehat{x}+\widehat{y} \vert \max(\vert\widehat{x}\vert , \vert\widehat{y}\vert) $
\end{center}
\paragraph*{Proposition 2 : } 
\subparagraph*{} Pour tout $(\widehat{x}, \widehat{y})\in \widehat{K}^{2}$ tel que $\widehat{x}\neq \widehat{y}$, les intersections avec $K$ des voisinages ouverts assez petits de $\widehat{x}$ et $\widehat{y}$ sont contenus dans deux boules ouvertes disjointes de $K$.
\paragraph*{Conséquences :}       
\subparagraph*{$(i)$} Un corps valué est complet si, et seulement, toute chaîne de boules emboîtées dont le rayon devient arbitrairement petit a une intersection non vide.              
\subparagraph*{$(ii)$} Le complété d'un corps valué s'obtient en comblant ses trous de diamètre nul.         
                                                                                                                  \paragraph*{Démonstration :}
                                                                                                                   \subparagraph*{} $\Longrightarrow )$ Soit $ (K, \vert . \vert )$ un corps valué complet et $(B_{n})_{n\in \mathbb{N}}$ une suite strictement décroissante de boules $ ( \forall n\in \mathbb{N}, B_{n+1} \subsetneq B_{n} ) $ et $\lim\limits_{\substack{n \rightarrow +\infty }} r_{n}=0 $ ;                                                                                                        \begin{center}                                                                                                         $\forall n\in \mathbb{N}$, $ B_{n}= B_{O}(x_{n}, r_{n})$.                                                                                                         \end{center}                                                                                                          Montrons que l'intersection de cette chaîne de boules est non vide.\\                                                                                                           La suite $(x_{n})_{n\in \mathbb{N}}$ est de Cauchy (pourquoi?) et donc convergente. Si on pose :
                                                                                                                  \begin{center}                                                                                                                   $x=\lim\limits_{\substack{n \rightarrow +\infty }} x_{n}$ alors $\lbrace x\rbrace=\bigcap _{n \in \mathbb{N}} B_{n}$.                                                                                                                   \end{center}
En effet, $x\in\bigcap _{n \in \mathbb{N}} B_{n}$ résulte de la décroissance de $(r_{n})_{n\in \mathbb{N}}$ vers 0 et du fait que $x=\lim\limits_{\substack{n \rightarrow +\infty }} x_{n}$. L'unicité se déduit de l'inégalité $d(x,y)\leq r_{n}$ pour tout $n$ que réaliserait tout élément $y$ de $\bigcap _{n \in \mathbb{N}} B_{n}$ et qui implique que $x=y$.\\
La réciproque est triviale.
\begin{flushright}
$\blacksquare$
\end{flushright}
                                                                                                                     
\section*{3. Propriétés topologiques des Corps valués }
\addcontentsline{toc}{section}{3. Propriétés topologiques des Corps valués }
\subsection*{3.1. Topologie générale (ou topologie de Hausdorff)}
\addcontentsline{toc}{subsection}{3.1. Topologie générale (ou topologie de Hausdorff)}
\subparagraph*{} On a défini dans ce qui précède la notion de distance induite par une valeur absolue sur un corps ainsi que la notion de boule dans un corps topologique. On s'intéresse dans cette sous-partie aux propriétés topologiques des corps valués ultra-métriques (ou non-archimédiens) considérés comme des espaces métriques non-archimédiens (n.a).
\subparagraph*{} Dans toute la suite $K$ est considéré comme étant un espace métrique muni de la métrique $d$. Si la valuation est triviale alors la topologie associée à $K$ est discrète.\\
Un tel corps topologique est un corps $K$ doté d'une topologie dite de Hausdorff qui rend l'addition, la multiplication et la division  $(x\rightarrow x^{-1},x \neq 0)$ continues.
\subparagraph*{} Pour tout point $a$ de $K$, la famille $(B_O(a, r))_{r>0}$ est un système fondamental de voisinages de $a$. On peut considérer les figures principales des espaces
métriques, triangles, boules, sphères, en particulier on a les propriétés suivantes :
\paragraph*{Proposition (Principe du triangle isocèle) :}

\subparagraph*{} Pour tous $x, y$ et $z$ dans $K$, si $d (x, y) \neq d (y, z)$, alors :
\begin{center}
 $d (x, z) = max (d (x, y) , d (y, z))$.
\end{center}

\paragraph*{Preuve :}
\subparagraph*{} Soit $(x, y, z) \in K^3$ tel que $d (x, y) \neq d (y, z)$, alors (par exemple) :
\begin{center}
 $d (x, y) < d (y, z)$
\end{center}
D'où :
\begin{center}
$max(d (x, y) , d (y, z)) = d (y, z)$
\end{center}
Or :
\begin{center}
 $d (x, z) \leqslant max (d (x, y) , d (y, z))$
\end{center}
Il en résulte :
\begin{center}
$d (x, z) \leqslant d (y, z)$
\end{center}
Et :
\begin{center}
 $d (y, z) \leqslant max (d (y, x) , d (x, z))$ $\leqslant d (x, z)$
 \end{center}
Car :
\begin{center}
$max (d (y, x) , d (x, z)) = d (x, z)$
\end{center}
Sinon $max(d (y, x) , d (x, z)) = d (y, x)$ entraîne que : \begin{center}
$d (y, z) \leqslant d (y, x)$
\end{center} 
Ce qui contredit l'hypothèse.
Donc 
\begin{center}
$d (x, z) = d (y, z)$
\end{center} 

\begin{flushright}
$\blacksquare$
\end{flushright}
\paragraph*{Remarques :}

\subparagraph*{$\ast$} Géométriquement, dans le plan non-archimédien, tout triangle est isocèle et la base en est le plus petit côté.
\subparagraph*{$\ast$} Pour tous $x$ et $y$ dans $K$, si $| x | \neq | y |$, alors $| x+y |= max (| x |, | y |)$
\paragraph*{Propositions :}
\begin{itemize}
\item[$i.$] Tout point d'une boule, ouverte ou fermée, en est un centre .
\item[$ii.$] Pour tous $a \in K$ et $r \in ]0, +\infty[ $ ; les boules ouvertes $B_o(a, r)$ et fermées $B_f(a, r)$ sont des ouverts et fermés (O.F.) à la fois dans $(K,d)$.
\item[$iii.$] Toute sphère est un O.F.
\end{itemize}
\paragraph*{Preuve :}
\subparagraph*{} $i.$ Soit $a \in K$ et $r \in ]0, +\infty[$ ; alors, pour tout $x \in B_o(a, r)$ on a :
\begin{center}
$\forall y \in K$ : $y \in B_o(x, r) \Rightarrow d(x, y) < r$
\end{center}
L'inégalité ultra-métrique nous donne :
\begin{center}
$d (a, y) \leqslant max (d (a, x) , d (y, x))$
\end{center}
Alors :
\begin{center}
$d (a, y) < r$
\end{center}
Donc :
\begin{center}
$y \in B_o (a, r)$
\end{center}
d'où :
\begin{center}
$B_o (x, r) \subset B_o(a, r)$
\end{center}
Par symétrie, on a aussi $B_o(a, r) \subset B_o(x, r)$, et par suite :
\begin{center}
$B_o(a, r) = B_o(x, r)$
\end{center}
On fait de même pour avoir :
\begin{center}
$\forall x \in B_f(a, r)$ : $B_f(a, r) = B_f(x, r)$
\end{center}
\subparagraph*{} $ii.$ Soit $a \in K$ et $r \in ]0, +\infty [$  alors :\\
$B_f(a, r)$ est un fermé de $(K, d)$, montrons que c'est un ouvert de $(K, d)$ .\\
Soit $x \in B (a, r)$ , alors :
\begin{center}
$B_f(a, r) = B_f(x, r)$
\end{center}
Or
$B_o(x, r) \subset B_f(x, r)$ , d'où :
\begin{center}
$B_o(x, r) \subset B_f(a, r)$
\end{center}
Par la même démarche, on montre que $B_o(a, r)$ est un fermé de $(K, d)$ .\\
\subparagraph*{} $iii.$ On a $S (a, r) = \{x \in K : d (x, a) = r \} = B_f(a, r) \cap (B_o(a, r))^c$.
\begin{flushright}
$\blacksquare$
\end{flushright}
\paragraph*{Proposition :}
\subparagraph*{} Une famille finie de boules ouvertes (resp. fermées) non disjointes est une famille de cercles, notamment son intersection et sa réunion sont des boules ouvertes (resp. fermées).
\paragraph*{Preuve :}
\subparagraph*{} Soit $(B_o(x_k, r_k))_{1\leqslant k\leqslant n}$ une famille de boules ouvertes telle que :
\begin{center}
$\displaystyle{\bigcap_{1\leqslant k\leqslant n}B_o(x_k, r_k)}$ $\neq \emptyset$
\end{center}
Soit $x \in \displaystyle{\bigcap_{1\leqslant k\leqslant n}B_o(x_k, r_k)}$, alors :
\begin{center}
$\forall k \in \{0, 1, ..., n\}$ : $B_o(x_k, r_k) = B_o(x, r_k)$ 
\end{center}
Soient $ r = \displaystyle{\sup_{1\leqslant i\leqslant n}r_i}$ et $ s = \displaystyle{\inf_{1\leqslant i\leqslant n}r_i}$, Alors :
\begin{center}
$\displaystyle{\bigcap_{1\leqslant k\leqslant n}B_o(x_k, r_k)} = B_o(x, s)$ et
$\displaystyle{\bigcup_{1\leqslant k\leqslant n}B_o(x_k, r_k)} = B_o(x, r)$
\end{center}
\begin{flushright}
$\blacksquare$
\end{flushright}
\paragraph*{Conséquence :}
\subparagraph*{} Deux boules sont ou bien disjointes ou bien l'une contient l'autre.
\paragraph*{Proposition :}
\subparagraph*{} La distance d'un point, en dehors d'une boule, aux points de cette boule ne dépend pas du choix de ces points.
\paragraph*{Preuve :} 
\subparagraph*{} Soit $B (a, r)$ une boule et soient $a_1$ et $a_2$ deux points de $B (a, r)$ et $x$ un élément extérieur à $B (a, r)$ , alors :
\begin{center}
$d (a, x) \geqslant r$ et $d (a, a_1) < r$ 
\end{center}
D'où le triangle $(xaa_1)$ est isocèle de base $a_1a$, d'où :
\begin{center}
$d (a, x) = d (a_1, x )$ 
\end{center}
Et on a aussi :
\begin{center}
$d (a, x) = d (a_2, x )$
\end{center}
D'où \begin{center}
$d (a_1, x) = d (a_2, x)$
\end{center}
\begin{flushright}
$\blacksquare$
\end{flushright}
\paragraph*{Proposition :}
\subparagraph*{} La distance des points de deux boules disjoints $B_1$ et $B_2$ ne dépendent pas du choix des points, et elle est égale à : $d (B_1, B_2)$.
\paragraph*{Preuve :}
\subparagraph*{} Soient $x \in B (a, r_1)$ et $y \in B (b, r_2)$, alors :
\begin{center}
$d (x, y) = d (x, b) = d (a, b)$
\end{center}
\begin{flushright}
$\blacksquare$
\end{flushright}
\paragraph*{Théorème :}
\subparagraph*{} Soit $U$ un ouvert non vide de $K$, alors il existe une partition de $U$ en des boules de $K$. Plus précisément, étant donné une suite $(r_i)_{i>1}$ strictement décroissante de réels à termes strictement positifs, $U$ peut être couvert par des boules disjointes de la forme :
\begin{center}
$B (a, r_n)$ tel que $a \in U$ et $n \in \mathbb{N}^*$.
\end{center}
\paragraph*{Preuve :}
\subparagraph*{} Pour chaque $a \in U$, considérons la boule $B_a$ définie par :
\begin{center}
$ B_a = \left\{ \begin{array}{rcl}
B (a, r_1) & \mbox{si} & B (a, r_1) \subset U \\ B (a, r_n) & \mbox{si} & B (a, r_n) \subset U \text{ et } B (a, r_{n-1}) \nsubseteq  U \end{array}\right.$
\end{center}

Et on a : pour tout $(a, b) \in U^2$, si $B_a \cap B_b \neq \emptyset$, alors d'après la conséquence précédente (par exemple) 
\begin{center}
$B_a \subset B_b$ (1)
\end{center}
D'où si $B_a = B (a, r_n)$ et $B_b = B (b, r_m)$, on a (d'après i. de la deuxième proposition : " Tout point d'une boule, ouverte ou fermée, en est un centre") : 
\begin{center}
$B (a, r_m) = B (b, r_m)$ 
\end{center}
Mais le rayon de $B_a$ est le plus grand parmi les $(r_i)_{i>1}$ , d'où :
\begin{center}
$B (a, r_m) \subset B (a, r_n)$ 
\end{center}
Ou encore $B (b, r_m) \subset B (a, r_n)$ c'est à dire que :
\begin{center}
$B_b \subset B_a$ (2)
\end{center}
Donc d'après (1) et (2):
\begin{center}
$B_a = B_b$
\end{center}
Il s'ensuit que la collection $(B_a)_{a\in U}$ est disjointe. 
\begin{flushright}
$\blacksquare$
\end{flushright}

\subsection*{3.2. V-topologie }
\addcontentsline{toc}{subsection}{3.2. V-topologie }

Nous avons vu que les valeurs absolues ainsi que les valuations (générales) induisent canoniquement sur leur corps une topologie pour laquelle toutes les opérations de corps sont continues, avec la propriété supplémentaire suivante :
\begin{center}
$\forall x, y \in K$ et soit $V \in \vartheta(0)$ on a $
xy \in V  \Rightarrow \exists W \in \vartheta(0)$ : $x \in W$ ou $y \in W$
\end{center}
C-à-d si le produit de deux éléments est un voisinage de $0$ forcément au moins l'un de ces deux éléments est voisin de $0$.\\
De telles topologies de corps sont appelées V-topologies. Dans cette sous-partie, nous allons étudier cette topologie et on va monter que, inversement, toute V-topologie sur un corps $K$ doit être induite par une valeur absolue ou une valuation de $K$.\\

Un corps topologique doté d'une topologie de Hausdorff qui rend l'addition, la multiplication et la division $(x\rightarrow x^{-1},x \neq 0)$ continues. Compte tenu de la continuité de l'addition, une telle topologie peut être spécifiée simplement en donnant un système fondamental $T$ de voisinage de $0 \in K$. En fait, à partir de $T$ on obtient, via la translation $x \rightarrow a + x$, un système de voisinages de chaque élément $a \in K$. Par conséquent, pour un corps $K$, nous nous référerons à un système fondamental $T$ de voisinages de $0 \in K$ comme topologie sur $K$, et à la paire $(K, T)$ comme un corps topologique.\\
Rappelons que pour deux sous-ensembles R et S de K nous utilisons les notations:
\subparagraph*{} $R \pm S = \{$ $x \pm y$ $|$ $x \in R, y \in S$ $\}$
\subparagraph*{}$RS = \{$ $xy$ $|$ $x \in R, y \in S$ $\}$
\subparagraph*{}$R^{-1} = \{ $ $x^{-1}$ $|$ $x \in R$ $\}$ si $0 \notin R$

\paragraph*{Définition (V-Topologie) :}
\subparagraph*{} Une topologie (non-discrète) T sur un corps $K$ est dite V-topologie si elle vérifie les axiomes suivantes :

\begin{itemize}
\item[$VT_1$]- $\displaystyle{\bigcap_{U \in T} U}$ = $\{0 \}$; $\{0 \} \notin T$
\item[$VT_2$]- $\forall$ $U, V \in T$ $\exists W \in T : W \subseteq U \cap V$
\item[$VT_3$]- $\forall U \in T$  $\exists V \in T$ : $V - V \subseteq U$
\item[$VT_4$]- $\forall U \in T$ $\exists x; y \in K$ $\exists V \in T$ : $(x + V) (y + V) \subseteq xy + U$
\item[$VT_5$]- $\forall U \in T$ $\exists x, y \in K^*$ $\exists V \in T$ : $(x + V)^{-1} \subseteq x^{-1} + U$
\item[$VT_6$]- $\forall U \in T$ $\exists V \in T$ $\forall x, y \in \mathbb{K}$ : $xy \in V \Rightarrow x \in U$ ou $y \in U$
\end{itemize}
\paragraph*{Exemples :}
\subparagraph*{} Un exemple typique d'une V-topologie est la topologie donnée par une valeur absolue $|.|$. Les boules ouvertes $\{ x \in K ; |x| < r \}$, avec $r > 0$ sont les éléments de $T$.
\subparagraph*{} Un autre exemple de V-topologie est la topologie générée par une valuation, prenons les voisinages $U(0)$ à la place de $W \in T$, on remarque que cette topologie satisfait aux axiomes précédentes. Inversement, on a le théorème suivant dû à \textsc{Kowaslky} et \textsc{Durbaum}.
\paragraph*{Théorème :}
\subparagraph*{} Soit $K$ un corps et $T$ une topologie sur $K$. Alors $T$ est une V-topologie si et seulement s'il existe une valeur absolue archimédienne ou une valuation sur $K$ dont la topologie induite coïncide avec $T$.

\section*{4. Structures définies sur des corps valués}
\addcontentsline{toc}{section}{4. Structures définies sur des corps valués}
\subsection*{4.1. Espaces normés sur un corps valué }
\addcontentsline{toc}{subsection}{4.1. Espaces normés sur un corps valué }
\paragraph*{} On suppose, pour ce qui suit, que le lecteur est familier avec les définitions d'une norme, distance associée à une norme... \\
On renvoie également à la partie 1.2.3. (Extensions de corps) des préliminaires où on a mentionné la structure d'espace vectoriel lié à un corps.  

\paragraph*{Définitions : }
\subparagraph*{$\ast$} On appelle \textbf{espace normé} sur un corps valué non discret $K$ un espace vectoriel $E$ sur le corps $K$, muni de la structure définie par la donnée d'une norme sur $E$.    
\subparagraph*{$\ast$} On appelle \textbf{espace normable} sur $K$ un espace vectoriel sur $K$, muni d'une topologie qui peut être définie par une norme. 
\paragraph*{Exemple : } Sur un corps valué non discret $(K, \vert . \vert )$, considéré comme espace vectoriel par rapport à lui-même, la valeur absolue $\vert . \vert $ est une norme.  
\paragraph*{Théorème :} 
\subparagraph*{} Soit $K$ un corps valué complet. Toutes les normes sur un $K$-espace vectoriel $E$ de dimension finie sont équivalentes. \\
 i.e. si $\Vert . \Vert$ et $\Vert . \Vert_{1}$ sont deux normes sur $E$, alors il existe $\alpha$ et $\beta > 0$ tels que :
\begin{center}
 $\alpha \Vert x\Vert_{1} \leq \Vert x\Vert \leq \beta \Vert x\Vert_{1}$, $\forall x\in E$.
\end{center}
\paragraph*{Démonstration :}
\subparagraph*{} Soit $ (e_{i})_{1\leq i\leq m} $ une base du $K$-espace vectoriel $E$ de dimension finie $=m$.\\
On pose pour $x=\sum_{i=1}^{m} \lambda_{i}e_{i} \in E$ : 
\begin{center}
 $\Vert x \Vert_{\infty}=\max\limits_{1\leq i\leq m} \vert \lambda_{i}\vert$.
\end{center} 
$\Vert .\Vert_{\infty}$ définit une norme sur $E$. \\
On va montrer que toute norme $\Vert . \Vert : E\longrightarrow \mathbb{R}_{+}$ est équivalente à $\Vert .\Vert_{\infty}$.\\
$(E, \Vert .\Vert_{\infty})$ est d'abord isométriquement isomorphe à $(K^{m}, \vert .\vert_{\infty})$ où 
\begin{center}
$\vert (\lambda_{i})_{1\leq i\leq m}\vert_{\infty}=\max\limits{1\leq i\leq m} \vert \lambda_{i}\vert$ 
\end{center}
Puisque $K$ est complet, $(K^{m}, \vert .\vert_{\infty})$ est complet et donc $(E, \Vert .\Vert_{\infty})$ est complet.
\subparagraph*{}
On a :
\begin{center}
 $ \Vert x\Vert =\Vert\sum_{i=1}^{m} \lambda_{i}e_{i}\Vert\leq \sum_{i=1}^{m} \vert\lambda_{i}\vert \Vert e_{i}\Vert \leq \max\limits_{1\leq i\leq m} \vert \lambda_{i}\vert \sum_{i=1}^{m} \Vert e_{i}\Vert = \beta \Vert x\Vert_{\infty}$
\end{center}
où $ \beta=\sum_{i=1}^{m} \Vert e_{i}\Vert > 0 $
  
\subparagraph*{} Pour l'existence de $ \alpha $, on procède par récurrence sur la dimension $m$ de $E$. 
\subparagraph*{$\ast$} Pour $\dim E=1$, on a $E=Ke_{1}$ et tout $x \in E$ est de la forme $x=\lambda e_{1}$, avec $\Vert x\Vert = \vert \lambda\vert \Vert e_{1}\Vert =\Vert x\Vert_{\infty}\Vert e_{1}\Vert$. Ainsi, $\Vert .\Vert_{\infty}$ et $\Vert .\Vert$ sont équivalentes.   
\subparagraph*{$\ast$} On suppose que sur tout espace vectoriel de dimension strictement inférieure à $ m $, toutes les normes sont équivalentes \textit{(H.R)}.  
\subparagraph*{$\ast$} On pose pour $1\leq j\leq m$, $E_{j}= \bigoplus \limits_{\substack{i \neq j}} K.e_{i}$ ; alors :
\begin{center}
 $E_{j}$ est un sous-espace vectoriel de $E$ de dimension $m-1$.
 \end{center} 
Par hypothèse de récurrence, toutes les normes sur $E_j$ sont équivalentes. En particulier, les restrictions de $\Vert .\Vert_{\infty}$ et $\Vert .\Vert$ à $E_{j}$ sont équivalentes. Puisque $(E_j ,\Vert .\Vert_{\infty})$ est complet, l'espace normé $(E_j ,\Vert .\Vert)$ est également complet. Ainsi, $(E_j ,\Vert .\Vert)$ est fermé dans $E$. Il vient que $e_j + E_j$ est une partie fermée de $(E ,\Vert .\Vert)$.\\
Posons $\alpha_{j}= \inf \limits_{\substack{y \in e_{j}+E_{j}}} \Vert y\Vert$. Puisque $0\notin e_{j}+E_{j}$ et $e_{j}+E_{j}$ fermée dans $E$, on a :
\begin{center}
$\alpha_{j}>0$. On pose $\alpha = \min\limits_{\substack{1\leq j\leq m}} \alpha_{j} > 0$
\end{center}
Soit $x=\sum_{i=1}^{m} \lambda_{i}e_{i} \in E$ ; pour $1\leq j\leq m$ fixé, on a :
\begin{center}
$x=\lambda_j e_j + \sum\limits_{\substack{i\neq j}} \lambda_{i}e_{i} $
\end{center}
- Si $\lambda_j = 0$, alors :
\begin{center}
$\alpha \vert\lambda_j\vert = 0 \leq \Vert x\Vert $
\end{center}
- Si $\lambda_j\neq0$, alors :
\begin{center}
 $ \lambda_{j}^{-1} x= e_j + \sum\limits_{\substack{i\neq j}} \lambda_{j}^{-1}\lambda_{i}e_{i} \in e_j + E_j $ et 
$\alpha \leq \alpha_{j} \leq \Vert e_j + \sum\limits_{\substack{i\neq j}} \lambda_{j}^{-1}\lambda_{i}e_{i}\Vert = \Vert\lambda_{j}^{-1} x\Vert $
\end{center} 
D'où :
\begin{center}
$\alpha \vert\lambda_{j}\vert\leq\Vert x\Vert$
\end{center}
Donc :
\begin{center}
$\forall 1\leq j \leq m$, $\alpha \vert\lambda_{j}\vert\leq\Vert x\Vert$
\end{center} 
Ce qui implique que :
\begin{center}
$\alpha \Vert x\Vert_{\infty}\leq\Vert x\Vert$
\end{center} 
On a également :
\begin{center}
$\Vert x\Vert\leq\beta \Vert x\Vert_{\infty}$
\end{center} 
D'où les deux normes sont équivalentes, ce qui clôt la démonstration par récurrence.
\begin{flushright}
$\blacksquare$
\end{flushright}

\paragraph*{Corollaire :} 
\subparagraph*{}Soit $K$ un corps valué complet. Si $E$ est un extension finie de $K$, il existe au plus une valeur absolue sur $E$ qui prolonge celle de $K$. 

\paragraph*{Esquisse de démonstration :}
\subparagraph*{} $E$ est un sur-corps de $K$ tel que $[E:K]=\dim_{K} E < +\infty $. On considère une valeur absolue quelconque qui prolonge celle de $K$, c'est une norme sur le $K$-espace vectoriel de dimension finie $E$. En considérant une norme quelconque $\Vert .\Vert$ sur $E$, il vient que notre valeur absolue appliqué à un élément $x$ de $E$ est égale à $\lim\limits_{\substack{n\rightarrow +\infty}} \Vert x^{n}\Vert^{\frac{1}{n}}$. D'où son unicité si elle existe.        
\begin{flushright}
$\blacksquare$
\end{flushright}
\subsection*{4.2. Algèbres normées sur un corps valué }
\addcontentsline{toc}{subsection}{4.2. Algèbres normées sur un corps valué }
\paragraph*{Définitions : }
\subparagraph*{$\ast$} Étant donné une algèbre $A$ sur un corps valué commutatif non discret $K$, on dit qu'une norme $\Vert .\Vert$ sur $A$ est compatible avec la structure d'algèbre de $A$ si elle vérifie la relation : $\Vert xy\Vert \leq \Vert x\Vert .\Vert y\Vert$ pour tous $x$ et $y$ dans $A$. Une algèbre sur $K$, munie de la structure définie par une norme compatible avec sa structure d'algèbre, est appelé \textbf{algèbre normée}. 
\subparagraph*{$\ast$} Lorsqu'une algèbre $A$ est munie d'une topologie pouvant être définie par une norme, et pour laquelle $(x, y) \longmapsto xy$ est continue, on dit que l'algèbre topologique $A$ est \textbf{normable}.

\newpage 

\vspace*{\stretch{1}}
\begin{center}
\includegraphics[scale=0.8]{2.png}  
\end{center}
	\begin{center}
		{\LARGE {\huge Deuxième partie :} }\\
		\textit{   }\\
		\textit{   }\\
		{\LARGE {\Huge Corps valués locaux}}
	\end{center}
\begin{center}
\includegraphics[scale=1]{3.png}  
\end{center}
\addcontentsline{toc}{part}{Deuxième partie : Corps valués locaux}
\vspace*{\stretch{1}}
\newpage

\subparagraph*{}
\paragraph*{Introduction :}
En théorie des nombres et en mathématiques en général, la localisation est une technique importante, qui permet de simplifier de nombreux problèmes. Pour reprendre une expression de \textsc{Neukirch}, "\textit{la localisation c'est la division}", et on divise pour mieux régner. Il y a derrière cette notion en fait beaucoup plus que ça, mais c'est l'idée primaire. Nous allons, dans cette partie, essayer de présenter l'idée de localisation des corps valués qui a été à l'origine de ce qu'on appelle \textbf{corps valués locaux}.\\
Tous les anneaux considérés dans cette partie - sauf mention expresse du contraire - sont supposés être commutatifs et posséder un élément unité. Tous les homomorphismes d'anneaux sont supposés transformer l'élément unité en l'élément unité. \subparagraph*{} Tout sous-anneau d'un anneau $A$ est supposé contenir l'élément unité de $A$. Si $A$ est un anneau local, son idéal maximal sera noté $m(A)$, son corps résiduel $A/m(A)$ sera noté $\kappa(A)$ ou simplement $\kappa$, et le groupe multiplicatif des éléments inversibles de $A$ sera noté $U(A)$; on a donc $U(A) = A \setminus m(A)$. 

\section*{1. Anneau de valuation et de valuation discrète}
\addcontentsline{toc}{section}{1. Anneau de valuation et de valuation discrète}
\subsection*{1.1. Valeur absolue et valuation, quelle relation ?}
\addcontentsline{toc}{subsection}{1.1. Valeur absolue et valuation, quelle relation ?}
\subparagraph*{}
 Nous avons déjà défini dans la partie précédente la notion de "valeur absolue" ainsi que celle de "valuation", on rappelle qu'une valuation est une application d'un anneau commutatif unitaire non nul $(A, + , . )$ vers un groupe abélien totalement ordonné $(G, + ,< )$ union l'infini $\upsilon : A \rightarrow G \cup \{+\infty\}$ vérifiant les trois axiomes suivants : 
\begin{itemize}
\item[•] $\forall x \in A$ : $\upsilon(x) = + \infty$ si et seulement si $x = 0$.
\item[•] $\forall a,b \in A$ : $\upsilon(a . b) = \upsilon(a)+ \upsilon(b)$
\item[•] $\forall a,b \in A$ : $\upsilon(a + b) \geqslant inf(\upsilon(a); \upsilon(b))$ pour tout $(a, b) \in K^2$ 
\end{itemize}
\paragraph*{Remarque :}
\subparagraph*{1. } On utilise les conventions classiques $a < \infty $ et $ a+ \infty = \infty$,  $\forall a \in G$.
\subparagraph*{2. } Certains auteurs se restreignent aux
valuations sur un corps commutatif. (C'est le cas intéressant en effet).
\subparagraph*{3. } Que $A$ soit un corps ou non, $\upsilon$ est un morphisme de monoïdes de $(A^*, .)$ dans $(G, +)$.
\subparagraph*{4. } Lorsque $A$ est un corps, $\upsilon$ est donc un morphisme de groupes de $(A^*, .)$ dans $(G, +)$, si bien que $\upsilon(A^*)$ est un sous-groupe de $G$.
\subparagraph*{5. } Lorsque $A$ est un corps, on demande
parfois à $\upsilon$ d'être surjective, mais on peut toujours se ramener à cette situation en remplaçant $G$ par $\upsilon(A^*)$.
\subparagraph*{6. } Deux valuations $\upsilon$ et $\upsilon$' sur $A$ sont dites équivalentes s'il existe un isomorphisme
de demi-groupes (un ensemble muni d'une loi de composition interne associative) ordonnés:
$\lambda : \upsilon(A^*) \rightarrow \upsilon'(A^*)$ tel que $\upsilon' = \lambda \circ \upsilon$.
\subparagraph*{Valuation triviale :} Sur un corps $A$ une valuation $ \upsilon $ est dite triviale si elle est définie par :
\begin{center}
$\begin{cases}
\upsilon(x) = 0 & \text{si } x \in A^* \\
\upsilon(0) = + \infty & \text{si      } x = 0 
\end{cases}$
\end{center}
\subparagraph*{Valuation discrète :} Lorsque $G = \mathbb{Z}$ muni de l'addition, $v$ est dite \textit{valuation de Dedekind} ou \textit{valuation discrète}. Deux valuations discrètes $\upsilon$ et $\upsilon$' sur $A$ sont équivalentes si et seulement si elles sont proportionnelles, c'est-à-dire
s'il existe un rationnel $k$ non nul tel que 
\begin{center}
$\forall x \in A^* , \upsilon'(x) = k \upsilon(x)$.
\end{center}
On note $V(K)$ l'ensemble des valuations de $K$ et $Val^u(K)$
l'ensemble des valeurs absolues ultramétriques de $K$.

\subparagraph*{} Soit $\alpha \in \mathbb{R}_+^*$ avec $\alpha < 1$. Alors l'ensemble $V(K)$ des valuations de $K$ est en correspondance
biunivoque avec l'ensemble $Val^u(K)$. En effet :
\begin{itemize}
\item[1 - ] À toute valeur absolue ultramétrique on associe une valuation telle que \begin{center}
$\upsilon(x)=$
$\left\{ \begin{array}{rcl} -log(|x|) & \mbox{pour tout} & x \in K^* \\+ \infty  & \mbox{pour} & x = 0 \end{array}\right.$
\end{center}
\item[2 - ] À toute valuation $\upsilon$ sur $K$ on associe la valeur absolue ultramétrique définie par :
\begin{center}
$ |.|_{\upsilon} : K \rightarrow \mathbb{R}_+^*$ \\
$x \mapsto \alpha^{\upsilon(x)}$
\end{center}
\end{itemize} 
\subparagraph*{} Quant aux valeurs absolues non ultra-métriques, on sait (Ostrowski) qu'elles sont de la forme : 
$|.| = |f(x)|_{\infty}^c$, avec $0< c \leqslant 1$ et $|.|_{\infty}$ est la valeur absolue (module) sur $\mathbb{C}$, 
où $f: K \mapsto \mathbb{C}$ est un isomorphisme de $K$ sur un sous-corps du corps des nombres complexes. 

\paragraph*{Théorème :} 
\subparagraph*{} Soit $\alpha \in ]0, 1[$. Alors la correspondance suivante :
\begin{center}
$ \alpha : V(K) \rightarrow Val^u(K)$\\
$\upsilon \mapsto |.|_{\upsilon}$
\end{center}
est une bijection.\\
En effet, puisque pour tout $0<\alpha <1$ l'application $ t\longmapsto \alpha^{t} $ est un isomorphisme du groupe ordonné $\mathbb{R}$ sur le groupe ordonné $\mathbb{R}^{*}_{+}$.
\paragraph*{Remarque :} D'après cette petite partie préliminaire à propos de la relation entre les valuations et les valeurs absolues, on conclut qu'il est équivalent de raisonner en termes de valeur absolue ultra-métrique ou en termes de valuation.

\subsection*{1.2. Anneaux de valuation}
\addcontentsline{toc}{subsection}{1.2. Anneaux de valuation}
\subsubsection*{1.2.1. Cas général :}
\addcontentsline{toc}{subsubsection}{1.2.1. Cas général }
\paragraph*{Théorème :} 
\subparagraph*{} Soit $K$ un corps, et $V$ un sous-anneau 
de $K$. Les conditions suivantes sont équivalentes : 
\begin{itemize}
\item[a)] $V$ est un élément maximal de l'ensemble des sous-anneaux locaux de $K$, cet ensemble étant ordonné par la relation " $B$ domine $A$ " entre $A$ et $B$. 
\item[b)] Il existe un corps algébriquement clos $L$, et un homomorphisme $h$ de $V$ dans $L$ qui est maximal dans l'ensemble des homomorphismes de sous-anneaux de $K$ dans $L$, ordonné par la relation " $g$ est un prolongement de $f$ " entre $f$ et $g$. 
\item[c)] Si $x \in K \setminus V$, alors $x^{-1} \in V$.
\item[d)] Le corps des fractions de $V$ est $K$, et l'ensemble des idéaux principaux de $V$ est totalement ordonné par la relation d'inclusion.
\item[e)]  Le corps des fractions de $V$ est $K$, et l'ensemble des idéaux de $V$ est totalement ordonné par la relation d'inclusion. 
\end{itemize}  
\paragraph{Définitions :}
\subparagraph*{Anneau de valuation pour un corps $K$ :} On dit que $V$ est un anneau de valuation pour le corps $K$ si les conditions équivalentes a), b), c), d), e) sont satisfaites.
\subparagraph*{Anneau de valuation :} On dit qu'un anneau est un anneau de valuation s'il est intègre et si c'est un anneau de valuation pour son corps des fractions. 
\paragraph{Exemples d'anneaux de valuation :}
\subparagraph*{$\ast$} Tout corps est un anneau de valuation. 
\subparagraph*{$\ast$} Pour $p$ un nombre premier, l'anneau local $\mathbb{Z}_{(p)}$  le sous-ensemble du corps $\mathbb{Q}$ des rationnels formé des fractions $\dfrac{r}{s}$, où $s$ n'est pas divisible par $p$ est un anneau de valuation. 

\subsubsection*{1.2.2. Cas des corps valués non-archimédiens :}
\addcontentsline{toc}{subsubsection}{1.2.2. Cas des corps valués non-archimédiens }
Soient $(K, |.|)$ un corps valué non-archimédien.
\subparagraph*{} On a défini dans la partie précédente la boule ouverte de centre $a$ et de rayon $r$ par : 
\begin{center}
 $B_o(a,r) = \{x \in K : |a - x| < r \}$.
 \end{center}
  et la boule fermée de centre $a$ et de rayon
$r$ par :
\begin{center}
 $B_f(a,r) = \{x \in K : |a - x| \leqslant r \}$
\end{center} 
\subparagraph*{} On a vu également que la boule ouverte $B_o (a, r)$ est un fermé de $(K, | . |)$ et la boule fermée $B_f (a, r)$ est un ouvert de $(K, | . |)$. Par conséquent, les boules $B_o (0, 1)$ et $B_f (0, 1)$ sont les boules unités ouverte et fermée respectivement.

\paragraph*{Proposition :\\}
 
\begin{itemize}
\item[$i\text{ }\text{ }-$] La boule $B_f(0, 1)$ est un sous anneau de $K$ qui contient $1_K$.
\item[$ii\text{ }-$] La boule $B_o (0, 1)$ est l'ensemble des éléments non inversibles de $B_f(0, 1)$ et est son unique idéal maximal.
\item[$iii-$]  Le quotient $B_f(0, 1) /Bo (0, 1)$ est un corps.
\end{itemize}

\paragraph*{Démonstration :}
\subparagraph*{i-} On a : $0_K \in B_f(0, 1)$ , d'où $B (0, 1) \neq \emptyset$, et pour tout $(x, y) \in (B_f(0, 1))^2$ 
\begin{center}
$|xy| \leqslant 1$ et $| x - y | \leqslant max (| x |, | y |) \leqslant 1$
\end{center}
d'où :
\begin{center}
$xy \in B_f(0, 1)$ et $x - y \in B_f(0, 1)$ 
\end{center}
Et il est bien évident que $1_K \in B_f(0, 1)$ .\\
Donc :
\begin{center}
$B_f(0, 1)$ est un sous anneau de $K$.
\end{center}
\subparagraph*{ii-} Il est évident que $B_o (0, 1)$ est un sous anneau de $K$.\\ 
Et pour tout $x \in B_f(0, 1)$ et $y \in B_o (0, 1)$ on a :
\begin{center}
$| xy | =| x |.| y | < 1$
\end{center}
Donc $xy \in B_o (0, 1)$ , et par suite $B_o (0, 1)$ est un idéal de $B_f(0, 1)$ .\\
Montrons que $B_o (0, 1)$ est l'ensemble des éléments non inversibles de $B_f(0, 1)$:\\
Soit $x \in B_o (0, 1)$ ; supposons que $x$ est inversible dans $B_f(0, 1)$ alors :
\begin{center}
$\exists y \in B_f(0, 1)$ tel que $ xy = 1$ et $y \neq 0$
\end{center}
Alors :
\begin{center}
$| x | . | y |= 1$
\end{center}
Ou encore :
\begin{center}
$| x |= \dfrac{1}{|y|}$
\end{center}
Or $| y | \leqslant 1$ d'où :
\begin{center}
$| x |\geqslant 1$
\end{center}
Ce qui est en contradiction avec le fait que 
\begin{center}
$x \in B_o (0, 1)$ $(| x |< 1)$
\end{center}
Inversement, soit $x \in B_f(0, 1)$ : $| x |= 1$ et $(\exists y \in K) : xy = 1$, alors :
\begin{center}
$| y | =| x | . | y |=| x.y |=| 1 |= 1$.
\end{center}
Donc $y \in B_f(0, 1)$ et $x$ est inversible dans $B_f(0, 1)$.\\ 
Le disque unité $B_f(0, 1) \setminus B_o (0, 1)$ de $B_f(0, 1)$ étant égal à l'ensemble des éléments inversibles de $B_f(0, 1)$ donc $B_o (0, 1)$ est un idéal maximal de $B_f (0, 1)$ ;
et il est bien évident que c'est l'unique ayant cette propriété. On le notera $P (| . |)$ ou tout simplement $P$.
\subparagraph*{iii-} 
La démonstration repose sur le théorème -qu'on admet- suivant :
\begin{center}
"Soit $A$ un anneau commutatif unitaire et $I$ un idéal de $A$
tel que $I \neq A$. On a les équivalences suivantes :\\
$I$ est un idéal premier de $A \Leftrightarrow A/I$ est un anneau intègre.\\
$I$ est un idéal maximal de $A \Leftrightarrow  A/I$ est un anneau corps."
\end{center}
\begin{flushright}
$\blacksquare$
\end{flushright}

\paragraph*{Définition (Anneau de valuation de la valeur absolue) :}
\subparagraph{} Soit $(K, | . |)$ un corps valué non-archimédien. $B_f(0, 1)$ s'appelle l'anneau d'intégrité de $(K, | . |)$ ou l'anneau de valuation de la valeur absolue $| . |$; on le note  $i (| . |)$ ou $i (K)$ ou $O$.
\subparagraph{} Le corps $O/P$ s'appelle le corps résiduel de $| . |$ on le note $\kappa$ ou $\overline{K}$.

\paragraph*{Définition (Anneau de valuation d'une valuation) :}
\subparagraph{} Soit $K$ un corps commutatif muni d'une
valuation $\upsilon$ . Les éléments de $K$ de
valuation positive ou nulle constituent un
sous-anneau $\mathcal{O}$ appelé l'anneau de
valuation associé à la valuation $\upsilon$ sur $K$ :
\begin{center}
$\mathcal{O} = \{ x \in K ; \upsilon(x) \geqslant 0 \}.$
\end{center}
le corps des fractions de $\mathcal{O}$ est $K$.
\subparagraph*{} On a $\upsilon(1/x) = -\upsilon(x)$ pour tout élément non
nul $x$ de $K$, et donc $x$ est un élément inversible de $\mathcal{O}$ si et seulement si $\upsilon(x) \neq 0$. Par conséquent, $\mathcal{O}$ est un anneau local dont l'unique idéal maximal $m(\mathcal{O})$ est constitué des éléments de valuation strictement positifs :
\begin{center}
$m(\mathcal{O}) = \{ x \in ; \upsilon(x) > 0 \}.$
\end{center}
\paragraph*{Caractérisations des anneaux de valuation :}
\subparagraph*{} Il existe diverses caractérisations des anneaux de valuation (on les énonce sans démonstrations) :\\
Soient $\mathcal{O}$ un anneau intègre et $K$ son corps des fractions. Les conditions suivantes sont équivalentes :
\begin{itemize}
\item[1.]$\mathcal{O}$ est un anneau de valuation (pour une certaine valuation sur $K$) ;
\item[2.] Pour tout élément $x$ de $K$ qui n'appartient pas à $\mathcal{O}$, l'inverse de $x$ appartient à $\mathcal{O}$ ;
\item[3.] Sur l'ensemble des idéaux principaux de $\mathcal{O}$, l'ordre défini par l'inclusion est total ; 
\item[4.] Sur l'ensemble des idéaux de $\mathcal{O}$, l'ordre défini par l'inclusion est total.
\end{itemize}

\paragraph*{Remarque :}
Deux valuations $\upsilon$ et $\upsilon$' sur $K$ sont équivalentes si et seulement si elles ont
le même anneau de valuation.

\subsection*{1.3. Anneaux de valuation discrète}
\addcontentsline{toc}{subsection}{1.3. Anneaux de valuation discrète}
\subsubsection*{1.3.1. Définition et exemple }
\addcontentsline{toc}{subsubsection}{1.3.1. Définition et exemple }
\paragraph{Définition (Anneaux de valuation discrète) :}

\subparagraph*{Première définition :} Un anneau de valuation discrète est un anneau de valuation dont la valuation est discrète mais non triviale.

\subparagraph*{} Autrement dit, $A$ est un anneau commutatif unitaire intègre, et il existe sur son corps des fractions $K$ une valuation $\upsilon$, à valeurs entières mais non toutes nulles, telle que :
\begin{center}
$ {\displaystyle A=\{x\in K\mid v(x)\geqslant 0\}}$
\end{center}

Par conséquent (comme tout anneau d'une valuation non triviale) $A$ est un anneau local mais pas un corps, et son unique idéal maximal $M$ est non nul, et constitué des éléments de valuation strictement positive :
\begin{center}
$ \displaystyle M=\{x\in K\mid v(x)>0\}$
\end{center}

De plus (comme la valuation est à valeurs entières) tout idéal est engendré par n'importe lequel de ses éléments de valuation minimum, si bien que $A$ est principal. En particulier, un générateur de $M$ est appelé uniformisante ou paramètre local de l'anneau.

La réciproque est claire : tout anneau local et principal qui n'est pas un corps est un anneau de valuation discrète. On pose $\upsilon(a)$ égal à l'entier naturel $n$ tel que $aA = M^n$. On obtient donc une définition équivalente :
 
\subparagraph{Seconde définition :} On appelle un anneau $A$ anneau de valuation discrète si c'est un anneau principal, et s'il possède un et un seul idéal premier non nul $m(A)$. 
{\small [On rappelle qu'un idéal $p$ d'un anneau commutatif $A$ est dit premier si l'anneau quotient $A/p$ est intègre.] }
\paragraph{Exemple d'anneau de valuation discrète :}
\subparagraph{} Soit $p$ un nombre premier, et soit $\mathbb{Z}_{(p)}$ le sous-ensemble du corps $\mathbb{Q}$ des rationnels formé des fractions $\dfrac{r}{s}$, où $s$ n'est pas divisible par $p$; c'est un anneau de valuation discrète de corps résiduel le corps $\mathbb{F}_p$ à $p$ éléments. Si $\upsilon_p$ désigne la valuation associée, 
$\upsilon_p (x)$ n'est autre que l'exposant de $p$ dans la décomposition de $x$ en facteurs premiers.

\subsubsection*{1.3.2. Caractérisations des anneaux de valuation discrète }
\addcontentsline{toc}{subsubsection}{1.3.2. Caractérisations des anneaux de valuation discrète}
\paragraph{Proposition 1 :} 
\subparagraph{} Soit $A$ un anneau commutatif. Pour que $A$ soit un anneau de valuation discrète, il faut et il suffit que ce soit un anneau local noethérien, et que son idéal maximal soit engendré par un élément non nilpotent. 

\paragraph{Proposition 2 :} 
\subparagraph{} Soit $A$ un anneau intègre noethérien. Pour que $A$ soit un anneau de valuation discrète, il faut et il suffit qu'il vérifie les deux conditions suivantes :
\begin{itemize}
\item[$i\textit{ }- $] $A$ est intégralement clos. 
\item[$ii- $] $A$ possède un idéal premier non nul et un seul. 
\end{itemize}

\section*{2. Corps valués locaux}
\addcontentsline{toc}{section}{2. Corps valués locaux}
\subsection*{2.1. Compacité des corps valués }
\addcontentsline{toc}{subsection}{2.1. Compacité des corps valués }
En considérant un corps valué comme étant un \textit{espace métrique} qui est un \textit{corps topologique} (muni d'une topologie pour laquelle toutes les opérations de corps sont continues), on peut parler de la notion de compacité de corps valués.

\paragraph*{Définition (Corps valué compact) :}
\subparagraph*{} Soit $K$ un corps valué séparé. On dit que $K$ est compact si de tout recouvrement ouvert de $K$, on peut extraire un sous-recouvrement fini.\\
Autrement dit, pour toute famille d'ouverts $(U_i)_{i\in I}$ de $K$ telle que :
\begin{center}
$K =\displaystyle \bigcup_{i\in I} U_i $
\end{center}
Il existe un sous-ensemble fini $J$ de $I$ tel que : \begin{center}
$K =\displaystyle \bigcup_{j\in J} U_j $
\end{center}
Avec $(U_i)$ est une réunion de boules ouvertes définies à partir de la distance qui provient de la valeur absolue (ou de la valuation).
\subsection*{2.2. Compacité locale des corps valués }
\addcontentsline{toc}{subsection}{2.2. Compacité locale des corps valués }
\subsubsection*{2.2.1. Définitions et exemples}
\addcontentsline{toc}{subsubsection}{2.2.1. Définitions et exemples}
\paragraph*{Définition (Corps valué localement compact) :}
\subparagraph*{} Soit $K$ un corps valué séparé. On dit que $K$ est localement compact si tout point de $K$ admet un voisinage compact.

\paragraph*{Définition (Corps valué local) :}
\subparagraph*{} Un corps valué $(K, |.|)$ est dit \textbf{local} s'il est localement compact pour une topologie non discrète. 
\subparagraph*{} Si la topologie provient d'une valeur absolue archimédienne, le corps $K$ est dit corps valué local archimédien. 
\subparagraph*{} Si la topologie provient d'une valeur absolue non-archimédienne, le corps $K$ est dit corps valué local non-archimédien.
\subparagraph*{} Notons que pour une valuation discrète on peut donner la définition suivante:
\subparagraph*{} Soit $(K,|.|)$ un corps valué complet par rapport à la valeur absolue $|.|$ induite par une valuation discrète $\upsilon$. On dit que $K$ est un corps local si son corps résiduel $\kappa(K)$ est fini.
 
\paragraph*{Exemples :} 
\subparagraph*{$\ast$} Pour tout nombre premier $p$ le corps des nombres p-adiques $\mathbb{Q}_p $ est localement compact, donc un corps valué local.
\subparagraph*{}En effet, la démonstration repose sur le fait que $\mathbb{Z}_p$ soit compact 
, et que les boules de $\mathbb{Q}_p$ soient isomorphes à $\mathbb{Z}_p $ donc compactes. Remarquons pour commencer que l'image
de $\mathbb{Q}^*_p$ par la valeur absolue p-adique est $p^{\mathbb{Z}}$ qui est un sous-groupe discret de $\mathbb{R}^*_+$. En conséquence toute boule fermée de rayon strictement positif est ouverte et toute boule ouverte de rayon strictement positif est fermée (on a déjà montrer ce résultat d'une autre façon dans la partie précédente) : l'espace topologique $\mathbb{Q}_p$ est totalement discontinu
(c'est la définition même, tout point possède une base de voisinages à la fois ouverts et fermés).\\
À cause de l'inégalité ultramétrique : $|x + y|_p \leqslant max\{|x|_p, |y|_p\}$, les boules de $\mathbb{Q}_p$ centrées en $0$ sont des sous-groupes additifs. La multiplicativité $|xy|_p = |x|_p|y|_p$ nous dit que toute boule de rayon inférieur ou égal à 1 centrée en 0 est stable par multiplication.\\
 Examinons en particulier la boule unité fermée de $\mathbb{Q}_p$ :
\begin{center}
$\mathbb{Z}_p = \{x \in \mathbb{Q}_p , |x|_p \leqslant 1\}$
\end{center}
C'est l'anneau des entiers $p$-adiques, formé des sommes
\begin{center}
$\sum\limits_{i\geqslant 1} a_ip^i$ $a_i \in \{0, ... , p - 1\}$ ;
\end{center}
 en particulier $\mathbb{Z}_p$ est l'adhérence de $\mathbb{Z}$ dans $\mathbb{Q}_p$.
La boule unité ouverte est l'idéal maximal $p\mathbb{Z}_p$ de $\mathbb{Z}_p$, soit l'ensemble des sommes 
$\sum\limits_{i\geqslant 1} a_ip^i$, $a_i \in \{0, ... , p - 1\}$. Enfin les boules fermées
\begin{center}
$p^i\mathbb{Z}_p = \{x \in \mathbb{Q}_p | \upsilon_p(x) > i\}$, $i \in \mathbb{Z}$
\end{center} forment un système fondamental
de voisinages de 0.\\
Nous allons montrer que la boule unité $\mathbb{Z}_p$ est compacte. Il en résultera que l'espace topologique $\mathbb{Q}_p$ est \textit{localement compact}.

\subparagraph*{} Pour montrer la compacité de $\mathbb{Z}_p$ : puisque c'est un espace métrique, il suffit de remarquer qu'il est séquentiellement compact (un espace séquentiellement compact est un espace topologique dans lequel toute suite possède au moins une sous-suite convergente), comme produit dénombrable d'espaces qui le sont : ${\displaystyle \mathbb{Z}_{p}\simeq \left(\{0,1,\dots ,p-1\}\right)^{\mathbb {N} }}$, d'après la représentation en série.
\subparagraph*{} La compacité locale de $\mathbb{Q}_p$ découle donc du fait que $\forall x \in \mathbb{Q}_p$, $ x+\mathbb{Z}_p$ est voisinage compact de $x$.
\subparagraph*{$\ast$} Si $K$ est un corps fini, le corps $K((X))$ des séries formelles de Laurent à coefficients dans $K$ est un corps local.\\
Avec :
\\ Une série formelle en $X$ sur le corps $K$ est une expression :
\begin{center}
$\sum\limits_{n\geqslant 1}a_nX^n = a_0 + a_1X =a_2X^2 + ...$
\end{center}
\subparagraph*{} Soit $K[[X]]$ l'anneau des séries formelles sur $K$. Le corps des fractions de $K[[X]]$ est appelé corps des séries de Laurent sur $K$ et noté $K((X))$.
L'anneau des entiers de $K((X))$ est $K[[X]]$ et son corps résiduel est $K$, qui est fini par hypothèse.

\subparagraph*{$\ast$} Le corps des nombres réels $\mathbb{R}$ muni de la valeur absolue usuelle est un corps valué local. 

\subparagraph*{$\ast$} L'ensemble des rationnels $\mathbb{Q}$, pour sa topologie de sous-espace de $\mathbb{R}$, n'est pas localement compact. En effet, comme $\mathbb{Q}$ est dénombrable, alors $\mathbb{R}$ est séparable et dénombrable à l'infini. Soit $x \in \mathbb{Q}$. Si $x$ possède un voisinage compact $K$ dans $\mathbb{Q}$, alors $K$ est un voisinage compact de $x$ dans $\mathbb{R}$, donc il existe $r > 0$ tel que $[x - r, x + r] \subset K \subset \mathbb{Q}$, ce qui est impossible. Par conséquent, $\mathbb{Q}$ n'est pas localement compact.
\subparagraph*{$\ast$} Tout espace compact est localement compact, car il est un voisinage compact de chacun de ses points. Donc tout corps valué compact est un corps local.

 \paragraph*{Remarque :}
\subparagraph*{$\ast$} Soit $(K,\upsilon)$ un corps muni d'une valuation discrète, pour que $K$ soit local, il faut et il suffit qu'il soit complet et que son corps résiduel $\kappa(K)$ soit un corps fini.
\subparagraph*{$\ast$} On a vu dans le début de cette partie qu'il y a une relation biunivoque entre l'ensemble des valuations et l'ensemble des valeurs absolues ultramétriques et cette relation repose sur le choix d'un nombre $\alpha < 1$ quelconque, on ajoute que lorsque le corps $K$ vérifie les conditions de la dernière définition (correspond à une valuation discrète), il y a une façon canonique de choisir le nombre $\alpha$ : on prend $a = q^{-1}$, où $q$ est le nombre d'éléments du corps résiduel $\kappa(K)$. La valeur absolue correspondante est dite normalisée. La proposition suivante en donne une caractérisation " analytique " : 
\paragraph*{Proposition :}
\subparagraph*{} Soit $K$ un corps vérifiant les conditions de les conditions de la dernière définition (muni d'une valuation discrète), et soit $\mu$ une mesure de \textsc{Haar} du groupe additif localement compact $K$. Pour toute partie mesurable $E$ de $K$, et pour tout $x \in K$, on a alors :
\begin{center}
$\mu(xE) = |x| \mu(E)$ 
\end{center} où $|x|$ désigne la valeur absolue normalisée de $x$. 
\subparagraph*{} Si on prend $E = A$ ( l'anneau de valuation de $\upsilon$ (discret)), on voit que $E$ est réunion de $(A : xA)$ classes modulo $xE,$ d'où $\mu(E) = (A : xA).\mu(xE)$, et 
$|x| = 1/(A : xA)$. Comme $(A : xA)$ est égal à $q^{-\upsilon(x)}$, on trouve bien :
\begin{center}
$q^{-\upsilon(x)} = |x|$
 \end{center} 
avec $(A : xA)$ est Le cardinal de l'ensemble des classes à gauche (modulo $xA$) qu'on appelle l'indice du sous-groupe $xA$ par rapport à $A$; il est égal au cardinal de l'ensemble des classes à droite. 

\begin{flushright}
$\blacksquare$
\end{flushright}
\paragraph*{Quelques rappels :}
\subparagraph*{$\ast$ Mesure :} Une application $\mu : \mathcal{F} \rightarrow [0; +\infty]$ est appelée mesure positive, ou simplement mesure, sur $(\Omega; \mathcal{F})$ si elle vérifie les deux propriétés suivantes :
\begin{itemize}
\item[1.] $\mu(\emptyset) = 0 $.
\item[1.] $\sigma$-additivité : Pour toute suite $\{A_n; n \geqslant 1 \}$ d'éléments de $\mathcal{F}$ deux-à-deux disjoints :
\begin{center}
$\mu$ $\left(  \bigcup\limits_{n \geqslant 1}^{} A_n \right) $ = $\sum\limits_{n\geqslant 1}^{} \mu (A_n)$
\end{center}
\end{itemize}
\subparagraph*{$\ast$ Mesure de Haar :} Soit G un groupe localement compact. On appelle mesure de Haar à gauche (resp. à droite) sur G une mesure positive non nulle sur G, invariante à gauche (resp. à droite). 
\begin{center}
$\lambda(gB)= \lambda(B)$.
\end{center}
\subparagraph*{} L'existence d'une mesure de Haar est assurée dans tout groupe localement compact. Elle est finie sur les parties compactes de $G$. 

\subsubsection*{2.2.2. Fonction module sur un corps local}
\addcontentsline{toc}{subsubsection}{2.2.2. Fonction module sur un corps local}
\paragraph*{Définitions :}
\subparagraph*{}
Soient $G$ un groupe localement compact (un groupe dont l'espace sous-jacent est localement compact), $\varphi$ un automorphisme de $G$, $\mu$ une mesure de Haar à gauche sur $G$. Il est clair que $\varphi^{-1}(\mu)$ est encore une mesure de Haar à gauche sur $G$. Il existe donc un nombre $a > 0$ et un seul tel que $\varphi^{-1}(\mu) = a\mu$. Ce nombre est indépendant du choix de $\mu$. Remarquons que, si l'on partait d'une mesure de Haar à droite, par exemple $\Delta_G^{-1}.\mu$ , on aboutirait au même scalaire $a$ : car, comme $\varphi^{-1}$ laisse $\Delta_G$ invariant, on a 
$\varphi^{-1}(\Delta^{-1}_G$.$\mu)$ = $\Delta_G^{-1}$.$\varphi^{-1}(\mu)= a\Delta_G^{-1}.\mu$. 
\paragraph*{Module d'un automorphisme :} Le nombre $a > 0$ tel que $\varphi^{-1}(\mu) = a\mu$ s'appelle le module de l'automorphisme $\varphi$ et se note $mod_{G\varphi}$ ou simplement $mod_\varphi$. 
\paragraph*{Fonction module :} Soit $K$ un corps localement compact (non nécessairement commutatif). On définit la fonction mod (ou $mod_K$) sur $K$ comme suit : on a $mod_K(0) = 0$ et pour $x\neq 0$ dans $K$, le nombre $mod_K(x)$ est le module de l'automorphisme $y \mapsto xy$ du groupe additif de $K$. 
\paragraph*{Proposition 1 :} 
\subparagraph*{} Si $K$ est un corps local,la fonction $mod_K$ appartient à $\vartheta(K)$. En outre :
\begin{itemize}
\item[i- ] Si $s > 0$ est tel que $(mod_K)^s = g$ soit une valeur absolue, alors $g$ définit la topologie de $K$.
\item[i- ] Si $K$ est non discret et si $mod_K$ est une valeur absolue ultramétrique, il existe une valuation discrète normée $v$ sur $K$, dont l'anneau est compact et le corps résiduel fini à $q$ éléments, de sorte que $mod_K = q^{-\upsilon}$. La topologie de $K$ est définie par $\upsilon$. 
\end{itemize} 
\paragraph*{Proposition 2 :}
\subparagraph*{} Soient $K$, $K'$ deux corps locaux (non nécessairement commutatifs) tels que $K$ soit un 
sous-corps topologique de $K'$ et que $K$ soit non discret. Alors : 
\begin{itemize}
\item[$i -$] $K'$ est un espace vectoriel à gauche (resp. à droite) de dimension finie sur $K$. 
\item[$ii-$] Si $K$ est contenu dans le centre de $K'$, on a, pour tout $x \in K'$
\begin{center}
 $mod_{K'}(x) = mod_{K'/K}(K(x))$.
 \end{center} 
\end{itemize}

\paragraph*{Théorème :} 

\subparagraph*{} Soit $K$ un corps valué complet non discret. Un espace vectoriel topologique séparé $E$ sur $K$ qui admet un voisinage de $0$ précompact $V$ est de dimension
finie. Si $E$ n'est pas réduit à 0, $K$ et $E$ sont alors localement compacts.
\paragraph*{Démonstration :}
\subparagraph*{} Pour démontrer la première assertion, on peut se limiter au cas où $E$ est complet, car $E$ est un sous-espace partout dense de son complété $\widehat{E}$ et l'adhérence $\overline{V}$ de $V$ dans $\widehat{E}$ est compacte et est un voisinage de 0 dans $\widehat{E}$. 
\subparagraph*{} On peut donc supposer qu'il y a dans $E$ un voisinage compact $V$ de $0$. Soit $a \in K$ tel que $0 < |a| < 1$ ; il y a donc des points $a_i \in V$ en nombre fini tels que
\begin{center}
$V \subset \bigcup\limits{i}(a_i+ \alpha V)$
\end{center}
Soit $M$ le sous-espace (de dimension finie) de $E$ engendré par les $a_i$ ; il est fermé dans $E$ ; dans l'espace vectoriel topologique séparé $E/M$, l'image canonique de $V$ est un voisinage compact $W$ de $0$ tel que $W \subset aW$ ; ceci s'écrit encore $a^{-1}W \subset W$, d'où par récurrence sur $n$, a $\alpha^{-n} W \subset W$ pour tout entier positif $n$. Comme $W$ est absorbant, on en déduit que $W = E/M$ ; autrement dit $E/M$ est compact. Pour prouver la première assertion, il suffit donc de démontrer le lemme suivant :
\paragraph*{Lemme :} Tout espace vectoriel topologique compact $E$ sur un corps valué non discret est réduit à $0$.
\subparagraph*{} En effet, comme $E$ est complet, on peut supposer qu'il en est de même pour $K$.
Si $E$ n'était pas réduit à $0$, il contiendrait une droite, fermée dans $E$, donc compacte,et isomorphe à $K_s$, et par suite $K$ serait compact ; mais cela est absurde, car l'application $\xi \mapsto |\xi|$ de $K$ dans $\mathbb{R}$ est continue, donc serait bornée, alors qu'il existe des $y \in K$ tels que $|y| > 1$, donc tels que $|y^n| = |y|^n$ soit
arbitrairement grand.
\subparagraph*{} Revenant au théorème, on voit que si $E$ admet un voisinage de $0$ précompact et n'est pas réduit à $0$, $E$ est de dimension finie sur $K$, donc isomorphe à un espace $K^n$ avec
$n > 0$ ; comme $K$ est complet, il en est de même pour $E$, qui est donc localement compact. Puisque $K_s$ est isomorphe à une droite de $E$ 
 , nécessairement
fermée dans $E$ 
, $K$ est localement compact.
\begin{flushright}
$\blacksquare$
\end{flushright}
\paragraph*{Remarque :} La conclusion du théorème ne subsiste plus lorsque $K$ est un corps discret, comme le montre l'exemple de $\mathbb{R}$ (muni de la topologie usuelle) considéré comme espace vectoriel topologique sur le corps $\mathbb{Q}$ discret.

\paragraph{Corollaire 1 :}
\subparagraph*{} Tout corps localement compact dont le 
centre est non discret est de rang fini sur son centre. 
\subparagraph*{} En effet, le centre $\mathbb{Z}$ d'un corps localement compact $K$ est fermé dans $K$, donc localement compact. 

\paragraph{Corollaire 2 :}
\subparagraph*{} Soient $K'$ un corps localement compact 
et $K$ un sous-corps fermé de $K'$ (non nécessairement commutatifs). Si $K'$ est un espace vectoriel à gauche (resp. à droite) de dimension finie $n$ sur $K$, on a 
\begin{center}
$mod_{K'}(x) = (mod_K(x))^n$ pour tout $x \in \mathbb{R}$. 
\end{center}

En effet, de façon générale, on sait que dans un espace vectoriel (à gauche ou à droite) de dimension finie $n$ sur $K$, l'homothétie de rapport $x \in K$ a un module égal à $(mod_K(x))^n$; il suffit d'appliquer cela à $K'$. 
\begin{flushright}
$\blacksquare$
\end{flushright}

\subsubsection*{2.2.3. Quelques propriétés topologiques des corps valués locaux}
\addcontentsline{toc}{subsubsection}{2.2.3. Quelques propriétés topologiques des corps valués locaux}
\paragraph*{Remarque :}
Il est clair que tout espace compact est localement compact, mais la réciproque n'est pas vraie ; par exemple, tout espace discret est localement compact, mais pas compact s'il est infini.

\paragraph*{Proposition 1 :}
\subparagraph*{} Soit $K$ un corps local, alors l'ensemble des voisinages fermés d'un point quelconque de $K$ est un système fondamental de voisinages de ce point.\\
En outre, un espace séparé qui vérifie cette propriété est dit \textbf{espace régulier}. 
\paragraph*{Démonstration :}
\subparagraph*{} Tout point $x$ d'un corps local (localement compact) $K$ possède un voisinage compact $V$.
\subparagraph*{$\ast$} Comme $K$ est séparé, $V$ est fermé :\\
Dans le corps $K$ qui est un espace topologique séparé, tout ensemble compact est fermé. En effet, si $V$ est une partie compacte de l'espace séparé $K$ et $x$ un point de $V^c$, il résulte alors de la proposition vue en S5 :
\begin{center}
"Soient $X$ un espace séparé, $A$ et $B$ deux parties compactes de $X$ sans point commun. Alors il existe un voisinage $V$ de $A$ et un voisinage $W$ de $B$ qui ne se rencontrent pas."
\end{center}
qu'il y a un voisinage de $x$ ne rencontrant pas $V$, puisque $\{x\}$ est fermé (car dans un espace séparé, tout ensemble fini est fermé.\footnote{Il suffit de remarquer que tout ensemble réduit à un point est fermé en vertu de l'axiome \textsc{Hausdorff} : L'intersection des voisinages fermés d'un point quelconque de $X$ est l'ensemble réduit à ce point.} 
Donc $V^c$ est ouvert, d'où $V$ est un fermé.
\subparagraph*{$\ast$} $V$ est un sous-espace régulier :\\
Tout espace compact est régulier. (On admet ce résultat car la démonstration repose sur des notion avancer qui dépasse le niveau du licence).
\subparagraph*{$\ast$} $K$ est régulier :\\
Si, dans un espace topologique $K$, tout point possède un voisinage fermé qui est un sous-espace régulier de $K$, $K$ est régulier. En effet, $K$ est séparé ; d'autre part, soit $x$ un point quelconque de $K$ et soit $V$ un voisinage fermé de $x$ dans $K$, qui est un sous-espace régulier de $K$. Pour tout voisinage $U$ de $x$ dans $K$ tel que $U \subseteq V$, $U$ est un voisinage de $x$ dans $V$, donc il existe par hypothèse un voisinage $W $ de $x$ dans $V$, fermé
dans $V$ et contenu dans $U$. Mais $W$ est aussi un voisinage de $x$ dans $K$ puisque $V$ est un voisinage de $x$ dans $K$, et $W$ est fermé dans $K$ puisque $V$ est fermé dans $K$.
 
\begin{flushright}
$\blacksquare$
\end{flushright}

\paragraph*{Proposition 2 :}
\subparagraph*{} Dans un corps local $K$, tout ensemble compact $X$ admet un système fondamental de voisinages compacts.
\paragraph*{Preuve :}
\subparagraph*{}En effet, l'intersection d'un voisinage fermé de $x$ et d'un voisinage compact de $x$ est un voisinage compact de $x$, Car dans un espace compact, tout ensemble fermé est compact, Il suffit d'appliquer l'axiome "\textit{ Toute famille d'ensembles fermés dans $X$, dont l'intersection est vide, contient une sous-famille finie, dont l'intersection est vide}",en remarquant que si $A$ est fermé dans un espace $X$, tout ensemble fermé dans $A$ est fermé dans $X$.

\begin{flushright}
$\blacksquare$
\end{flushright}

\paragraph*{Proposition 3 :}
\subparagraph*{} Tout corps local est un espace de Baire, c'est à dire que le théorème de \textsc{Baire} s'y applique : toute intersection dénombrable d'ouverts denses est dense ou, autrement dit, toute union dénombrable de fermés d'intérieur vide est d'intérieur vide.
\paragraph*{Démonstration :}
\subparagraph*{} Dans ce qui suit, int($A$) désigne l'intérieur d'une partie $A$ de $K$.\\
Soit $K$ un corps localement compact. Nous utiliserons que dans $K$, tout ouvert non vide contient un compact d'intérieur non vide. En effet, tout ouvert contenant un point $x$ contient un voisinage compact de $x$, puisque $x$ possède un système fondamental de voisinages compacts.\\
Soient $(U_n)_{n\in \mathbb{N}}$ une suite d'ouverts denses dans $K$ et $V$ un ouvert non vide quelconque ; nous voulons montrer que l'intersection des $U_n$ rencontre $V$.\\
Puisque $U_0$ est dense, il rencontre $V$. L'ouvert $U_0 \cap V$ étant non vide, il contient un compact $K_0$ d'intérieur non vide. Une fois $K_0$ choisi, $U_1 \cap \text{int} (K_0)$ est un ouvert non vide, donc contient un compact $K_1$ d'intérieur non vide. En itérant cette construction, on obtient une suite décroissante de compacts non vides $K_n$ tels que $K_0 \subset V$ et $\forall n\in \mathbb{N}$, $K_n \subset U_n$.\\
L'intersection de ces compacts est donc incluse dans $ V\cap \bigcap_{n\in \mathbb{N}} U_n $, or cette intersection des $K_n$ est non vide d'après la propriété de Borel-Lebesgue. En effet, les $K_n$ sont des fermés du compact $K_0$ et toute intersection d'un nombre fini d'entre eux est non vide (puisqu'ils forment une suite décroissante de parties non vides). Finalement, $ V\cap \bigcap_{n\in \mathbb{N}} U_n  \neq \varnothing$, ce qui prouve le résultat.
 \begin{flushright}
$\blacksquare$
\end{flushright}
\paragraph*{Remarque :} On peut donner une preuve express que $\mathbb{Q}$ est non localement compact car $\mathbb{Q}$ n'est pas de Baire.


\subsection*{2.3. Corps valués locaux archimédiens}
\addcontentsline{toc}{subsection}{2.3. Corps valués locaux archimédiens}
\paragraph*{Deuxième théorème d'Ostrowski :}
\subparagraph*{} Soit $K$ un corps local complet pour une valeur absolue archimédienne $|.|$, donc $K$ est isomorphe à $\mathbb{R}$ ou $\mathbb{C}$ munis de leurs valeurs absolues archimédiennes usuelles.
\paragraph*{Démonstration :}
\subparagraph*{} On admet que si $K$ est un corps valué archimédien donc la caractéristique de $K$ est nulle.\\

On a donc $K$ est de caractéristique égale à $0$, donc $\mathbb{Q} \subset K$, et la restriction de $|.|$ à $\mathbb{Q}$ est évidemment non triviale et en plus elle est ultra-métrique.
Donc, selon le premier théorème d'Ostrowski, on peut remplacer $|.|$ par $|.|_{\infty}$ sur $\mathbb{Q}$. Et puisque $K$ est complet, on a nécessairement $\mathbb{R} \subset K$.\\
Si on montre que tout élément de $K$ est solution d'une équation du second degré sur $\mathbb{R}$, on peut conclure que $K=\mathbb{R}$ ou $\mathbb{C}$, et sachant que la seule extension archimédienne de la valeur absolue de $\mathbb{R}$ est le module de $\mathbb{C}$ - qui représente la valeur absolue usuelle de $\mathbb{C}$ -. 
Soit $\alpha \in K$, on considère la fonction $f : \mathbb{C}\rightarrow \mathbb{R}$ définie par :
\begin{center}
$f(z)= |\alpha^2 -(z+\overline{z})\alpha + z\overline{z}|$
\end{center}
On a $f$ est continue et elle tend vers $+\infty$ lorsque $|z|_{\infty}$ tend vers $+\infty$, donc $f$ admet un minimum global $m$. L'objet est de montrer que $m = 0$, auquel cas $\alpha$ est le zéro d'une équation quadratique sur $\mathbb{R}$.\\
Posons $S =\{ z \in \mathbb{C} \mid f (z) = m \}$, on voit que $S$ est fermé puisque $f$ est continue, et il est borné car $f(z)\rightarrow \infty$ lorsque $|z|_{\infty}\rightarrow \infty$. Puisque $S$ est compact et $|.|$ est continu sur $\mathbb{C}$, on peut choisir un élément $z_{0}\in S$ tel que $|z_0|_{\infty} \geq |z|_{\infty}$ pour tout $z \in S$.\\
Supposons que $m > 0$, choisissons $\epsilon$ tel que $0 < \epsilon < m$, considérons le polynôme défini sur $\mathbb{R}$ par :
\begin{center}
$g(x) = x^2 - (z_0 + \overline{z}_0)x + z_0\overline{z}_0 + \epsilon$ ;
\end{center}
et supposons que $z_1, \overline{z}_1 \in \mathbb{C}$ sont les racines de $g(x)$. Alors $z_1\overline{z}_1 = z_0\overline{z}_0 + \epsilon$, et ainsi $|z_1|_{\infty} > |z_0|_{\infty}$. Par le choix de $z_0$, nous devons avoir $f (z_1)> m$.\\
Pour tout entier positif $n$, considérons le polynôme
$G_n (x) = (g(x) - \epsilon)^n - (-\epsilon)^n$
de sorte que $G_n (z_1) = 0$ puisque $g (x)$ est un facteur de $G_n (x)$. Soit $z_1 = \omega_1, \omega_2, ..., \omega_{2n}$ sont les zéros de $G_n (x)$ en $\mathbb{C}$. Puisque cet ensemble de zéros est le même que l'ensemble de ses conjugués, nous avons
\begin{center}
$G_n(x) = \prod\limits_{i=1}^{2n}(x - \omega_i) = \prod\limits_{i=1}^{2n}(x - \overline{\omega}_i)$ et
$G_n(x)^2= \prod\limits_{i=1}^{2n}(x^2 - (\omega_i + \overline{\omega}_i )x + \omega \overline{\omega}_i).$
\end{center}
En substituant $\alpha$ dans le polynôme $G_n (x)^2$ et en appliquant la valeur absolue $|·|$, on obtient
\begin{center}
$|G_n(x)|^2= \prod\limits_{i=1}^{2n}f(\omega_i)\geqslant f(z_1)m^{2n-1}$
\end{center}
On a aussi par la définition de $G_n$ et l'inégalité triangulaire 
\begin{center}
$|G_n(x)|^2\leqslant (f(z_0)^n + \epsilon^n)^2= (m^n + \epsilon^n)^2$
\end{center}
On a $|g(\alpha) - \epsilon|_{\infty} = f(z_0)$, donc, $f(z_1)m^{2n-1} \leqslant (m^n +\epsilon^n)^2$, d'où :
\begin{center}
 $\dfrac{f(z_1)}{m} \leqslant \left(1 + \left(\dfrac{\epsilon}{m}\right)^n\right)^2$
 \end{center} 
Puisque $\epsilon < m$, et $n$ est un entier positif arbitraire, faisons $n \rightarrow +\infty$ donne $f(z_1) \leqslant m$,
contredisant $f (z_1)> m$. Ainsi $m = 0$.

\begin{flushright}
$\blacksquare$
\end{flushright}
\paragraph*{Remarque :} Il est clair qu'on a pas utilisé le fait que le corps $K$ soit local, on en conclut que \textbf{tout corps valué archimédien est isomorphe à $\mathbb{R}$ ou $\mathbb{C}$.}

\subsection*{2.4. Corps valués locaux non-archimédiens}
\addcontentsline{toc}{subsection}{2.4. Corps valués locaux non-archimédiens}
On note $O$ la boule ouverte $B(0,1)$, c'est l'anneau de valuation pour la valeur absolue.
\paragraph*{Proposition 1 :}
\subparagraph*{} Si $(K, |.|)$ un corps valué non-archimédien local, alors :
\begin{itemize}
\item[$i   -$] K est complet.
\item[$ii  -$] Le corps résiduel de K est fini.
\item[$iii -$] La valuation de K est discrète.
\end{itemize}
\paragraph*{Démonstration :} 
\subparagraph*{$i$  -} Soit $C$ un voisinage compact de 0. Soit $\pi \in K^*$ tel que : $| \pi |< 1$.\\
Alors $\exists n \in \mathbb{N}$ suffisamment grand tel que $\pi^n.O \in C$ car sinon :
\begin{center}
 $\exists (x_n)_n \subseteq O$ telle que $\forall n \in \mathbb{N}$, $\pi^n.x_n \notin C$
 \end{center}
ce qui est en contradiction avec le fait que $\lim\limits_{n \to +\infty}$ $\pi^n.x_n = 0$.\\
Or $\pi^n.O$ est un fermé de $C$, donc $\pi^ n.O$ est compact. D'où $O$ est compact, il est donc complet. \\
Donc, si ($x_n)_n$ est une suite de Cauchy dans $K$, alors il existe $N \in \mathbb{N}$ tel que $(x_n - x_N)_{n>N}$ est de Cauchy dans $O$, elle est donc convergente dans $O$, et par suite $(x_n)_n$ est convergente dans $K$. 
 
\subparagraph*{$ii$  -} Si $R = \{a\lambda \in O : \lambda \in \Lambda \}$  ; alors pour tout $\lambda \in \Lambda$,
\begin{center}
 $O_{\lambda} = {a \in O : |a - a_{\lambda}| < 1}$
 \end{center} est un ouvert de $O$ et on a : $O = \bigcup\limits_{\lambda \in \Lambda} o_{\lambda} $, d'où, $O$ étant compact,\\
il existe $\lambda_1$, $\lambda_2$, ..., $\lambda_n$, tel que $O = \bigcup\limits_{i=1}^{n} o_{\lambda_i} $ donc $R = \{a_{\lambda_1}, a_{\lambda_2}, ..., a_{\lambda_n} \}$ . Donc le corps résiduel O/P est fini.
\subparagraph*{$iii$ -} Soit $\pi$ l'élément de $R$ vérifiant : $| \pi |= max {| x |: x \in \mathbb{R} et | x |< 1}$.\\
montrons que $P = \pi.O$
\\Soit $x \in \pi.O$, alors il existe $y \in O$ : $x = \pi.y$ et par suite :
\begin{center}
$| x | =| \pi.y | = | \pi | . | y |\leqslant | \pi | < 1$
\end{center}
Donc $x \in P$.\\
Inversement, soit $x \in P$ alors : il existe $a \in \mathbb{R}$ tel que $a - x \in P$ d'où $a \in P$
et par suite $| a | \leqslant | \pi |$ et donc $| x | \leqslant | \pi |$; car si non on aura :
 \begin{center}
$| a | = max (| a - x |, | x |) = | x | > | \pi |$ ce qui est absurde.
\end{center}
{\small $\dfrac{\pi}{x}$} $\in O$ et $x = x$.{\small $\dfrac{\pi}{x}$}
c'est à dire que $x \in \pi.O$. Donc $P = \pi.O$; c'est à dire
que $P$ principal ; donc, d'après la proposition :
\begin{center}
" Si ($K, | . |$) est un corps valué non archimédien alors
$K$ est discret si, et seulement si, $P$ est un idéal principal."  
\end{center}
la valuation de $K$ est discrète.
 
\begin{flushright}
$\blacksquare$
\end{flushright}

\paragraph*{Proposition 2 :}
\subparagraph*{} Soit ($K, |.|$) un corps valué non archimédien vérifiant les trois propriétés suivantes :
\begin{itemize}
\item[$i$   -] Le corps $K$ est complet.
\item[$ii$  -] Le corps résiduel de $K$ est fini.
\item[$iii$ -]La valuation de $K$ est discrète.
\end{itemize}
Alors $K$ est local.
\paragraph*{Démonstration :} 
Montrons que $O$ est compact ; soit $(V_{\lambda})_{\lambda} \in \Lambda $ un recouvrement de $O$ par des ouverts $K$, montrons qu'on peut en extraire un recouvrement fini : par l'absurde, supposons que ce n'est pas le cas.\\
Comme $O = \bigcap\limits_{a\in R} (a + \pi.O)$ où $R$ est fini et $\pi$ est l'élément défini dans 3.5 (le corps résiduel est fini et la valuation est discrète). D'où il existe au moins $a_0 \in R$ tel que $a_0 + \pi.O$ ne peut être recouvert par un nombre fini de $V_{\lambda}$.\\
Ensuite de la même façon il existe $a_1 \in R$ tel que $a_0 + a_1.\pi + \pi^2.O$ ne peut être recouvert par un nombre fini de $V_{\lambda}$. \\Et ainsi de suite il existe une
suite $(a_n)_n$ d'éléments de $R$ telle que pour tout $n \in \mathbb{N}$ $\sum\limits{i=0}^{n-1} a_i \pi^i + \pi^n O$
n'est pas recouvert par un nombre fini de $V_{\lambda}$. Soit l'élément de $O$ défini par : $a = \sum\limits{i=0}^{\infty} a_i^{\pi^i}$ ; alors $\exists \lambda_0 \in \Lambda : a \in V_{\lambda_0}$.\\ Comme $V_{\lambda_0}$ est ouvert, alors
$\exists N \in \mathbb{N}$ : $a + \pi^N.O \subseteq V_{\lambda_0}$. Ce qui est une contradiction parce que nous
avons construit a de sorte qu'aucun des ensembles $a + \pi n.O$, $n \in \mathbb{N}$, n'est couvert par une famille finie de $V_{\lambda_0}$.\\
Donc, pour tout $x \in K$, le sous-ensemble $C = x + O$ est un compact de $K$,car c'est limage directe par l'application continue $y \mapsto x + y$ du compact $O$, de plus $C$ est un voisinage de $x$ car il contient le voisinage ouvert $x+\pi.O$
de $x$. Donc $K$ est localement compact.
\begin{flushright}
$\blacksquare$
\end{flushright}
\paragraph*{Corollaire :}
\subparagraph*{} Soit ($K, |.|$) un corps valué non archimédien $K$ est local si, et seulement si, les trois propriétés suivantes sont satisfaites :
\begin{itemize}
\item[$i$   -] Le corps $K$ est complet.
\item[$ii$  -] Le corps résiduel de $K$ est fini.
\item[$iii$ -] La valuation de $K$ est discrète.
\end{itemize}

\subparagraph*{Remarque :}
\subparagraph*{1 -} On peut maintenant donner une démonstration au premier exemple des corps valués locaux : le corps des nombre p-adiques $\mathbb{Q}_p$ :
\subparagraph*{•} On a vu dans la partie précédente que $\mathbb{Q}_p$ est le complété de $\mathbb{Q}$, donc il est évident que $\mathbb{Q}_p$ est complet.
\subparagraph*{•} Le corps résiduel de $\mathbb{Q}_p$ est $\mathbb{F}_p = \mathbb{Z}/p\mathbb{Z}$ et il est de plus fini.\\
En effet, ce résultat provient de la proposition suivante :
\begin{center}
"L'application naturelle de $\mathbb{Z}/p^n\mathbb{Z}$ dans $\mathbb{Z}_p/p^n\mathbb{Z}_p$ est un isomorphisme."
\end{center} 
Donc il est fini de cardinal $p$.
\subparagraph*{•} La valuation p-adique est une valuation discrète  
$\upsilon_p (\mathbb{Q}^*_p) = \upsilon_p (\mathbb{Q^*}) = \mathbb{Z}$.
\subparagraph*{2 -} Soit ($K, | . |$) un corps local, alors deux cas de figure se présentent, selon la caractéristique de $K$ et sa caractéristique résiduelle :

\subparagraph*{i -} Si la caractéristique de $K$ est différente de sa caractéristique résiduelle, le corps $K$ est nécessairement de caractéristique nulle, et il est isomorphe à une extension finie du corps des nombres $p$-adiques, où $p$ désigne la caractéristique résiduelle de $K$; (qui est un nombre premier, le corps résiduel étant fini).
\subparagraph*{ii -} En cas de caractéristiques égales, le corps $K$ est isomorphe au corps des séries formelles de Laurent à coefficients dans son corps résiduel.

\section*{3. Extensions de corps locaux}
\addcontentsline{toc}{section}{3. Extensions de corps locaux}
\subsection*{3.1. Définitions et exemples}
\addcontentsline{toc}{subsection}{3.1. Définitions et exemples}
\paragraph*{Définitions : }
\subparagraph*{}
Soit $\overline{\mathbb{Q}}_{p}$ une clôture algébrique de $\mathbb{Q}_{p}$. On définit $\mathbb{C}_{p}$ comme étant le complété de $\overline{\mathbb{Q}}_{p}$.
i.e. $\mathbb{C}_{p}=\widehat{\overline{\mathbb{Q}}}_{p}$.
\subparagraph*{}
On rappelle qu'un corps local est un corps complet pour une valuation discrète. En particulier, si $K$ est un corps local muni de la valuation $v$, alors il existe $a \in \mathbb{Z}$ tel que $v(K^{*}) = aZ$. On appelle \textbf{uniformisante} de $K$ un élément $\pi$ de $K$ tel que $v(\pi) = a$.
\paragraph*{Exemples :}
\subparagraph*{$\ast$} $\mathbb{Q}_{p}$ muni de $v_{p}$ est un corps local et $p$ en est une uniformisante.
\subparagraph*{$\ast$} Si $K$ est un corps, alors $K((T))$, corps des séries de Laurent à coefficient dans $K$, muni de la valuation $\upsilon_T$, est un corps local, et $T$ en est une uniformisante.

\paragraph*{Proposition :}
\subparagraph*{} L'idéal maximal $m_K$ de l'anneau $O_K$ de $K$ est principal, et un élément de $K$ en est un générateur si et seulement si c'est une uniformisante.
\subsection*{3.2. Écriture et corps résiduel}
\addcontentsline{toc}{subsection}{3.2. Écriture et corps résiduel}
Nous allons voir que les éléments de corps locaux possèdent tous une écriture particulière, sous forme de série, à partir de leur corps résiduel.
\paragraph*{Théorème de représentation :}
\subparagraph*{} Soit $K$ un corps local, $\pi$ une uniformisante de $K$ et soit $R$ un système de représentants de $K$ contenant $0$. Alors tout $x \in K^*$ s'écrit de façon unique comme une somme :
\begin{center}
$x= \sum\limits_{n\geqslant 0} s_n\pi^n$ $ s_n \in R,s_n \neq 0,n=v(x)$.
\end{center}
En particulier, $n \geqslant 0$ ssi $x \in O$ .
\paragraph*{Démonstration :}
\subparagraph*{} Prenons $x \in O^*$. On sait qu'il existe un unique $s_0 \in R$ non nul tel que 
\begin{center}
 $x = s_0 + \pi x_1$
\end{center} avec $x_1 \in  O$. En itérant cette procédure pour $x_1$, on obtient par récurrence 
\begin{center}
$x = s_0 + s_1\pi + ... + s_{k-1}\pi^{k-1} + x_k\pi^k$
\end{center}
avec $s_0,...,s_{k-1} \in R$ et $x \in O$. On peut écrire alors $x = r_k + x_k\pi^k$ où $r_k$ est la somme partielle
$\sum\limits_{i=0}^{k-1}s_i\pi^i$. Or, comme $|\pi| < 1$, il est clair que $|x_k\pi^k| \leqslant |\pi|^k \rightarrow 0$ et alors la suite $(r_k)$ converge vers $x$ lorsque $k$ tend vers l'infini. Et en notant que la somme commence en $n = 0 = v(x)$, on a le résultat pour $x \in O^*$.
\subparagraph*{} Maintenant, si $x \in K$, il suffit de noter qu'il existe un unique $n = v(x) \in \mathbb{Z}$ tel que $v(\pi^{- n} x) = 0$, donc $\pi^{- n} x \in  O^*$. Par ce qu'on vient de montrer, on voit immédiatement qu'on obtient une
série pour $x$ commençant à l'indice $k = n$.
\begin{flushright}
$\blacksquare$
\end{flushright}
\subsection*{3.3. Ramification et inertie}
\addcontentsline{toc}{subsection}{3.3. Ramification et inertie}
En théorie algébrique des nombres, on parle de ramification d'un idéal premier, lorsque le prolongement de cet idéal à un sur-corps admet au moins un facteur premier ayant une multiplicité plus grande que $1$. Plus précisément, si $L/K$ est une extension finie de corps de nombres et $P$ un idéal premier de $O_K$ (l'anneau des entiers de $K$) alors l'idéal $PO_L$ de $O_L$ se décompose sous la forme :
\begin{center}
$PO_L = P_1^{e(1)}... P_k^{e(k)}$ 
\end{center}
\subparagraph*{} Nous allons maintenant voir que la théorie de la ramification pour les corps locaux (de corps résiduel de caractéristique non nulle) se trouve bien plus simple
que celle des corps de nombres. Dans toute cette partie, sauf mention du contraire, $F$ désigne un corps complet
pour une valuation discrète, dont le corps résiduel $k_F$ est de caractéristique $p$,(i.e. un corps local).


\subsubsection*{3.3.1. Extensions de valuations}
\addcontentsline{toc}{subsubsection}{3.3.1. Extensions de valuations}
\subparagraph*{} Après avoir vu l'équivalence des normes en dimension finie, nous allons en déduire qu'il n'y a qu'une seule manière de prolonger une valuation sur un corps complet à une extension finie de celui-ci.
\paragraph*{Définition :}
\subparagraph*{} Soit $L$ une extension finie d'un corps $K$, et $x \in L$. Alors on définit la norme de $x$ de $L$ sur $K$, $N_{L/K}(x)$ comme étant le déterminant de l'endomorphisme $y \mapsto xy$ du $K$-espace vectoriel $L$. Si $P = X^d + ... + a_0$ est le polynôme minimal unitaire de $x$ sur $K$, alors $N_{L/K}(x) = ((-1)^da_0)^{\dfrac{[L:K]}{d}}$.
\paragraph*{Théorème :}
\subparagraph*{} Soit $K$ un corps complet pour une valuation $\upsilon$, et soit $L$ une extension finie de $K$. Alors il existe une unique manière de prolonger $\upsilon$ en une valuation de $L$. De plus, si $x \in L$, alors :
\begin{center}
$\upsilon(x) = \dfrac{1}{[L : K]}\upsilon(N_{L/K}(x))$.
\end{center}
\paragraph*{Démonstration :}
\subparagraph*{} Considérons d'abord l'unicité : $L$ peut être vu comme un $K$-espace vectoriel de dimension finie $[L : K]$. Si $\upsilon_1$ et $\upsilon_2$ sont deux valuations sur $L$ prolongeant $\upsilon$ sur $K$, alors par l'équivalence des normes en dimension finie et la caractérisation de cette équivalence en terme de valuations : il existe $s \in \mathbb{R}^*_+$ tel que $\upsilon_2(x) = s\upsilon_1(x)$. En prenant $x \in K$, on obtient $s = 1$, ce qu'on voulait démontrer.
\subparagraph*{} Maintenant, pour l'existence, montrons que la formule donnée définit bien une valuation sur $L$. La seule chose non immédiate à vérifier est l'inégalité ultra-triangulaire : $\upsilon(\alpha + \beta) > inf(\upsilon(\alpha), \upsilon(\beta))$. En soustrayant $\upsilon(\alpha)$ ou $\upsilon(\beta)$ de chaque côté
de l'égalité (et en supposant $\alpha$, $\beta \neq 0$, sans ça, il n'y a rien à montrer), on est
ramené à $\upsilon(1 + x) > inf(0, \upsilon(x))$ et $\upsilon(1 + \dfrac{1}{x}) > inf(0, \upsilon(\dfrac{1}{x}))$ Ceci revient alors à montrer que si $\upsilon(x) > 0$, alors $\upsilon(1 + x) > 0$ (en effet, montrer cette implication couvre bien le second cas de $\upsilon(1 + x) > inf(0, \upsilon(x))$, qui se ramène bien l'inégalité ultra-triangulaire).
\subparagraph*{} Soit $x \in L$, supposons que $\upsilon(N_{L/K}(x)) > 0$ et montrons que $\upsilon(N_{L/K}(1 + x)) > 0$. Soit $f(X) = X^d + ... + a_0$ le polynôme minimal de $x$ sur $K$. On a alors (par multiplicativité des degrés) $d|[L : K]$ et $N_{L/K}(x) = ((-1)^da_0)^{\dfrac{[L:K]}{d}}$. Ainsi, $\upsilon(N_{L/K}(x)) > 0$ implique $a_0 \in O_K$ et l'irréductibilité de $f$ implique alors $f$ à coefficients dans $O_K$, d'après le corollaire 1.2.14. De plus, le polynôme minimal
de $1 + x$ est $f(X - 1)$ et donc $N_{L/K}(1 + x) = ((-1)^df(-1))^{\dfrac{[L:K]}{d}} \in O_K$, ce qui
conclut.
\begin{flushright}
$\blacksquare$
\end{flushright}
\subsubsection*{3.3.2. Définitions}
\addcontentsline{toc}{subsubsection}{3.3.2. Définitions}
\subparagraph*{Lemme :}
 Soient $K$ un corps valué complet et $L$ une extension finie de $K$, alors $k_L$ est une extension algébrique de $k_K$ de degré $\leqslant [L : K]$.
\subparagraph*{Indice d'inertie de l'extension :} Si $K$ est une extension finie de $F$, on a vu que $k_K$ est une extension finie de $k_F$. Le degré de $k_K$ sur $k_F$ est l'indice d'inertie de l'extension $K/F$, et sera noté $f = f(K/F)$.
\subparagraph*{Indice de ramification d'une extension :} Si $x \in K$, on a $\upsilon(x) = \dfrac{1}{[L : K]}\upsilon(N_{L/K}(x))$, donc $\upsilon(K^*) \subset \dfrac{1}{[K:F ]} \upsilon(F^*)$ et $\upsilon(F^*)$ est un sous groupe d'indice fini $e = e(K/F)$ de $\upsilon(K^*)$ (en tant que sous-groupes, discrets, de $\mathbb{R}$). $e$ est appelé l'indice de ramification de l'extension $K/F$.
\subparagraph*{Divers type de ramification :} On dit que l'extension $K/F$ est non ramifiée si $e(K/F) = 1$, et si $k_K/k_F$ est séparable. Elle est \textit{totalement ramifiée} si $e(K/F) = [K : F]$. Elle est
\textit{modérément ramifiée} si $e(K/F)$ est premier à la caractéristique du corps résiduel et si $k_K/k_F$ est séparable, et \textit{sauvagement ramifiée} dans le cas contraire.
\paragraph*{Lemme :}
Si $L/K$ et $K/F$ sont deux extensions finies, alors $e(L/F ) = e(L/K)e(K/F )$ et $f(L/F ) = f(L/K)f(K/F )$.
\subparagraph*{} En effet, pour la seconde égalité, c'est simplement la multiplicativité des degrés. Pour la première, la multiplicativité des degrés, la propriété : \begin{center}
"$e(K/F)f(K/F) = [K : F]$"
\end{center}
et la seconde égalité permettent de conclure.

\subsection*{3.4. Quelques structures des extensions}
\addcontentsline{toc}{subsection}{3.4. Quelques structures des extensions}
\subsubsection*{3.4.1. Extensions totalement ramifiées}
\addcontentsline{toc}{subsection}{3.4.1 Extensions totalement ramifiées}
\paragraph*{Définition (Polynôme d'Eisenstein) : }
\subparagraph*{} Si $K$ est un corps local, on appelle polynôme d'Eisenstein de degré $d \in \mathbb{N}$ un polynôme $P (X) = X^d + a_{d-1}X^{d-1} + ... + a_0 \in K[X]$ tel que $a_0$ soit une uniformisante de $K$ et $a_1,..., a_{d-1} \in m_K$.
\paragraph*{Remarque :} Un polynôme d'Eisenstein est irréductible (c'est une forme du critère d'Eisenstein).
\paragraph*{Lemme :}
 Soit $K$ un corps complet pour une valuation discrète, et $P$ un polynôme d'Eisenstein de degré $d$. Soit $L = K[X]/P$ et $x$ l'image de $X$ dans $L$. alors :
 \begin{itemize}
 \item[$(i)$]Si $\upsilon(K^*) = uZ$ avec $u > 0$, alors $\upsilon(x) = \dfrac{u}{d}$ et $L/K$ est totalement ramifiée.
 \item[$(ii)$]$x$ est une uniformisante de $L$.
 \item[$(ii)$]Si $y = \sum\limits_{i=0}^{d-1} a_ix^i$ , avec les $a_i$ dans $K$, alors $\upsilon(y) = inf_i(\upsilon(a_i) + i\dfrac{u}{d})$
 \end{itemize}
 \paragraph*{Démonstration :}
\subparagraph*{} Pour le $(i)$, on a la formule $\upsilon(x) = \dfrac{1}{[L : K]}\upsilon(N_{L/K}(x))$ et $P$ étant
d'Eisenstein, $a_0$ est une uniformisante, donc $\upsilon(N_{K/F} (x)) = \left(\dfrac{[L:K]}{u}\right) u = u$. Le $(i)$ donne alors le $(ii)$, comme on connait $\upsilon(L^*)$. Pour le $(iii)$, cela vient directement du fait que chaque terme a une valuation différente. Pour voir le fait qu'elles soient différentes, il suffit de constater que $\upsilon(a_i) + i\dfrac{u}{d} = u(\upsilon_i + \dfrac{i}{d} )$ avec $\upsilon_i \in \mathbb{Z}$, d'après la définition de $u$, et tout les $\upsilon_i + \dfrac{i}{d}$ ont une partie fractionnaire distincte (comme $i \leqslant d$). 
\begin{flushright}
$\blacksquare$
\end{flushright}
\subsubsection*{3.4.2. Extensions non ramifiées}
\addcontentsline{toc}{subsubsection}{3.4.2. Extensions non ramifiées}
\paragraph*{Théorème :}
\subparagraph*{} Soit $k$ une extension finie séparable de $k_F$, alors il existe une extension non ramifiée $F (k)$ de $F$ dont le corps résiduel est $k$.\\
Si $L$ est une extension finie de $F$, et si $k_L/k_F$ est une extension séparable, alors $F (k_L) \subset L$, l'extension $L/F (k_L)$ est totalement ramifiée, et $F (k_L)$ est l'unique extension non ramifiée de $F$ ayant ces deux propriétés.
 
\subsubsection*{3.4.3. Extensions galoisiennes}
\addcontentsline{toc}{subsubsection}{3.4.3. Extensions galoisiennes}
\subparagraph*{} Galois -et beaucoup de chercheurs après lui- a créé un outil essentiel pour l'analyse des extensions $L$ d'un corps $K$, il correspond aux automorphismes de $L$ laissant le corps $K$ invariant. L'ensemble de ces automorphismes, munis de la loi interne de composition des applications, forme un groupe. Cet outil est particulièrement efficace dans le cas des extensions finies, par exemple sur le corps des rationnels dans le cas d'un corps de décomposition. Un élément de ce groupe restreint à un ensemble de racines d'un polynôme correspond à une permutation de cet ensemble de racines. Dans le cas des extensions finies, il correspond à un groupe fini appelé groupe de Galois.
\subparagraph*{} Pour que cet outil soit pleinement pertinent, il faut en fait que les polynômes minimaux de l'extension n'aient pas de racines multiples. Ce qui est toujours le cas pour des extensions sur les corps des rationnels ou plus généralement sur un corps de caractéristique nulle. Dans ce cadre, il est par exemple possible de montrer qu'il existe un élément $a$ dit primitif tel que l'extension soit une extension simple égale à $K(a)$. Il faut de plus que l'extension contienne suffisamment de racines. Il faut en fait que le cardinal du groupe soit égal à la dimension de l'extension. Si ces deux hypothèses sont vérifiées, on parle alors d'extension de Galois.
\subparagraph*{} Le groupe de Galois permet de comprendre finement la structure de l'extension. Par exemple, il existe une bijection entre ses sous-groupes et les sous-corps de l'extension. Il est utilisé pour la détermination des polygones constructibles à la règle et au compas ou pour le théorème d'Abel sur la résolution d'équations polynômiales par radicaux.

\paragraph*{Définitions :}
\subparagraph*{} Soient $K$ un corps, $L$ une extension algébrique de $K$, et $\Omega$ la clôture algébrique de $K$. $L$ est identifié à un sous-corps de $\Omega$, ce qui, dans le cas où l'extension est finie, ne nuit en rien à la généralité.

\subparagraph*{Extension normale :} L'extension est dite normale si tout morphisme de $L$ dans $\Omega$ laissant $K$ invariant est un automorphisme de $L$.
\subparagraph*{Remarque :} Un morphisme de corps est toujours injectif. Le morphisme est aussi un morphisme d'espaces vectoriels car $L$ dispose d'une structure d'espace vectoriel sur $K$. Donc, si $L$ est une extension finie, alors il suffit que le morphisme ait une image incluse dans $L$ pour qu'un argument de dimension prouve la surjectivité.
\subparagraph*{Extension galoisienne :} L'extension est dite de \textsc{Galois} ou galoisienne si elle est \textit{normale} et \textit{séparable}.
\subparagraph*{Remarque :} L'extension $L$ de $K$ est dite séparable si le polynôme minimal sur $K$ de tout élément de $L$ n'a aucune racine multiple dans $\Omega$. Si $K$ est un corps parfait (par exemple s'il est de caractéristique $0$ ou si c'est un corps fini), alors $L$ est toujours séparable sur $K$.
\subparagraph*{Groupe de Galois :} L'ensemble des automorphismes de $L$ qui laissent $K$ invariant, muni de la loi de composition des applications, forme un groupe appelé le groupe de \textsc{Galois} de l'extension et est souvent noté $Gal(L/K)$.

\paragraph*{Exemples :}
\subparagraph*{$\ast$} Le corps des nombres complexes est une extension de Galois du corps des nombres réels. C'est une extension simple (c'est-à-dire engendrée par le corps des nombres réels et un seul élément supplémentaire) dont le groupe de \textsc{Galois} est le groupe cyclique d'ordre $2$.
\paragraph*{Démonstration :}
Le corps des nombres complexes est une extension de dimension deux sur le corps des nombres réels. C'est donc une extension simple engendrée par l'imaginaire pur $i$. Comme le corps des nombres complexes est algébriquement clos, tout morphisme de ce corps vers une extension quelconque laissant les nombres réels invariants est un automorphisme. Soit $\sigma$ un automorphisme différent de l'identité. Un automorphisme permute les racines d'un polynôme à coefficients dans le corps de base. Donc $\sigma(i)$ est une racine du polynôme $X^2+1$ au même titre que $i$. Et $\sigma(i)$ est différent de $i$ car sinon $\sigma$ serait l'identité sur une base des nombres complexes: $(1, i)$ et serait donc égal à l'identité. Le polynôme précédent n'admet que deux racines car il est de degré deux. On vérifie que les deux racines sont $i$ et $-i$. L'image de la base $(1, i)$ par $\sigma$ est donc $(1, -i)$. Cela signifie que $\sigma$ est l'application conjuguée. Il est simple de vérifier que $\sigma$ est effectivement un automorphisme d'ordre deux. Le groupe de \textsc{Galois} est donc bien un groupe à deux éléments isomorphe au groupe cyclique d'ordre deux.
\begin{flushright}
$\blacksquare$
\end{flushright}
\subparagraph*{$\ast$} L'extension simple engendrée par la racine cubique de deux sur le corps des rationnels n'est pas une extension de \textsc{Galois}. En effet, ce corps ne contient pas toutes les racines, il existe donc un morphisme de $L$ dont l'image n'est pas $L$.
\paragraph*{Démonstration :}
Il est aisé de vérifier que $L$ possède pour base la famille des trois éléments un, la racine cubique de deux et la racine cubique de quatre. La démonstration est donnée par la proposition :
\begin{center}
"Si $L$ est une extension de $K$ et $l$ un élément de $L$ algébrique sur $K$, le $K$-espace vectoriel $K(l)$ a pour base $(1, l, l^2,... , l^{n-1})$, où $n$ désigne le degré du polynôme minimal de $l$."
\end{center}
appliquée au polynôme $P(X)=X^3-2$. Or si l'on note
 $r = j. \sqrt[3]{2}$, où $j$ est une racine cubique complexe de l'unité, $r$ est aussi une racine du polynôme. On peut
garantit l'existence d'un morphisme de $L$ dans l'extension de $r$ sur les nombres rationnels. L'extension $L$ n'est donc pas normale, ce n'est pas une extension galoisienne.
\begin{flushright}
$\blacksquare$
\end{flushright}
\section*{4. Corps locaux non commutatifs}
\addcontentsline{toc}{section}{4. Corps locaux non commutatifs}
\subsection*{4.1. Définitions }
\addcontentsline{toc}{subsection}{4.1. Définitions }
\paragraph*{Le centre d'un corps :} Soit  $( K , + , . )$ , Le centre de $( K , + , . )$ est le sous-ensemble de $K$ formé par tous les éléments $x$ de $K$ tels que $x . r = r . x$ pour tout $r$ de $K$. Le centre de $K$ est un sous-corps commutatif de $K$ .
\paragraph*{Le centre d'une algèbre :} Le centre d'une algèbre $E$ est constitué de tous les éléments $x$ de $E$ tels que $x . a = a . x$ pour tout $a$ de $E$.
\paragraph*{Algèbre de quaternions sur un corps :}On appelle algèbre de quaternions sur le corps $K$ toute algèbre (unitaire et associative) $A$ de dimension 4 sur $K$ qui est simple (c'est-à-dire que $A$ et $\{0\}$ sont les seuls idéaux bilatères) et dont le centre est $K$.

\paragraph*{Corps de quaternions :} On appelle corps de quaternions sur $K$ toute algèbre de quaternions sur $K$ dont l'anneau sous-jacent est un corps.
\paragraph*{Exemple :} Le seul corps de quaternions sur $\mathbb{R}$ (à isomorphisme près) est le corps $\mathbb{H}$ des quaternions de \textsc{Hamilton}.
\subsection*{4.2. Structures des corps locaux non commutatifs}
\addcontentsline{toc}{subsection}{4.2. Structures des corps locaux non commutatifs}
\subparagraph*{} Il est possible d'élargir la définition équivalente ci-dessus d'un corps local non archimédien en autorisant les corps non commutatifs. Les notions d'anneau d'entiers, de corps résiduel et de caractéristique résiduelle s'étendent sans difficulté à ce cadre.
\subparagraph*{} Le centre d'un corps local non commutatif est un corps local.
\subparagraph*{} Un corps local non commutatif est localement compact si et seulement si son centre l'est, et si et seulement s'il a un corps résiduel fini. Dans ce cas, il est de dimension finie sur son centre.

\subparagraph*{} Inversement, les corps non commutatifs localement compacts non discrets sont classés comme suit 
:
\begin{itemize}
\item[$\ast$] Si leur valeur absolue est archimédienne, ils sont isomorphes au corps des quaternions de Hamilton ;
\item[$\ast$] Sinon, ce sont des algèbres cycliques sur leur centre, paramétrées par un élément de $\mathbb{Q}/\mathbb{Z}$.
\end{itemize}

\newpage
\vspace*{\stretch{1}}
\begin{center}
\includegraphics[scale=0.8]{2.png}  
\end{center}
\begin{center}

{\LARGE {\huge Troisième partie :}}\\
\textit{   }\\
\textit{   }\\
{\LARGE {\Huge Applications des corps valués locaux à la théorie des nombres}}
\addcontentsline{toc}{part}{Troisième partie : Applications des corps valués locaux}\\
\textit{  }
\end{center}
\begin{center}
\includegraphics[scale=1]{3.png}  
\end{center}
\vspace*{\stretch{1}}
\newpage

\paragraph*{Introduction.}
Il est toujours bien, dans l'abstraction et l'idéalité des maths, de donner des applications et des exemples concrets d'utilisation de notions mathématiques avancées, parce que, quoique cette abstraction soit suffisante en elle-même pour le cercle des adeptes de mathématiques, les non mathématiciens ne partagent généralement pas cet avis. On présente donc, dans cette dernière partie, quelques applications -mathématiques- des corps valués locaux, et surtout du corps local $\mathbb{Q}_{p}$, qui est la grande nouveauté des corps valués.\\
Les corps valués locaux ont plusieurs applications en cryptographie, et particulièrement dans l'analyse des registres à décalage à rétroaction linéaire. On les retrouve également dans plusieurs preuves de théorèmes et lemmes en théorie des nombres, mais également dans les autres secteurs des mathématiques, on citera en l'occurrence le théorème de Monsky en géométrie carré triangles, qui affirme qu'il est impossible de partitionner un carré en un nombre impair de triangles de même aire, et dont la démonstration utilise le corps $\mathbb{Q}_{p}$ et surtout la valuation $2$-adique. Des applications de cette notion sont également retrouvées en physique, et notamment dans la fameuse théorie des cordes.\\
Dans cette partie, nous allons parler de polynômes de Newton, qui est une théorie assez avancée dans l'étude des corps valués locaux, et qui convient en continuité avec la partie précédente, pour en déduire un critère d'irréductibilité purement 'arithmétique' des polynômes dans $\mathbb{Q}[X]$ et $\mathbb{Z}[X]$. On discutera ensuite des applications des $p$-adiques aux équations diophantiennes et à la théorie des nombres élémentaires. On essaiera d'aller du plus compliqué vers le plus simple, terminant ainsi notre projet par les nombres $p$-adiques expliqués à un lycéen. 

\section*{1. Critère d'Eisenstein et polygones de Newton}
\addcontentsline{toc}{section}{1. Critère d'Eisenstein et polygones de Newton}
La théorie des polygones de Newton est un outil très puissant pour l'étude des polynômes, qui nous permettra, en l'occurrence, de donner une généralisation du critère d'Eisenstein. Un polygone de Newton est, informellement, un polygone du plan euclidien que l'on peut associer à un polynôme, lorsque les coefficients de ce dernier sont des éléments d'un corps valué. Il est particulièrement utile lorsqu'il s'agit d'un corps local non archimédien, i.e. un corps des nombres $p$-adiques ou un corps de séries de Laurent sur un corps fini.\\
Pour ce qui suit, on considère $K$ un corps valué muni de la valuation $v$ et on fixe $\overline{K}$ une clôture algébrique de $K$.       
\subsection*{1.1. Construction d'un polygone de Newton}
\addcontentsline{toc}{subsection}{1.1. Construction d'un polygone de Newton}
Soit $P(X) \in K[X]$ un polynôme. On peut supposer, sans perte de généralité, que $P(0)=1$ et écrire $P(X)$ sous la forme :\\
$P(X)=1+a_{1}X+...+ a_{n}X^{n}$ où $a_{1},...,a_{n} \in K$ et $a_{n}\neq 0$. On considère l'ensemble $S$ des points du plan définis par :\\
 $A_{0}=(0,0)$ et $A_{i}=(i,v(a_{i}))$ pour $1\leq i\leq n$, en ignorant les indices $i$ pour lesquels $a_{i}=0$. Le polygone de Newton de $P$ est alors la frontière inférieure de l'enveloppe convexe de cet ensemble de points. Rappelons que l'enveloppe convexe d'un ensemble est l'ensemble convexe le plus petit parmi ceux qui le contiennent.\\
Une construction plus explicite consiste à considérer l'axe des ordonnées, et le faire tourner autour de l'origine $A_{0}$ dans le sens inverse des aiguilles d'une montre, jusqu'à ce qu'il rencontre l'un des points $A_{i_{1}}$ de $S$ ; on obtient alors le premier segment $[A_{i_{0}},A_{i_{1}}]$ du polygone de Newton. Si l'on continue à faire tourner l'axe, autour du point $A_{i_{1}}$ cette fois, il finit par rencontrer un point  $A_{i_{2}}$, et on obtient ainsi le second segment $[A_{i_{1}},A_{i_{2}}]$. En répétant cette opération autant de fois que possible, on finit par obtenir le polygone de Newton.\\

\begin{center}
\includegraphics[scale=0.5]{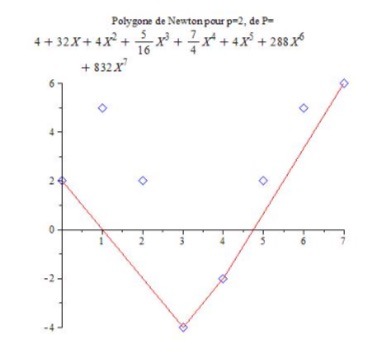}  
\end{center}   

\subsection*{1.2. Quelques résultats dans les corps locaux non archimédiens}
\addcontentsline{toc}{subsection}{1.2. Quelques résultats dans les corps locaux non archimédiens}
Dans ce qui suit, $(K, v)$ est un corps local non archimédien ($\mathbb{Q}_{p}$ ou $\mathbb{F}_{q}((X))$).
\paragraph*{1.2.1. Définition (Polynôme pur) : } 
\subparagraph*{}
On dit qu'un polynôme $P$ est \textbf{pur} de pente $m$ si son polygone de Newton est un unique segment dont la pente est $m$. Dans ce cas, on a nécessairement $ m=\dfrac{v(a_{n})}{n}$.  
\paragraph*{1.2.2. Théorème (Factorisation d'un polynôme) :  }
\subparagraph*{}
Le polynôme $P$ admet une factorisation de la forme :
 \begin{center}
 $P(X)=Q_{1}(X)...Q_{r}(X)$
 \end{center}
Où chaque $Q_{i}$ est un polynôme de $K[X]$ de degré $l_{s}$ et de pente $m_{s}$.   
\paragraph*{Conséquence :}
Si $P$ est irréductible dans $K[X]$, alors il est pur de pente $\dfrac{v(a_{n})}{n}$. 

\subsection*{1.3. Cas de $\mathbb{Q}_{p}$ : Critère de Dumas et critère d'Eisenstein}
\addcontentsline{toc}{subsection}{1.3. Cas de $\mathbb{Q}_{p}$ : Critère de Dumas et critère d'Eisenstein}
On considère maintenant $K$ comme étant un corps de nombres $p$-adiques $\mathbb{Q}_{p}$, $v$ coïncide donc avec la valuation $p$-adique $v_{p}$. \\
En utilisant ce qui précède et le théorème donnant la localisation des racines dans ce cas, on retrouve un critère d'irréductibilité, dit \textbf{critère de Dumas} : 
\paragraph*{Critère de Dumas :}
Si le polynôme $P$ est pur, et $v_{p}(a_{n})$ est premier avec $n$, alors $P$ est irréductible dans $\mathbb{Q}_{p}[X]$. 
\subparagraph*{}
Ce critère peut être démontré par contraposée. La réciproque est fausse ($P(X)=1+X+X^{2}$). Ce qui nous intéresse est le théorème suivant, qui est un cas particulier du critère de Dumas, appelé \textbf{critère d'Eisenstein}.
\paragraph*{Critère d'Eisenstein :} 
Si $P(X)=a_{0}+...+a_{n-1}X^{n-1}+X^{n}$ vérifie $v(a_{i})\geq 1$ pour tout $0\leq i\leq n-1$ et $v(a_{0})=1$, alors $P$ est irréductible.
\subparagraph*{} 
On voit toujours que cet énoncé utilise toujours les valuations et n'est probablement pas accessible aux lycéens ou à ceux qui ne sont pas familiers avec le corps $\mathbb{Q}_{p}$, mais le critère peut être énoncé de la façon suivante, qui est sa forme la plus connue d'ailleurs : 
\paragraph*{Critère d'Eisenstein :} 
\subparagraph*{}
On considère un polynôme à coefficients entiers qu'on écrit comme suit : 
\begin{center}
$P(X)=a_{0}+...+a_{n-1}X^{n-1}+a_{n}X^{n}$
\end{center}
Si on suppose qu'il existe un nombre premier $p$ tel que :\\
\\
* $p$ divise $a_{0},...,a_{n-1}$\\
* $p$ ne divise pas $a_{n}$\\
* $p^{2}$ ne divise pas $a_{0}$
\subparagraph*{}
Alors $P(X)$ est irréductible dans $\mathbb{Q}[X]$. Si de plus les $a_{i}$ sont premiers entre eux, alors $P(X)$ est irréductible dans $\mathbb{Z}[X]$.
\begin{flushright}
$\blacksquare$
\end{flushright}
\subparagraph*{}
Ce dernier énoncé n'est pas vraiment difficile à appliquer, on le retrouve par exemple dans le livre '\textit{Olympiades de mathématiques, réflexes et stratégies}' du marocain \textsc{Tarik} \textsc{Belhaj} \textsc{Soulami}, qui est principalement destiné à des élèves de première et de terminale.  
\subparagraph*{ Exemple :} $X^{n}-2$ est irréductible dans $\mathbb{Q}[X]$ pour tout entier $n\geq 1$ (critère d'Eisenstein pour $p=2$), ce qui prouve qu'il y a dans $\mathbb{Q}[X]$ des irréductibles de tout degré.

\section*{2. Nombres p-adiques et équations diophantiennes}
\addcontentsline{toc}{section}{2. Nombres p-adiques et équations diophantiennes}
Malgré l'origine théorique des nombres $p$-adiques et leur construction presque indépendante de la théorie des nombres élémentaire, ils jouent très bien leur rôle dans le monde de l'arithmétique, et spécialement dans une théorie centrale de cette branche : la théorie des \textbf{équations diophantiennes}. La résolution de ces équations a même constitué la motivation de \textsc{Hensel} pour définir les nombres $p$-adiques. Une équation diophantienne, nommée après \textsc{Diophante} quoique existant depuis l'antiquité, est un problème du type : 
\begin{center}
   $F(x_{1},...,x_{n})=0$ 
   \end{center}   
où $F$ est un polynôme à une ou plusieurs variables $x_{1},...,x_{n}$ et où on recherche les solutions entières. Ce problème, très difficile, peut être attaqué en affaiblissant l'équation en considérant l'ensemble des congruences pour tout $m$ : 
\begin{center}
      $F(x_{1},...,x_{n})\equiv 0\mod m $
     \end{center}     
Ce qui équivaut, grâce au \textit{Théorème des restes chinois}, de considérer l'ensemble des congruences :  
\begin{center}
      $F(x_{1},...,x_{n})\equiv 0\mod p^{v} $
     \end{center} 
Pour toutes les puissances de nombres premiers. On espère alors que l'existence ou la non-existence de solutions pour cette dernière équation pourrait nous donner des indices à propos des solutions de l'équation d'origine. Ce sont les nombres $p$-adiques qui nous donnent la réponse.  
\paragraph*{Théorème :}
Soient $F(x_{1},...,x_{n})$ un polynôme à coefficients entiers et $p$ un nombre premier fixé. L'équation 
\begin{center}
$F(x_{1},...,x_{n})\equiv 0\mod p^{v} $
\end{center} 
est résoluble (admet des solutions) en entiers pour $ v\geq 1 $ quelconque si et seulement si l'équation  
 \begin{center}
   $F(x_{1},...,x_{n})=0$ 
   \end{center}
 est résoluble en entiers $p$-adiques (dans $\mathbb{Z}_{p}$).   
\paragraph*{Équivalence (Somme de deux carrés de Fermat) :}
Soit $n$ un entier. L'équation $n=x^{2}+y^{2}$ admet une solution dans $\mathbb{Z}$ si et seulement si, pour tout $p$ premier, elle a une solution dans $\mathbb{Z}_{p}$.\\ 
Ce qui donnerait probablement le fameux \textit{théorème des deux carrés de Fermat}, qui affirme que l'équation ci-dessus admet des solutions si et seulement si pour tout $p$ figurant dans la décomposition en facteurs premiers de $n$ : $p\equiv 3\mod 4 \Leftrightarrow v_{p}(n)$ est pair. 
\paragraph*{Somme de trois carrés :} 
Soit $n$ un entier. L'équation $n=x^{2}+y^{2}+z^{2}$ admet une solution dans $\mathbb{Z}$ si et seulement si, pour tout $p$ premier, elle a une solution dans $\mathbb{Z}_{p}$.\\         
Ce qui se ramène à $n>0$ et $n\neq 4^{r}(8m+7)$.
\subparagraph*{}
On voit donc que les $p$-adiques sont fortement reliés aux équations diophantiennes, qui constituent un sujet de recherche toujours actif pour les mathématiciens contemporains. Mieux analyser ce lien nécessite des connaissances assez avancées telles que la théorie du corps de classes local, les groupes de Brauer et le principe local-global notamment. Notre objectif dans cette partie étant de simplifier les choses, on laisse le soin au lecteur intéressé d'approfondir et explorer le sujet.

\section*{3. Valuations p-adiques et théorie des nombres élémentaire}
\addcontentsline{toc}{section}{3. Valuations p-adiques et théorie des nombres élémentaire}
Dans cette sous-partie, on présente les valuations $p$-adiques d'entiers naturels, ce sont des définitions, théorèmes et applications accessibles à toute personne en connaissance des notions de base des entiers naturels et d'arithmétique, en l'occurrence les élèves de lycée. \\
Dans toute la suite, on considère $p$ un nombre premier. 
\subsection*{3.1. Définitions et théorèmes}
\addcontentsline{toc}{subsection}{3.1. Définitions et théorèmes}
\paragraph*{Définitions : }
Soit $n\geq 2$ un entier. On appelle \textbf{valuation $p$-adique} de $n$ la puissance de $p$ qui apparaît dans la décomposition de $n$ en facteurs premiers. On la note $v_{p}(n)$. \\
$v_{p}(n)$ est donc le plus grand entier naturel $k$ tel que $p^{k}\mid n$ (donc $p^{k+1}\nmid n$).\\
On définit la valuation $p$-adique d'une fraction comme suit : \\
$v_{p}(\dfrac{a}{b})=v_{p}(a)-v_{p}(b)$. 
\paragraph*{Théorème 1 :}
\begin{center}
  $v_{p}(ab)=v_{p}(a)+v_{p}(b)$.
 \end{center} 
\paragraph*{Preuve :} 
Posons $v_{p}(a)=\alpha$ et $v_{p}(b)=\beta$, on a alors : $ a=p^{\alpha}a_{1} $ et $ b=p^{\beta}b_{1}$ avec $a_{1}$ et $b_{1}$ premiers avec $p$.\\
 $ab=p^{\alpha +\beta}a_{1}b_{1}\Longrightarrow v_{p}(ab)=\alpha +\beta=v_{p}(a)+v_{p}(b) $
\paragraph*{Théorème 2 :}
 Si $v_{p}(a)>v_{p}(b)$, on a $v_{p}(a+b)=v_{p}(b)$. 
 \paragraph*{Exemple d'application : }
\begin{center}
Montrer que $\sum\limits_{{i=1}}^n \dfrac{1}{i}$ n'est pas un entier pour $n\geq 2$.
\end{center}
\paragraph*{Solution.} 
Remarquons que $\sum\limits_{{i=1}}^n \dfrac{1}{i} = \dfrac{1}{n!}\sum\limits_{{i=1}}^n\dfrac{n!}{i}$.\\
En utilisant le théorème 2, on trouve que : $v_{2}(\dfrac{n!}{2i-1}+\dfrac{n!}{2i})=v_{2}(\dfrac{n!}{2i})$ puis que $v_{2}(\dfrac{n!}{4i-2}+\dfrac{n!}{4i})=v_{2}(\dfrac{n!}{4i})$, on itère pour finalement trouver l'égalité : 
\begin{center}
 $v_{2}(\sum\limits_{{i=1}}^n \dfrac{n!}{i})=v_{2}(\dfrac{n!}{2^{\lfloor\log_{2}n\rfloor}})$
\end{center}
On suppose maintenant -par l'absurde- que $\sum\limits_{{i=1}}^n \dfrac{1}{i}$ est un entier, donc \begin{center}
$v_{2}(\sum\limits_{{i=1}}^n \dfrac{1}{i}) = v_{2}(\dfrac{1}{n!}\sum\limits_{{i=1}}^n\dfrac{n!}{i}) \geq 0$
\end{center}
D'où : $v_{2}(\dfrac{n!}{2^{\lfloor\log_{2}n\rfloor}})\geq v_{2}(n!)$, ce qui implique que  $\lfloor\log_{2}n\rfloor \leq 0$, absurde puisque $n\geq 2$.
\paragraph*{Autre exemple : }
On procède presque de la même manière pour montrer que $\sum\limits_{{i=0}}^n \dfrac{1}{2i+1}$ n'est pas un entier pour $n\geq 1$.

\subsection*{3.2. L'art d'utiliser les valuations élémentairement }
\addcontentsline{toc}{subsection}{3.2. L'art d'utiliser les valuations élémentairement }
\paragraph*{Exemple 1 :}
Soient $d$ et $n$ deux entiers $\geq 1$. On suppose que $ d^{2} \mid n^{2} $, montrer que $d \mid n$.
\paragraph*{Solution.} Il s'agit de montrer que $v_{p}(d)\leq v_{p}(n)$ pour tout nombre premier $p$. \\
On a : $v_{p}(d^{2})\leq v_{p}(n^{2}) \Longrightarrow 2v_{p}(d)\leq 2v_{p}(n) \Longrightarrow v_{p}(d)\leq v_{p}(n)$.

\paragraph*{Exemple 2 :} 
Soient $a, b, c \geq 1$ trois entiers tels que 
\begin{center}
$\dfrac{1}{a}=\dfrac{1}{b}+\dfrac{1}{c}$
\end{center}
On note $d$ le $pgcd$ des entiers $a, b, c$. Montrer que les entiers $abcd$, $d(c-a)$ et $d(b-a)$ sont des carrés parfaits.    
\paragraph*{Solution.}
L'équation étant homogène en $a, b, c$ et le problème invariant par multiplication de $a, b, c$ par un même entier non nul, on peut supposer que $a, b, c$ sont premiers entre eux. Par symétrie de $b$ et $c$, il s'agit finalement de montrer que $abc$ et $c-a$ sont des carrés parfaits.\\
  L'équation du problème devient $bc=ac+ab$ puis $b(c-a)=ac$. On en déduit que $abc=b^{2}(c-a)$, ce qui réduit le problème à montrer que $c-a$ est un carré parfait donc que toutes ses valuations $p$-adiques sont paires.\\
    Soit $p$ un premier divisant $c-a$, donc il divise $b(c-a)=ac$, ce qui implique que $p\mid a$ ou $ p\mid c$. Puisque divisant leur différence, $p$ divise les deux. $a, b, c$ sont premiers entre eux, donc $p$ ne divise pas $b$, d'où : 
    \begin{center}
        $v_{p}(c-a)=v_{p}(b(c-a))=v_{p}(ac)=v_{p}(a)+v_{p}(c)=\alpha + \gamma$
        \end{center}    
En particulier, $v_{p}(c-a)\geq \alpha =v_{p}(a)$ donc $v_{p}(c-a+a)\geq\alpha$ d'où $\gamma\geq\alpha$. Par symétrie, $v_{p}(c-a)=2\alpha$ pair. CQFD.  

\paragraph*{Exemple 3 :}
Soient $a, b, c \geq 1$ trois entiers tels que la somme $\dfrac{a}{b}+\dfrac{b}{c}+\dfrac{c}{a}$ soit entière. Montrer que le produit $abc$ est un cube.  

\paragraph*{Solution.} On procède de la même façon que pour l'exemple 2 pour supposer que $a, b, c$ sont premiers entre eux.\\
 Soit $p$ un diviseur de $a$. Puisque $p\mid abc\mid a^{2}b+b^{2}c+c^{2}a$, on a $p\mid b^{2}c $, d'où ou bien $p\mid b$ ou bien $p\mid c$. \\
 On suppose donc que $p$ divise $a$ et $b$ (pas $c$) et on montre par l'absurde que $v_{p}(b)=2v_{p}(a)$ (pas difficile) ce qui nous permet de conclure que $v_{p}(abc)=3v_{p}(a)$, d'où $abc$ est un cube.

\subsection*{3.3. Valuations p-adiques des factorielles }
\addcontentsline{toc}{subsection}{3.3. Valuations p-adiques des factoriels }
La factorielle et la valuation $p$-adique sont de très bons amis. En effet, \textsc{Legendre} avait découvert un lien entre les valuations et les nombres factoriels des entiers, donnant ainsi, grâce à la formule qui porte son nom, des calculs assez simples de la valuation d'une factorielle. 

\paragraph*{Formule de Legendre :}
\subparagraph*{} Pour tout entier naturel $n$ et $p$ premier, on a : 
\begin{center}
$v_{p}(n!)=\sum\limits_{{i=1}}^\infty \lfloor \dfrac{n}{p^{i}}\rfloor$
\end{center}
\paragraph*{Ébauche de démonstration :} 
Une idée de preuve repose sur la combinatoire en considérant l'ensemble $ \lbrace 1, 2, ..., n \rbrace $ et en comptant le nombre de fois où une puissance donnée de $p$ apparaît dans cet ensemble.  
\paragraph*{Remarque :} Quoique la somme ci-dessus ait l'air infinie, elle est bel et bien finie. En effet, $(\lfloor \dfrac{n}{p^{i}}\rfloor)_{i}$ est nécessairement nulle à partir d'un certain $i$.  
\paragraph*{Autre formule de la valuation d'une factorielle :}
\subparagraph*{} Pour tout entier naturel $n$ et $p$ un nombre premier, on a : 
\begin{center}
$v_{p}(n!)=\dfrac{n-s_{p}(n)}{p-1}$
\end{center}
Où $s_{p}(n)$ représente la somme des chiffres de $n$ écrit en base $p$.
\paragraph*{Preuve :}
En base $p$, $n$ s'écrit sous la forme : 
\begin{center}
$n=\sum\limits_{{i=0}}^k a_{i}p^{i}$ avec $ 1\leq a_{k}\leq p-1 $ et $0\leq a_{i}\leq p-1$ pour $0\leq i\leq k-1$
\end{center}
On utilise la formule précédente pour trouver : 
\begin{center}
$v_{p}(n!)=\sum\limits_{{i=1}}^\infty \lfloor \dfrac{n}{p^{i}}\rfloor = \sum\limits_{{j=1}}^k a_{j}(p^{j-1}+p^{j-2}+...+1)= \sum\limits_{{j=1}}^k a_{j}(\dfrac{p^{j}-1}{p-1}) $
\end{center}
On évalue maintenant le terme $ \dfrac{n-s_{p}(n)}{p-1} $. On sait que $s_{p}(n)=\sum\limits_{{i=0}}^k a_{i}$, ce qui implique que : 
\begin{center}
$ \dfrac{n-s_{p}(n)}{p-1} = \dfrac{1}{p-1} (\sum\limits_{{i=0}}^k a_{i}p^{i} - \sum\limits_{{i=0}}^k a_{i})= \dfrac{1}{p-1} \sum\limits_{{i=0}}^k a_{i}(p^{i}-1) $
\end{center}
D'où l'égalité. 
\paragraph*{Applications :}
\subparagraph*{1. } Trouver tous les entiers naturels $n$ vérifiant $2^{n-1}\mid n!$
\subparagraph*{} On a $v_{2}(n!)=n-s_{2}(n)\geq n-1 \Longrightarrow s_{2}(n)\leq 1$, ce qui veut dire que $n$ est une puissance de $2$.
\subparagraph*{2. } On note $\pi(x)$ le nombre des nombres premiers au plus égaux à $x$ et $P$ l'ensemble des nombres premiers. 
\begin{itemize}
 \item[(a)] Montrer que $\binom{2n}{n}$ divise $\prod\limits_{{p\in P, p\leq 2n}}{p^{\lfloor\dfrac{\ln (2n)}{\ln p}\rfloor}}$
 \item[(b)] Montrer que $\binom{2n}{n}\leq (2n)^{\pi(2n)}$
  \item[(c)] Montrer que $\dfrac{x}{\ln x}=O(\pi(x))) $ quand $x$ tend vers $+ \infty$.
 \end{itemize} 
\subparagraph*{(a)} 
On a $\binom{2n}{n}=\dfrac{(2n)!}{(n!)^{2}}$ et :
\begin{center}
$v_{p}(\dfrac{(2n)!}{(n!)^{2}})= \sum\limits_{{k=1}}^\infty (\lfloor \dfrac{2n}{p^{k}}\rfloor - 2\lfloor \dfrac{n}{p^{k}}\rfloor) $
\end{center}
Or, $\lfloor 2x\rfloor - 2\lfloor x\rfloor = 0$ ou 1 donc  
\begin{center}
$\sum\limits_{{k=1}}^\infty (\lfloor \dfrac{2n}{p^{k}}\rfloor - 2\lfloor \dfrac{n}{p^{k}}\rfloor)\leq Card\lbrace k\in\mathbb{N}^{*} \mid \lfloor \dfrac{2n}{p^{k}}>0\rfloor\rbrace \leq \lfloor\dfrac{\ln (2n)}{\ln p}\rfloor $
\end{center}
De plus, les nombres premiers diviseurs de $\binom{2n}{n}=\dfrac{(2n)!}{(n!)^{2}}$ sont inférieurs à $2n$ (lemme d'Euclide). Il en découle que  
$\binom{2n}{n}\mid \prod\limits_{{p\in P, p\leq 2n}}{p^{\lfloor\dfrac{\ln (2n)}{\ln p}\rfloor}}$.

\subparagraph*{(b)} 
On a
\begin{center}
 $\binom{2n}{n} \leq \prod\limits_{{p\in P, p< 2n}}{p^{\lfloor\dfrac{\ln (2n)}{\ln p}\rfloor}} \leq \prod\limits_{{p\in P, p< 2n}}{p^{\dfrac{\ln (2n)}{\ln p}}} \leq \prod\limits_{{p\in P, p< 2n}}{(2n)}=(2n)^{\pi(2n)} $
\end{center}

\subparagraph*{(c)}
On passe au logarithme et on utilise une comparaison série intégrale pour déduire que   $\dfrac{2n}{\ln(2n)}=O(\pi(2n))$. Puis on a $\dfrac{x}{\ln(x)}\backsim \dfrac{2\lfloor x/2\rfloor}{\ln(2\lfloor x/2\rfloor)}$ et $\pi(x)\backsim\pi(2\lfloor x/2\rfloor) $, ce qui nous permet de conclure.

\subsection*{3.4. Lifting The Exponent (LTE)}
\addcontentsline{toc}{subsection}{3.4. Lifting The Exponent (LTE)}
\paragraph*{Théorème 1 :}
\subparagraph*{} Si $p$ est un nombre premier impair, premier avec les entiers $a$ et $b$ avec $p \mid a-b$ alors :
\begin{center}
  $v_{p}(a^{n}-b^{n})=v_{p}(a-b)+v_{p}(n)$
  \end{center}  
\paragraph*{Preuve :}
On procède par récurrence sur $v_{p}(n)$.\\
\textbf{Cas $v_{p}(n)=0$ :} On a $v_{p}(a^{n}-b^{n})=v_{p}(a-b)+v_{p}(a^{n-1}+a^{n-2}b+...+b^{n-1})$. Remarquons maintenant que : 
\begin{center}
$ a^{n-1}+a^{n-2}b+...+b^{n-1}\equiv n.a^{n-1}\mod p$ 
\end{center}
Car $a\equiv b\mod p$. On en déduit que  $v_{p}(a^{n}-b^{n})=v_{p}(a-b)$.\\
 \textbf{Cas $v_{p}(n)=1$ :} On pose $n=pn_{1}$ avec $p$ et $n_{1}$ premier entre eux. On a alors : 
\begin{center}
 $v_{p}(a^{pn_{1}}-b^{pn_{1}})=v_{p}((a^{p})^{n_{1}}-(b^{p})^{n_{1}})=v_{p}(a^{p}-b^{p})$
\end{center}
On pose maintenant $a=b+kp$ puisque $p \mid a-b$. On trouve : 
\begin{center}
  $(b+kp)^{p}-b^{p}=\binom{p}{1} (kp)+\binom{p}{2} (kp)^{2}+...+\binom{p}{p}(kp)^{p} $
  \end{center}  
En utilisant le fait que $ p\mid\binom{p}{i} $ pour tout $1\leq i\leq p-1$, on trouve que (Thm 1 et 2 de 3.1.) :  
 \begin{center}
 $ v_{p}((b+kp)^{p}-b^{p})= v_{p}(\binom{p}{1}p)+v_{p}(k)=2+v_{p}(a-b)-1$
 \end{center}
 D'où le résultat.\\
 \textbf{H.R.} On suppose que notre théorème est vrai pour $v_{p}(n)=k$ et on montre qu'il est vrai pour $v_{p}(n)=k+1$.\\
   On pose $n=p^{k+1}n_{1}$. On aura alors :
  \begin{center}
  $v_{p}((a^{p^{k}})^{pn_{1}}-(b^{p^{k}})^{pn_{1}})= v_{p}(a^{p^{k}}-b^{p^{k}})+1 = v_{p}(a-b)+k+1$ 
  \end{center}
Selon le cas $=1$ et H.R. 
 \paragraph*{Corollaire :}
 \subparagraph*{} Si $p$ est un nombre premier impair, premier avec les entiers $a$ et $b$ avec $p \mid a-b$ et $n$ est un entier \textbf{impair} alors :
 \begin{center}
  $v_{p}(a^{n}+b^{n})=v_{p}(a+b)+v_{p}(n)$
  \end{center} 

\paragraph*{Théorème 2 :}
\subparagraph*{} Si $p=2$, $n$ est pair et : 
\begin{itemize}
 \item  $4\mid x-y$ alors $v_{2}(x^{n}-y^{n})=v_{2}(x-y)+v_{2}(n)$
 \item $4\mid x+y$ alors $v_{2}(x^{n}-y^{n})=v_{2}(x+y)+v_{2}(n)$
 \end{itemize} 
 
 \paragraph*{Application 1 :} 
 \subparagraph*{} Soient $a, n$ deux entiers strictement positifs et $p$ un nombre premier impair tel que $a^{p}\equiv 1 \mod p^{n}$. Montrer que $a\equiv 1 \mod p^{n-1}$.
 \paragraph*{Solution.} Il est clair que $a$ et $p$ sont premiers entre eux. D'après le petit théorème de Fermat, $a^{p-1}\equiv 1 \mod p$. Puisque $a^{p}\equiv 1 \mod p$, on en déduit que $a\equiv 1 \mod p$. On peut donc utiliser LTE et on obtient :
\begin{center}
 $v_{p}(a-1)+1= v_{p}(a-1)+v_{p}(p)=v_{p}(a^{p}-1)$
\end{center}
 Le dernier terme étant supérieur ou égal à $n$, il en découle que  $v_{p}(a-1)\geq n-1$.
 
 \paragraph*{Application 2 :}
\subparagraph*{} Soit $k$ un entier strictement positif. Trouver tous les entiers strictement positifs $n$ tels que $3^{k}$ divise $2^{n}-1$.

\paragraph*{Solution.} Soit $k$ tel que $3^{k}$ divise $2^{n}-1$. En raisonnant modulo 3, on voit que $n$ est pair. Écrivons donc $n=2m$ avec $m>0$. Alors $3^{k}$ divise $4^{m}-1$. Comme 3 divise $4-1$, on applique LTE : 
\begin{center}
$v_{3}(4-1)+v_{3}(m)= v_{3}(4^{m}-1)\geq k$
\end{center}
On en déduit que $v_{3}(m)\geq k-1$. Ainsi $2\times 3^{k-1}$ divise $n$.\\
Réciproquement, le même raisonnement donne que $3^{k}$ divise $2^{n}-1$ si $2\times 3^{k-1}$ divise $n$.

\section*{4. Les nombres $p$-adiques au lycée}
\addcontentsline{toc}{section}{4. Les nombres $p$-adiques au lycée}
Certes, la partie précédente étant constituée majoritairement de problèmes et d'applications aux olympiades des mathématiques, elle est techniquement accessible aux élèves de lycée. Mais elle ne parle que des valuations $p$-adiques pour les nombres entiers, on essaiera d'élargir notre champ de vulgarisation dans cette partie, en trouvant une façon d'introduire tout l'ensemble des nombres $p$-adiques $\mathbb{Q}_{p}$, ainsi que celui des entiers $p$-adiques $\mathbb{Z}_{p}$, aux lycéens. On conclut donc notre projet par ce qu'on croit être la façon la plus simple d'aborder les nombres $p$-adiques et donc l'un des deux corps valués locaux non archimédiens existants.  
\subsection*{4.0. Construction des nombres réels}
\addcontentsline{toc}{subsection}{4.0. Construction des nombres réels}
Si on oublie qu'on sait ce qu'un nombre réel est, on peut présenter $\mathbb{R}$ autrement. En effet, fixons un entier naturel $p\geq 2$ et considérons toutes les expressions éventuellement \textbf{infinies à droite} :  
\begin{center}
$...+a_{2}\times p^{2}+a_{1}\times p+a_{0}+a_{-1}\times\dfrac{1}{p}+a_{-2}\times\dfrac{1}{p^{2}} $\textbf{+...}
\end{center}
où les $a_{i}$ sont des entiers naturels strictement inférieurs à $p$. On peut construire tout l'ensemble des réels ainsi.\\
En effet, n'importe quel réel peut s'écrire avec une infinité de décimales
\begin{center}
$\dfrac{1}{3}=0,3333...$, $\sqrt{2}=1,4142...$, $\pi=3,14159...$
\end{center}
Inversement, on peut placer n'importe quel nombre avec une infinité de décimales sur la droite réelle. Tout réel positif s'écrit donc comme somme éventuellement infinie à droite :
\begin{center}
$...+a_{2}\times 100+a_{1}\times 10+a_{0}+a_{-1}\times\dfrac{1}{10}+a_{-2}\times\dfrac{1}{100} $\textbf{+...}
\end{center}
où les $a_{i}$ sont des entiers naturels entre 0 et 9.\\
Ceci reste valable si on remplace $10$ par un entier $p\geq 2$. D'où notre première définition.
\subparagraph*{} On peut faire des opérations sur ces nombres en opérant sur leur partie \textbf{gauche.}\\
Prenons par exemple $p=10$ et calculons $x^{2}$ pour $x=1,4142...$ $(\sqrt{2})$\\
On fait 
\begin{center}
$(1,4)^{2}=1,96$ puis $(1,414)^{2}=1,999396$ etc
\end{center}
pour trouver $1,99999...$ qu'on identifie avec 2. On a donc $x^{2}=2$.\\
On signe ensuite ces nombres pour construire tout $\mathbb{R}$. 
\subsection*{4.1. Les nombres $p$-adiques}
\addcontentsline{toc}{subsection}{4.1. Les nombres $p$-adiques}
On fixe toujours $p\geq 2$ et on considère l'ensemble $\mathbb{Q}_{p}$ de toutes les expressions éventuellement \textbf{infinies à gauche }:
\begin{center}
\textbf{...+ }$a_{2}\times p^{2}+a_{1}\times p+a_{0}+a_{-1}\times\dfrac{1}{p}+a_{-2}\times\dfrac{1}{p^{2}}+... $
\end{center}
où les $a_{i}$ sont des entiers naturels strictement inférieurs à $p$.\\
On définit les opérations en considérant la partie \textbf{droite} de ces nombres. Prenons par exemple $p=10$ et calculons $x^{2}$ pour $x=...99999999$. On fait 
\begin{center}
 $(99)^{2}=9801$ puis $(9999)^{2}=99980001$ etc
 \end{center} 
pour trouver $...0000001$ que l'on identifie avec 1. Donc $x^{2}=1$. \\
On peut aussi calculer $x+1=0$ si bien qu'en fait $x=-1$.\\\\

On obtient ainsi l'ensemble $\mathbb{Q}_{p}$ des \textbf{nombres $p$-adiques}. Ceux qui n'ont pas de virgule sont \textbf{les entiers $p$-adiques} dont l'ensemble est noté $\mathbb{Z}_{p}$.  

\subsection*{4.2. Les entiers $p$-adiques}
\addcontentsline{toc}{subsection}{4.2. Les entiers $p$-adiques}
$\mathbb{Z}_{p}$ est donc l'ensemble des expressions éventuellement infinies à gauche : 
\begin{center}
\textbf{...+ }$a_{3}\times p^{3}+a_{2}\times p^{2}+a_{1}\times p+a_{0} $
\end{center}
où les $a_{i}$ sont des entiers naturels strictement inférieurs à $p$.\\
$\mathbb{Z}_{p}$ est un anneau. En particulier, tous les entiers $p$-adiques ont alors un opposé. On a vu par exemple que $-1=...99999999$ dans $\mathbb{Z}_{10}$.\\
On peut diviser dans $\mathbb{Z}_{p}$, par exemple, dans $\mathbb{Z}_{2}$ : 
\begin{center}
 $...10101011=1+\sum\limits_{{i=0}}^\infty 2^{2i+1} =1+\dfrac{2}{1-4}=\dfrac{1}{3}$
 \end{center} 

\subsection*{4.3. La valuation $p$-adique}
\addcontentsline{toc}{subsection}{4.3. La valuation $p$-adique}
La valuation $p$-adique $v_{p}(x)$ d'un nombre $x$ est le nombre de zéros à droite dans l'écriture de $x$ en base $p$ si $x$ est entier, ou l'opposé du nombre de chiffres après la virgule si $x$ est 'décimal'.\\
Plus formellement, si on écrit $x$ : 
\begin{center}
$x=...+a_{2}\times p^{2}+a_{1}\times p+a_{0}+a_{-1}\times\dfrac{1}{p}+a_{-2}\times\dfrac{1}{p^{2}}+...+a_{r}\times\dfrac{1}{p^{r}} $
\end{center}
avec $a_{r}\neq 0$, on aura alors $v_{p}(x)=r$.\\
Par exemple, on a $v_{10}(2400)=2$ et $v_{10}(1,625)=-3$. On peut aussi calculer $v_{2}(2400)=5$ et $v_{2}(1,625)=-3$ (pourquoi ?). 

\subsection*{4.4. La valeur absolue $p$-adique}
\addcontentsline{toc}{subsection}{4.4. La valeur absolue $p$-adique}
Si $x\in \mathbb{Q}_{p}$ alors la valeur absolue $p$-adique de $x$ est : 
\begin{center}
$\vert x\vert_{p}=\dfrac{1}{p^{v_{p}(x)}}$
\end{center}
Si $x, y\in \mathbb{Q}_{p}$ alors la distance $p$-adique entre $x$ et $y$ est définie par :
\begin{center}
$d_{p}(x,y)=\vert x-y\vert_{p}$
\end{center}
Par exemple, $7\times 11\times 13$ est $10$-adiquement très proche de $1$, car :\\$d_{10}(7\times 11\times 13,1)=\vert 1000\vert_{10}=0,001$.\\
Ou encore, 3100 est 2-adiquement proche de 28. On peut calculer :\\
$d_{2}(3100,28)=\vert 3072\vert_{10}=\dfrac{1}{1024}\simeq 0,001$.\\
Car $3072 = 3\times 2^{10}$ si bien que si bien que $v_{2}(3072) = 10$ et $2^{10} = 1024$.

\subsection*{4.5. Géométrie $p$-adique}
\addcontentsline{toc}{subsection}{4.5. Géométrie $p$-adique}
Comme on dispose d'une distance sur l'ensemble des nombres $p$-adiques, on peut faire de la géométrie $p$-adique : 
\paragraph*{Théorème :}
  \subparagraph*{} La valeur absolue $p$-adique vérifie \textbf{l'inégalité ultramétrique} 
 \begin{center}
  $\vert x+y\vert_{p} \leq \max (\vert x\vert_{p} , \vert y\vert_{p})$
 \end{center}
\paragraph*{Conséquence :} Pour tout $x, y, z$, on a 
\begin{center}
$d_{p}(x,z)\leq  \max (d_{p}(x,y), d_{p}(y,z))$
\end{center}
Cela veut dire que dans un triangle la longueur d'un côté quelconque est toujours inférieure à la longueur d'un des deux autres.
\subparagraph*{} On peut démontrer d'autres résultats déconcertants en géométrie $p$-adique :  
\paragraph*{Théorème 1 :} Tous les triangles sont isocèles.
\paragraph*{Démonstration :} 
Notons $a, b, c$ les longueurs des côtés du triangle. On peut
supposer que $a$ est plus petit que (ou égal à) $b$. L'inégalité ultramétrique nous dit que $c$ est plus petit que $a$ ou plus petit que $b$. De toutes façons, comme $a$ est plus petit que $b$, $c$ sera plus
petit que $b$. Donc, le maximum de $a$ et $c$ est plus petit que $b$.
Mais l'inégalité ultramétrique nous dit qu'il est aussi plus grand
que $b$ : il doit être égal à $b$. Il y a donc bien deux côtés qui ont
même longueur. 
\paragraph*{Théorème 2 :}
N'importe quel point d'un disque est au centre du disque.
\paragraph*{Démonstration :}
Soit $y$ un point du disque $D(x, a)$ de centre $x$ et de rayon $a$. Alors,
la distance de $x$ à $y$ est plus petite que $a$. Si $z$ est un point
quelconque du disque $D(x, a)$, alors la distance de $x$ à $z$ sera aussi
plus petite que $a$. Et l'inégalité ultramétrique implique que la
distance de $y$ à $z$ est plus petite que $a$. Ça veut dire que $z$ est dans
le disque $D(y, a)$ de centre $y$ et de rayon $a$. Symétriquement, si $z$
est un point quelconque du disque $D(y, a)$, alors $z$ est dans le
disque $D(x, a)$. Les deux disques sont donc identiques. Et $y$ est
donc "un" centre du disque $D(x, a)$.

\newpage

\end{document}